\documentclass[12pt]{book}
\usepackage{amsmath,amssymb,amsthm}
\usepackage[all]{xy}

\newcommand{\N}{\mbox{$\mathbb N$}}

\newcommand{\Q}{\mbox{$\mathbb Q$}}
\newcommand{\R}{\mbox{$\mathbb R$}}
\newcommand{\V}{\mbox{$\mathbb V$}}
\newcommand{\Hilb}{\mbox{${\rm Hilb}$}}
\newcommand{\TMseq}{\mbox{$T\{(M_k)_{k\in\N}\}$}}
\newcommand{\ATree}{\mbox{${\rm ATree}$}}
\newcommand{\ATreeM}{\mbox{${\rm ATree}\{(M_k)_{k\in\N}\}$}}
\newcommand{\PTree}{\mbox{${\rm PTree}$}}
\newcommand{\PTreeM}{\mbox{${\rm PTree}\{(M_k)_{k\in\N}\}$}}

\newcommand{\PRTrx}{\mbox{${\rm PRTree}\{M^X\}$}}
\newcommand{\PRTrprimex}{\mbox{${\rm PRTree}'\{M^X\}$}}
\newcommand{\PRTrupnz}{\mbox{${\rm PRTree}^n\{M^Z\}$}}
\newcommand{\PRTree}{\mbox{${\rm PRTree}$}}
\newcommand{\YTree}{\mbox{${\rm YTree}$}}
\newcommand{\emp}{\mbox{$\emptyset$}}
\newcommand{\sh}{\mbox{${\sqcup\!\sqcup}$}}
\newcommand{\vup}{\mbox{$\vee\!\!\!\nearrow\!$}}
\newcommand{\As}{\mbox{${\mathcal As}$}}
\newcommand{\Com}{\mbox{${\mathcal Com}$}}
\newcommand{\Cmg}{\mbox{${\mathcal Cmg}$}}
\newcommand{\Mag}{\mbox{${\mathcal Mag}$}}
\newcommand{\Minf}{\mbox{${\mathcal Mag}_{\omega}$}}
\newcommand{\Cminf}{\mbox{${\mathcal Cmg}_{\omega}$}}
\newcommand{\Dend}{\mbox{${\mathcal Dend}$}}
\newcommand{\Pois}{\mbox{${\mathcal Pois}$}}
\newcommand{\Lie}{\mbox{${\mathcal Lie}$}}
\newcommand{\Brace}{\mbox{${\mathcal Brace}$}}
\newcommand{\Prim}{\mbox{${\rm Prim}$}}

\theoremstyle{definition}
\newtheorem{defn}{Definition}[section]
\newtheorem{remark}[defn]{Remark}
\newtheorem{example}[defn]{Example}
\theoremstyle{plain}
\newtheorem{lemma}[defn]{Lemma}
\newtheorem{proposition}[defn]{Proposition}
\newtheorem{corollary}[defn]{Corollary}
\newtheorem{thm}[defn]{Theorem}
\DeclareMathOperator{\im}{im}
\DeclareMathOperator{\id}{id}

\def\boxtext#1{%
\vbox{%
\hrule
\hbox{\strut \vrule{} #1 \vrule}%
\hrule
}%
}

\def\Y{\setlength{\unitlength}{.4pt}\begin{picture}(60,40)(0,0)
\put(30,0){\line(0,1){10}} \put(30,10){\line(-1,1){30}}
\put(30,10){\line(1,1){30}}
\end{picture}}

\def\arbredown{\setlength{\unitlength}{.7pt}\begin{picture}(60,40)(0,0)
\put(30,0){\line(0,1){10}} \put(30,10){\line(-1,1){30}}
\put(30,10){\line(1,1){30}} \put(15,25){\line(1,1){15}}
\end{picture}}

\def\arbreup{\setlength{\unitlength}{.7pt}\begin{picture}(60,40)(0,0)
\put(30,0){\line(0,1){10}} \put(30,10){\line(-1,1){30}}
\put(30,10){\line(1,1){30}} \put(45,25){\line(-1,1){15}}
\end{picture}}

\def\ladderone{\begin{picture}(15,20)(0,0)
 \put(11.7,2){\line(0,1){10}}\put(10,10){$\bullet$}
\end{picture}}

\def\ladderonopen{\begin{picture}(15,20)(0,0)
 \put(11.7,2){\line(0,1){10}}\put(10,10){$\circ$}
\end{picture}}

\def\laddertwo{\begin{picture}(15,20)(0,0)
\put(11.7,-8){\line(0,1){10}}\put(10,0){$\bullet$}\put(11.7,2){\line(0,1)
{10}} \put(10,10){$\bullet$} \end{picture}}

\def\ladderthree{\begin{picture}(15,30)(0,0)
\put(11.7,-8){\line(0,1){10}}\put(10,0){$\bullet$}\put(11.7,2){\line(0,1)
{10}} \put(10,10){$\bullet$}\put(11.7,12){\line(0,1){10}}
\put(10,20){$\bullet$} \end{picture}}

\def\coratwo{\begin{picture}(30,45)(0,0)
\put(16,2){\line(0,-1){10}}
 \put(15,0){$\bullet$}
\put(16.7,2){\line(1,2){5}} \put(16.7,2){\line(-1,2){5}}
\put(10,10){$\bullet$} \put(20,10){$\bullet$}
\end{picture}}

\def\arbreun{\begin{picture}(50,45)(0,0)
\put(15,0){\line(0,10){10}}\put(15,10){\line(-1,1){20}}
\put(15,10){\line(1,1){20}}\put(10,15){\line(1,1){15}}
\put(5,20){\line(1,1){10}}
\end{picture}}

\def\arbredeux{\begin{picture}(50,45)(0,0)
\put(15,0){\line(0,1){10}} \put(15,10){\line(-1,1){15}}
\put(15,10){\line(1,1){15}} \put(10,15){\line(1,1){10}}
\put(15,20){\line(-1,1){10}}
\end{picture}}

\def\arbretrois{\begin{picture}(50,45)(0,0)
\put(15,0){\line(0,1){10}} \put(15,10){\line(-1,1){20}}
\put(15,10){\line(1,1){20}} \put(28,20){\line(-1,1){10}}
\put(2,20){\line(1,1){10}}
\end{picture}}

\def\arbrequatre{\begin{picture}(50,45)(0,0)
\put(15,0){\line(0,1){10}} \put(15,10){\line(-1,1){15}}
\put(15,10){\line(1,1){15}} \put(20,15){\line(-1,1){10}}
\put(15,20){\line(1,1){10}}
\end{picture}}

\def\arbrecinq{\begin{picture}(50,45)(0,0)
\put(15,0){\line(0,1){10}} \put(15,10){\line(-1,1){20}}
\put(15,10){\line(1,1){20}} \put(20,15){\line(-1,1){15}}
\put(25,20){\line(-1,1){10}}
\end{picture}}

\begin{document}

\begin{titlepage}
\begin{center}
\
\bigskip
\bigskip
\bigskip
\bigskip
\bigskip
\bigskip

 {\LARGE On Hopf algebra structures over operads}

\bigskip
\bigskip


\bigskip
\bigskip
\bigskip
\bigskip
\bigskip
\bigskip
\bigskip
\bigskip
\bigskip
\bigskip
\bigskip


\bigskip
\bigskip

{\Large Ralf Holtkamp}

\bigskip
\bigskip
\bigskip
\bigskip
\bigskip
\bigskip
\bigskip
\bigskip
\bigskip
\bigskip
\bigskip
\bigskip
\bigskip
\bigskip

{\large Fakult\"at f\"ur Mathematik}

\bigskip
\bigskip

{\large Ruhr-Universit\"at Bochum}

\bigskip
\bigskip

{\large 2004}

\end{center}
\end{titlepage}

\ \bigskip
\vfill
\centerline{11/2005}

\tableofcontents

\chapter*{Introduction}

\addcontentsline{toc}{chapter}{Introduction}

The observation that many objects possess naturally both a
multiplication and a comultiplication, which are compatible, lead
to the theory of Hopf algebras (cf.\  \cite{liabe},
\cite{limohopf}).

 Hopf\ \cite{lihopf}
discovered that a multiplication $M\times M\to M$ and the diagonal
$M\to M\times M$ on a manifold $M$ induce a multiplication
$H\otimes H\to H$ and a comultiplication $H\to H\otimes H$ on the
homology $H_{*}(M;K)$, and he proved a structure theorem for
$H_{*}(M;K)$.

 At first, commutative or cocommutative Hopf
algebras were studied. Over a field of characteristic 0, any
connected cocommutative Hopf algebra $H$ is of the form $U(\Prim
H)$, where $\Prim H$ is the space of primitive elements viewed as
a Lie algebra, and where $U$ is the universal enveloping functor.
This is the theorem of Milnor and Moore (proven for graded
Hopf algebras in \cite{limimo}).
Commutative Hopf algebras are $K$-linear duals of cocommutative
Hopf algebras.

Only a few types of examples of noncommutative non-cocommutative
Hopf algebras were constructed (cf.\ \cite{litaf}), when
 the study of quantum groups
(cf.\ \cite{likas}), emerged in the last decades of the twentieth
century. Faddeev, Reshetikhin, Takhtajan (cf.\ \cite{litak}) and
others obtained quantum groups as noncommutative deformations of
algebraic groups.
 The notion of a quantum group was
introduced by Drinfel'd\ \cite{lidri}. Quantum groups are given by
noncommutative Hopf algebras, the category of quantum groups being
opposite to the category of noncommutative Hopf algebras (see
\cite{limanin}).

\bigskip

Several combinatorial Hopf algebras (see\ \cite{liabs}) were
discovered, which provided new links between algebra, geometry,
combinatorics, and theoretical physics.

For example, noncommutative Hopf algebras of permutations and of
quasi-sym\-metric functions (cf.\ \cite{limare}) allow a new
approach to the representation theory of symmetric groups.

Connes and Kreimer discovered (cf.\ \cite{lick}, \cite{likr}) that
a commutative non-cocommutative Hopf algebra structure on rooted
trees encodes renormalization in quantum field theory. The same
formulas define a noncommutative Hopf algebra on planar rooted
trees.
 The use of binary planar rooted trees has been
proposed by Brouder and Frabetti \cite{libfa} for renormalization.

\bigskip

Recent developments in Hopf algebra theory show that there are
important objects which could be called "non-classical" Hopf
algebras.
 Typical examples are
dendriform Hopf algebras \cite{liroa,lirob}. Dendriform  algebras
are equipped with two operations $\prec$, $\succ$ whose sum is an
associative multiplication.

Other examples are magma Hopf algebras (considered in
\cite{lighmag}), or infinitesimal Hopf algebras (cf.\
\cite{liasinf}). Here the usual Hopf algebra axioms have to be
changed.

\bigskip

Instead of the classical types of algebras (like associative
algebras, commutative algebras or Lie algebras) one can consider
new types of algebras as well. The idea to model these types by
their multilinear operations, and to compute with spaces of
operations like monoids with a composition-multiplication, goes
back to Lazard's analyzers \cite{lilaloi} and is now known as the
theory of operads. The concept of Koszul
duality, being well-known for quadratic algebras,
 exists also for operads (see \cite{ligk}). As for the types
$\Lie$ and $\Com$, the theory of algebras over one operad is
strongly connected with the theory of algebras over the dual
operad.

The first applications of operads, introduced by Boardman and Vogt
 and named by May (\cite{libv}, \cite{limay}), were in
algebraic topology. In the 1990's the theory of
operads developed rapidly. Operads are important tools for
deformation quantization and formality, and further influences came from
mathematical physics, where operads of
 moduli spaces of Riemann surfaces with marked
points are relevant for string field theories (cf.\cite{likm},
\cite{likon}).

 Strong homotopy algebras and other types of algebras
 with infinitely many different operations (like $A_{\infty},
B_{\infty}, C_{\infty}$) appear in the context of operads (cf.\
\cite{ligejo}).

\bigskip

Operations modeled by an operad have multiple inputs but only
one output. Thus they can be represented by rooted trees, where
leaves correspond to inputs and the root corresponds to the
output. If there are no relations between the operations (the case
of free operads), we can compute in vector spaces of trees, and a
natural grafting of trees corresponds to the operad-composition
(cf.\ \cite{limss}).

\bigskip

The two operations $\prec, \succ$ of the dendriform operad $\Dend$
define operations on the set of planar binary trees which provide
the space of all planar binary trees with the structure of a free
dendriform algebra (see \cite{lilodia}). Loday and Ronco
\cite{lilrhopftree} introduced a natural
 Hopf algebra structure of the free dendriform algebra on one
 generator.
The Loday-Ronco Hopf algebra is in fact isomorphic to the
noncommutative planar Connes-Kreimer Hopf algebra and the
Brouder-Frabetti Hopf algebra (\cite{licompa}, \cite{lifo}). Thus
all three have the structure of a dendriform Hopf algebra.

\bigskip
\bigskip
\bigskip

In this work, we consider examples of "non-classical" Hopf
algebras. Just as $\mathcal P$-algebras, i.e.\ algebras over
different operads $\mathcal P$, make sense, $\mathcal P$-Hopf
algebras are defined (see Chapter 3). Since we deal with one
(coassociative) cooperation, there are coherence conditions which
the operad $\mathcal P$ has to fulfill, such that the tensor
product of $\mathcal P$-algebras is again a $\mathcal P$-algebra.

We do not include objects with various not necessarily associative
operations and not necessarily coassociative cooperations, but the
given definition is general enough to include dendriform Hopf
algebras (compare \cite{lilosci}).

To avoid problems with antipodes, we restrict to the case
 of filtered or graded Hopf algebras.

\bigskip

 Apart from dendriform Hopf algebras, we describe
  Hopf algebra structures over free operads.
Therefore combinatorial operations and admissible labelings for
various sorts of rooted trees are needed. We develop these tools
in Chapter 2. As an application, we describe the Loday-Ronco
dendriform Hopf algebra. We also describe its isomorphisms with
the noncommutative planar Connes-Kreimer Hopf algebra and with a
Hopf algebra of Brouder and Frabetti (in Chapter 3).

\bigskip

 Let $\Mag$ be the operad freely generated by a
non-commutative non-associative binary operation $\vee^2(x_1,x_2)$
also denoted by $x_1\cdot x_2$. More exactly, the free operad
$\Mag$ is the symmetrization of the free non-$\Sigma$ operad
generated by $\vee^2$.

Let similarly $\Minf$ be freely generated by operations $\vee^n$,
one for each $2\leq n\in \N$. A basis for the space of $n$-ary
operations is given by reduced planar
rooted trees with $n$ leaves. (This is the operad of Stasheff
polytopes, see \cite{listaphys}).

The suboperad $\Mag$ contains only the binary trees. Thus
free algebras over $\Mag$ have planar binary rooted trees (with labeled
leaves) as a
vector space basis. These free algebras are called magma algebras,
they are free non-associative.

\bigskip

For any operad $\mathcal P$ that fulfills the coherence conditions
noted above, the free $\mathcal P$-algebra generated by a set $X$
of variables is always a $\mathcal P$-Hopf algebra with the
diagonal or co-addition $\Delta_a$ as a comultiplication. The
$\mathcal P$-algebra homomorphism $\Delta_a$ is determined by
$\Delta_a(x)=x\otimes 1+1\otimes x,$ for all $x\in X$.

It is natural to pose the question of primitive elements.
 In the classical case, where
the operad is the operad $\As$ of associative algebras, Lie
polynomials occur as the primitive elements. The theorem of
Friedrichs (cf.\ \cite{lireu}) characterizes Lie polynomials as
the primitive elements of the free associative algebra equipped
with the coaddition.

In a natural way, there exists an operad $\Prim{\mathcal P}$, and
$\Prim\As=\Lie$. In the commutative case, only homogeneous
elements of degree 1 are primitive, i.e.\ $\Prim\Com={\rm Vect}$.
A table of pairs ${\mathcal P}, \Prim{\mathcal P}$ (triples in fact,
if coassociativity of the comultiplication is replaced by something
different) is given by Loday and Ronco in (\cite{lilr04}),
and we are going to add some new cases to the list.

\medskip

We look at free $\Mag$- and $\Minf$-algebras equipped with the
co-addition, i.e.\ we are interested in the operads
$\Prim\Mag$ and $\Prim\Minf$.

Although these algebras are free non-associative,
results of Drensky-Gerritzen \cite{lidg} and Gerritzen
\cite{ligeplaexp} on a canonical exponential function in
non-associative variables show that they have a rich structure.

In such free $\mathcal P$-algebras, primitive elements for
the co-addition are also constants for the partial derivatives
with respect to the variables. Whereas the primitive elements do
not form a $\mathcal P$-algebra, the constants do form a $\mathcal
P$-algebra. For example, the algebra of constants of the free
associative algebra $K\langle X\rangle$ is generated by Lie
polynomials (\cite{lifalk}).

We consider
Taylor expansion for polynomials in non-associative variables,
generalizing classical Taylor expansion or Taylor expansion in
associative algebras (cf.\ \cite{lidrecodim, lidre2}, and
\cite{ligea, ligeb}). These Taylor expansions yield projectors on constants.

We examine the first generators and relations needed to describe
$\Prim\Mag$ and $\Prim\Minf$. The first non-trivial relation is
 a non-associative Jacobi identity.

Reviewing a result of \cite{lighmag}, also obtained by Shestakov
and Umirbaev \cite{lisu}, we show that it does not suffice to use
$n$-ary operations for $n\leq 3$ to generate the operads
$\Prim\Mag$ and $\Prim\Minf$.

The same is true for the (operad of) primitive elements of dendriform
algebras. The results of Ronco (see \cite{liroa}, \cite{lirob},
\cite{lilosci}) show that $\Prim\Dend$ is the operad $\Brace$ of
brace algebras.

 In the case of dendriform Hopf algebras, the
question of primitive elements is the question of primitive
elements in the Loday-Ronco Hopf algebra, or equivalently in the
noncommutative Connes-Kreimer Hopf algebra.

\bigskip

In order to describe the operads $\Prim\Minf$ and $\Prim\Mag$, we
describe the graded duals of the given co-addition $\Mag$-Hopf
algebras. These duals are equipped with a commutative
multiplication $\sh$, which generalizes the shuffle multiplication
known for free associative algebras. We show in Chapter 4, that
these commutative algebras are freely generated by the primitive
elements (and dually the given co-addition Hopf algebras are
connected co-free). This is an analogon of the
Poincar\'e-Birkhoff-Witt theorem.

We determine the generating series for the operads $\Prim\Mag$ and
$\Prim\Minf$. We show that the dimension of ${\Prim\Mag}(n)$ is
related to the $\log$-Catalan numbers $c'_1=1,\ c'_2=1,\ c'_3=4,\
c'_4=13,\ c'_5=46,\ c'_6=166,\ \ldots\ \ $ by
\begin{equation*}
\dim{\Prim\Mag}(n)=(n-1)! c'_n.
\end{equation*}

By a recursive method, we show how the spaces $\Prim_n$ of
homogeneous elements of degree $n$ can be described as
$\Sigma_n$-modules. As in the case of $\Lie_n$, this method only
works for small $n$. For dendriform algebras, Ronco proved a
Milnor-Moore type theorem (see \cite{liroa, lirob}) using Eulerian
idempotents as projections onto $\Prim\Dend$. Similar projectors
may be defined for $\Mag$ and $\Minf$.

\bigskip

In the last section of Chapter 4, we sketch a generalized
Lazard-Lie theory (see \cite{lilaloi}, \cite{lipseudo}) for
$\mathcal P$-Hopf algebras.

Tuples of primitive elements occur as elements
of the first cohomology group $H^1$.
Elements of $H^2$ and $H^3$ are obstructions (for
uniqueness and existence) that are relevant for a
 desired classification of (complete) $\mathcal P$-Hopf algebras.

\bigskip
\goodbreak

In more detail, the content of this work can be outlined as
follows.

In Sections 1.1-1.3 we recall the definitions, together with first
examples, for $\Sigma$-vector spaces, operads and non-$\Sigma$
operads, $\circ_i$-operations, generating series, algebras (and
coalgebras) over operads, and especially free (free complete)
algebras and cofree coalgebras.

\bigskip

In Section 2.1, we introduce the notions of admissibly labeled
(planar or abstract) rooted trees, where a sequence of sets $M_k$
contains the allowed labels for vertices of arity $k$, see
Definition (\ref{defadmiss}). This notion is useful to present
operads (especially free operads), several examples of dendriform
Hopf algebras and $\Mag$- or $\Minf$-Hopf algebras in a consistent
way.

Among other integer sequences, we discuss the sequences of Catalan
numbers and $\log$-Catalan numbers, see Example
(\ref{exlogcatalan}), which will be relevant in Chapter 4.

We also recall some important operations on trees, like grafting
and de-grafting operators.

\bigskip

In Section 2.2, the connection between trees and parenthesized
words is explained. We describe cuts of trees, and we introduce
leaf-restrictions, see Lemma (\ref{lemleafrestrict}), and
leaf-splits. Shuffles for trees are also defined.

\bigskip

In Section 2.3, we first focus on planar binary trees, and a
(right) comb representation for them, which is used for a
bijection between such trees and forests of not necessarily binary
planar trees, see Lemma (\ref{lemtreebij}). Then we consider
Stasheff polytopes and super-Catalan numbers.

\bigskip

Section 2.4 is a continuation both of Sections 2.1-3 and Sections
1.1-3. We present free (non $\Sigma$-)operads using admissibly
labeled trees, see Lemma (\ref{lemfreeoperad}). Especially we
consider the operads $\Mag$ and $\Minf$ which are the main topic
of Chapter 4.

\bigskip

After a brief introduction to classical Hopf algebra theory, given
in Section 3.1, we consider in Sections 3.2 and 3.3 $\mathcal
P$-Hopf algebras for operads $\mathcal P$ with a coherent unit
action.

Coherent means that the tensor product of unitary $\mathcal
P$-algebras $A, B$ is equipped with the structure of a $\mathcal
P$-algebra in a way such that elements $a\otimes 1$ (for $a\in A$)
generate a $\mathcal P$-algebra isomorphic to $A$, and elements
$1\otimes b$ (for $b\in B$) generate a $\mathcal P$-algebra
isomorphic to $B$. Operads like $\Com$, $\As$, $\Mag$, $\Cmg$,
$\Pois$ and $\Dend$ are equipped with a coherent unit action.

 In the definition of a $\mathcal P$-Hopf algebra, see Definitions (\ref{defpbialg}) and (\ref{defphopf}),
 the comultiplication $\Delta$ is required to be a coassociative map
which is also a morphism of unitary $\mathcal P$-algebras.

Given a graded  $\mathcal P$-Hopf algebra $A$, we also consider
the graded dual $A^{*g}$, which is a graded $\As$-algebra with
respect to the operation $\Delta^*$, see Lemma
(\ref{lemdualhopf}).

\bigskip

In Sections 3.4-3.6 we focus on $\Dend$-Hopf algebras. The Hopf
algebra structures on rooted trees introduced by Loday and Ronco,
Brouder and Frabetti, and the (non-commutative) Connes-Kreimer
Hopf algebra are described. We review the isomorphisms between
them in an elementary way, see Propositions (\ref{propcompa}) and
(\ref{propcompb}).

\bigskip

In Section 3.7 we show that there exists an operad $\Prim\mathcal
P$, see Lemma (\ref{lemprimoperad}), with $\Prim\mathcal P(n)$
corresponding to the homogeneous multilinear primitive elements of
degree $n$ in the free $\mathcal P$-algebra (on $n$ variables). We
give the classical examples, and we recall the  Milnor-Moore
theorem.

 We sketch results of Ronco and Loday-Ronco, which show that the operads
 of brace algebras
 and (non-dg) $B_{\infty}$-algebras are of the form $\Prim\mathcal
P$.

\bigskip

In Section 4.1 we collect some information about the
 operads $\Mag$ and
$\Minf$. Both are equipped with a (canonical) coherent unit
action. We view the elements of the free $\Mag$-algebras $K\{X\}$
and $\Minf$-algebras $K\{X\}_{\infty}$ as polynomials, while
monomials correspond to admissibly labeled planar trees.

\bigskip

We consider formal derivatives for $\Mag$- and $\Minf$-algebras in
Section 4.2. We give a concrete description how to compute the
partial derivatives $\partial_k$ of tree monomials, using the
concept of leaf-restriction (introduced in Section 2.2), see
Proposition (\ref{proppartial}).

\bigskip

We show in Section 4.3 that the free unitary $\Minf$-algebra
$K\{X\}_{\infty}$ together with the co-addition map $\Delta_a$ is
a strictly graded $\Minf$-Hopf algebra, containing the free
unitary $\Mag$-algebra $K\{X\}$ as a sub-Hopf algebra, see
Proposition (\ref{propminfhopf}). We describe how the co-addition
acts on monomials, see Proposition (\ref{propdelta}).

We introduce generalized differential operators $\partial_T$ and
prove properties of these operators,
 see Proposition (\ref{propdiffop}).

\bigskip

In Section 4.4, see Proposition (\ref{proptaylor}), we introduce
the concept of Taylor expansions for $\Minf$-algebras
$K\{X\}_{\infty}/I$, generalizing classical Taylor expansions.

\bigskip

While in Sections 4.1 to 4.4 we prove properties of (free)
$\Minf$-algebras, which may immediately be translated into
properties of (free) $\Mag$-algebras, we start in Section 4.5 to
give a description of constant and primitive elements of
$K\{X\}_{\infty}$ and $K\{X\}$ separately.

The concept of Taylor expansion is used to describe the spaces of
primitive elements $\Prim\Mag(n)$ for degrees $n\leq 3$, see
Proposition (\ref{propmagdegthree}). In Proposition
(\ref{propminfdegthree}) the spaces of primitive elements
$\Prim\Minf(n)$ for degrees $n\leq 3$ are described.

We also prove that $\Prim\Mag$ and $\Prim\Minf$ can not be
generated by quadratic and ternary operations, see Corollary
(\ref{cornotquadtern}).

In Section 4.6, see Propositions (\ref{propdualshuffle}) and
(\ref{propdualshmag}), we explicitly describe the graded duals of
the $\Minf$-Hopf algebra $(K\{X\}_{\infty},\Delta_a)$ and of the
$\Mag$-Hopf algebra $(K\{X\},\Delta_a)$. Therefore we consider a
 shuffle multiplication which is a sum of the shuffles (introduced
 in Section 2.2).

\bigskip

For the operads $\Mag$ and $\Minf$, we prove an analogon of the
Poincar\'e-Birkhoff-Witt theorem, see Theorem (\ref{thmpbw}) and
Theorem (\ref{thmpbwdual}), in Section 4.7. Since the proof is the
same for the $\Mag$- and the $\Minf$-case, we treat both cases
together, see Proposition (\ref{propshuffortho}) and Lemma
(\ref{lemcomag}).

\bigskip

In Section 4.8 we have a closer look at the operad $\Prim\Mag$.
Using the explicit generating function of the operad $\Mag$, we
obtain from Theorem (\ref{thmpbw}) the Corollary
(\ref{corlogcatalan}), which states that $\dim\Prim\Mag(n)$ is
given by (a multiple of) the $n$-th $\log$-Catalan number, for all
$n$.

We describe $\Prim\Mag(4)$ by means of representation theory, see
Proposition (\ref{propmagdegfour}), and we demonstrate how to
compute a basis of highest weight vectors. We finish with the
analogues of Corollary (\ref{corlogcatalan}) and Proposition
(\ref{propmagdegfour}) for the case of $\Prim\Minf$.

\bigskip

The aim of the last section, Section 4.9, is to show that the
question for primitive elements is only the first step in a
generalized Lazard-Lie theory for $\mathcal P$-Hopf algebras. The
goal is to obtain classification results for (filtered, or graded,
or complete) $\mathcal P$-Hopf algebras.

\bigskip
\bigskip
\bigskip

 I wish to thank my teacher L.\ Gerritzen for his support, many
discussions, and his constant interest and encouragement.

This research was in part carried out at the Institut
Mittag-Leffler in Stockholm, the Institute of Mathematics and
Informatics in Sofia, and the Erwin Schr\"odinger International
Institute for Mathematical Physics in Vienna. I thank these
institutes for their hospitality. I am also grateful to the
Mathematical Sciences Research Institute in Berkeley for the
opportunity to spend the research year 1999/2000 there. A lot of
people I met at these institutes had a great influence on this
work. I would especially like to thank J.-L.\ Loday, M.\ Ronco,
V.\ Drensky, D.\ Kreimer, G.M.\ Bergman,
 S.\ Montgomery, and E.\ Taft.

\begin{chapter}{Basic Concepts. }

\begin{section}{Words and permutations}\label{sectconventions}

Let $K$ be a field of characteristic 0, and let Vect$_K$ be the
category of vector spaces over $K$. By $X$ we usually denote a
non-empty countable set of variables. If $X$ is infinite, we
usually consider $X=\{x_1,x_2,x_3,\ldots\}$. If the cardinality
$\# X$ of $X$ is finite, $\# X=m$, we use the variables
$x_1,\ldots,x_m$. We also consider $X$ as a basis of the vector
space $V_X=K X={\rm Span}_K(X)$.

Let $\N=\{0,1,2,\ldots\}$, $\N^*=\{1,2,\ldots\}$.

\begin{defn}

For $m\in \N^*$, we define the set $\underline m:=\{1,\ldots,m\}$.
We also define $\underline 0 := \emptyset$.
 For $m \in \N,$ let $\Sigma_m$ be the symmetric
group, acting from the left on the set $\underline m$ via
$\sigma.i=\sigma(i)$.
As a Coxeter group, we generate $\Sigma_m$ by $\{\tau_1,\tau_2,\ldots,
\tau_{m-1}\}$,
where $\tau_i$ is the transposition exchanging $i$ and $i+1$.

\end{defn}

\begin{remark}

To fix notation, we recall the free objects of the categories of
abelian semigroups, (not necessarily abelian) semigroups, and
magmas. We also consider the associated categories of unitary
semigroups. Words occur as elements of these free objects.

Elements of the free abelian semi-group $W_{\rm Com}(X)$ over $X$
are commutative words $x_{i_1}^{\nu_1}x_{i_2}^{\nu_2}\cdots
x_{i_r}^{\nu_r}$, $i_1 < i_2 < \ldots < i_r$, $r\geq 1$, $\nu_i\in
\N$. Adjoining a unit $1$ (empty word) we get the
 free abelian
semi-group $W^1_{\rm Com}(X)$ with unit.

The free semigroup over $X$ is denoted by $W_{\rm As}(X)$, the
free semigroup with unit by $W^1_{\rm As}(X)$. We denote the
concatenation $W^1_{\rm As}(X)\times W^1_{\rm As}(X)\to W^1_{\rm
As}(X), (v,w)\mapsto v.w$ by a lower dot. The elements of
$W^1_{\rm As}(X)$ are words $w=w_1.w_2\ldots w_r$, $w_i\in X$ for
all $i$. Here $r\in \N$ is the length of $w$.
Elements $\sigma$ of $\Sigma_m$ are also written in one-line notation, i.e.\ as
words $\sigma(1).\sigma(2)\ldots\sigma(m)$.
(The length of a permutation is
the number of inversions $i<j$ with $\sigma(i)>\sigma(j)$, though.)

A magma is just a set $M$ equipped with a binary operation
(usually denoted by $\cdot:M\times M\to M$).
 The elements of the free magma $W_{\rm
Mag}(X)$ over $X$ are parenthesized words. We are going to
identify these words with binary rooted trees, see Section
\ref{sectstring}.

\bigskip

 For any $r$, we have a left action of $\Sigma_r$ on the set $W_{r,\rm As}(X)$ of
words $w=w_1.w_2\ldots w_r$ of length $r$ by place permutation:
for $i=1,\ldots,r$, we move the $i$-th letter to position
$\sigma(i)$, thus $\sigma.(w_1.w_2\ldots w_r)= w_{\sigma^{-1}(1)}.
w_{\sigma^{-1}(2)}\ldots w_{\sigma^{-1}(r)}$.

\end{remark}

\begin{defn}

Let $\sigma\in\Sigma_n$. Let $(m_1,\ldots,m_n)$ be an ordered
partition of $m$, i.e.\ any tuple $(m_1,\ldots,m_n)$ with $m_i\in
\N^*$ (all $i$) and  $m=m_1+\ldots+m_n$. We denote by
$(m'_1,\ldots,m'_n)$ the ordered partition
$(m_{\sigma^{-1}(1)},\ldots,m_{\sigma^{-1}(n)})$ of $m$.

The block permutation $\sigma_{(m_1,\ldots,m_n)}\in\Sigma_m$ acts
on the set $\underline m$
 in the following way:
For each $1\leq i\leq n$ it
 maps the interval $\{j: m_1+\ldots+m_{i-1}<
j \leq m_1+\ldots+m_i\}$ strictly monotonic onto the interval
$\{j': m'_1+\ldots+m'_{\sigma(i)-1}< j' \leq m'_1+\ldots
+m'_{\sigma(i)}\}$. In other words, $\sigma_{(m_1,\ldots,m_n)}$
acts on the intervals of length $m_1,\ldots,m_n$ in the same way
as $\sigma$ acts on $\underline n$. If $n=3, m=9$,
$(m_1,m_2,m_3)=(2,4,3)$ and $\sigma$ is the transposition exchanging 3 and 1,
then $\sigma_{(m_1,m_2,m_3)}$ maps $\{7,8,9\}$, i.e.\ the third
interval with respect to $(m_1,m_2,m_3)$, onto $\{1,2,3\}$. The
interval $\{ 3,4,5,6 \}$ becomes the interval $\{4,5,6,7\}$, and
$\{1,2\}$ is mapped onto $\{8,9\}$.

\bigskip

For every $n\in\N, m_1,\ldots ,m_n\in \N$, $m:=m_1+\ldots+m_n,$
there are maps
\begin{equation*}
\Sigma_n\times \Sigma_{m_1}\times\ldots\times\Sigma_{m_n}\to
\Sigma_m
\end{equation*}
 defined by
\begin{equation*}
(\sigma,\gamma_1,\ldots,\gamma_n)\mapsto\sigma_{(m_1,\ldots
,m_n)}\circ(\gamma_1\times\ldots\times\gamma_n)
\end{equation*}
(with $\gamma_1\times\ldots\times\gamma_n$ acting on $\{1,\ldots
,m_1\}\times\ldots\times\{m-m_n+1,\ldots,m\}\subset \underline
m$.)

 By $K$-linear extension, we get maps
\begin{equation*}
\mu_{n;m_1,\ldots,m_n}: K\Sigma_n\otimes
K\Sigma_{m_1}\otimes\ldots\otimes K\Sigma_{m_n}\to K\Sigma_m.
\end{equation*}
\end{defn}

\begin{remark}

The maps $\mu_{n;m_1,\ldots,m_n}$ turn the sequence
$(K\Sigma_n)_{n\in \N}$  into a $K$-linear operad, see the next
section.

\end{remark}

\begin{remark}

The group ring $K\Sigma_n$
 corresponds to the regular
representation of $\Sigma_n$. We refer to \cite{lijk} for the
representations of the symmetric groups.

The irreducible representations of $\Sigma_n$  are symbolized by
Young diagrams. For example the regular representation of
$\Sigma_2$ is the sum of the trivial representation ${\rm
 \vbox{\offinterlineskip \hbox{\boxtext{\phantom{X
}}\boxtext{\phantom{X }}}}}\ $ and the sign representation ${\rm
 \vbox{\offinterlineskip \hbox{\boxtext{\phantom{X
}}} \hbox{\boxtext{\phantom{X }}}
 } \rm}\ .$

In degree 3, the regular representation is given by
\begin{equation*}
 {\rm \vbox{\offinterlineskip
\hbox{\boxtext{\phantom{X }}\boxtext{\phantom{X
}}\boxtext{\phantom{X }}}
 } \rm}\ \oplus\ 2
 \ {\rm
 \vbox{\offinterlineskip \hbox{\boxtext{\phantom{X
}}\boxtext{\phantom{X }}} \hbox{\boxtext{\phantom{X }}}
 } \rm}
\oplus\  {\rm
 \vbox{\offinterlineskip \hbox{\boxtext{\phantom{X
}}} \hbox{\boxtext{\phantom{X }}} \hbox{\boxtext{\phantom{X }}} }
\rm}\ .
\end{equation*}

\medskip

 The symmetric group $\Sigma_n$ is embedded into $GL_m(K)$, $m\geq
 n$, and there is a similar theory for $GL$-modules, cf.\
\cite{liw}. For any group $G$, a left (right) $G$-module is a vector
space $V$ over $K$ together with a
 left (right) action of $G$.

\end{remark}

\bigskip

\end{section}

\begin{section}{Operads}

 The concept of operads and its
precursors "analyseurs" and "compositeurs" (introduced by Lazard, see
\cite{lilaloi}) is a good device to handle several types of
algebras at once. From this point of view, the ability to form
compositions is most important. On the one hand, looking at the
corresponding free algebras, it is possible to insert elements
into another. On the other hand, looking at the operations which
characterize the algebra structure, we might combine operations to
get new ones (possibly of a higher number of arguments).

An elegant way to introduce operads is to define an operad as a
triple or monad (\cite{liml}, chap.\ VI, cf.\ also \cite{lifra})
in the monoidal category $({\rm End}(\mbox{Vect}_K),\circ)$, cf.\
\cite{likm}. One usually requires to have an action of the
symmetric groups. The analogously defined object without the
additional structure related to $\Sigma_n$-actions is called a
non-$\Sigma$-operad.

\begin{defn}\label{sigmaspace}
Let $\Sigma\mbox{-Vect}_K$ be the following category: Objects
$\mathcal V$ are sequences $({\mathcal V}(n))_{n\in \N}$ of vector
spaces ${\mathcal V}(n)$ with a right $\Sigma_n$-action.

 Morphisms are given by
homomorphisms compatible with the $\Sigma_n$-action.

\end{defn}

To every $\Sigma$-vector space ${\mathcal V}$ there is associated
an endofunctor on Vect$_K$ given by
\begin{equation*}
\begin{split}
&F_{\mathcal V}(V):= \bigoplus_{n=0}^{\infty} \bigl({\mathcal
V}(n)\otimes V^{\otimes
n}\bigr)_{\Sigma_n}=\bigoplus_{n=0}^{\infty} {\mathcal
V}(n)\otimes_{\Sigma_n}V^{\otimes n}.
\\
\end{split}
\end{equation*}
Here $\Sigma_n$ acts from the left on $V^{\otimes n}$ by
 place
permutation, and  $\bigl({\mathcal V}(n)\otimes V^{\otimes
n}\bigr)_{\Sigma_n}$ denotes the space of coinvariants for the
diagonal action of $\Sigma_n$. (Thus $p.\sigma\otimes
q_1\otimes\ldots \otimes q_n$ and $p\otimes
\sigma.(q_1\otimes\ldots \otimes q_{n})$ are identified.)

\bigskip

 The category $\Sigma\mbox{-Vect}_K$ can be
identified with the full subcategory of the category
End$(\mbox{Vect}_K)$ of functors $\mbox{Vect}_K\to \mbox{Vect}_K$
consisting of functors of the form $F_{\mathcal V}$. The
composition $\circ$ of functors induces an associative bifunctor
$\odot$ on End$(\mbox{Vect}_K)$ such that $({\rm
End}(\mbox{Vect}_K),\odot)$ and $(\Sigma\mbox{-Vect}_K, \odot)$
are monoidal categories in the sense of \cite{liml}, chap.\ VII.

Unit object is the identity functor Id$_{\mbox{Vect}_K}$, which is
given by the $\Sigma$-vector space $\mathcal I$ with ${\mathcal
I}(1) =K, {\mathcal I}(n)=0$ for $n\neq 1$.

\begin{defn}\label{defoperad}

 Let a $\Sigma$-space ${\mathcal P}$,
together with morphisms of functors
\begin{equation*}
\mu: F_{\mathcal P}\odot F_{\mathcal P} \to F_{\mathcal P} \text{
 and  } 1:\mbox{Id}\to F_{\mathcal P}
\end{equation*}
 be given such that the following diagrams are commutative:
\begin{equation*}
\xymatrix{
          F_{\mathcal P}\odot F_{\mathcal P}\odot F_{\mathcal P}\
          \
          \ar[d]_{\mu\odot \id}
          \ar[r]^{\ \ \id\odot \mu}
          & F_{\mathcal P}\odot F_{\mathcal P}
           \ar[d]^{\mu}&  & {{\rm Id}\odot F_{\mathcal P}}\ar[dr]_= \ar[r]^{1\odot \id}
            & F_{\mathcal P}\odot F_{\mathcal P}\ar[d]^{\mu}
           & \ar[l]_{\id\odot 1}F_{\mathcal P}\odot {\rm Id}\ar[dl]^=
           \\
             F_{\mathcal P}\odot F_{\mathcal P} \ar[r]_{\mu} & F_{\mathcal P}&
            &    & F_{\mathcal P}
           \\
            }
\end{equation*}
Then $({\mathcal P}, \mu, 1)$, or just ${\mathcal P}$, is called a
($K$-linear) operad. The morphism $\mu$ is called
(operad-)composition.

\end{defn}

\begin{remark}
The elements of ${\mathcal P}(n)$ are often viewed as abstract
operations with $n$ inputs and one output.

From the description of ${\mathcal P}\odot{\mathcal P}$ given by
$F_{{\mathcal P}\odot{\mathcal P}}=F_{\mathcal P}\odot F_{\mathcal
P}$ it follows that the composition $\mu=\Bigl(\mu_V: F_{\mathcal
P}\bigl(F_{\mathcal P}(V)\bigr)\to F_{\mathcal P}(V)\Bigr)$ is
explicitly given by maps
\begin{equation*}
\mu_{n;m_1,\ldots,m_n}: {\mathcal P}(n)\otimes {\mathcal
P}(m_1)\otimes \ldots \otimes {\mathcal P}(m_n)\to {\mathcal
P}(m_1+\ldots+m_n) , \text{ all } n, m_1,\ldots,m_n
\end{equation*}
called composition maps.

The composition maps can be interpreted as follows: The map
$\mu_{n;m_1,\ldots,m_n}$ combines $n$ operations $q_1,\ldots,q_n$
of $m_i (i=1,\ldots,n)$ arguments via an operation $p$ of $n$
arguments. The result, also denoted by $p(q_1,\ldots,q_n),$ is an
operation of $m_1+\ldots+m_n$ arguments.

\medskip

The unit map $1:\mbox{Id}\to F_{\mathcal P}$, or equivalently $1:
K\to{\mathcal P}(1)$, is given by an element $\id\in{\mathcal
P}(1)$ such that $\mu(\id \otimes p)=p=\mu(p\otimes
\id\otimes\ldots\otimes\id)$ for all $p\in{\mathcal P}$.

\medskip

The associativity, unity and $\Sigma$-invariance conditions
following from (\ref{defoperad}), like
\begin{equation*}
\begin{split}
\mu_{n;m_1,\ldots,m_n} \bigl(p.\sigma\otimes q_1\otimes\ldots
\otimes q_n) &=
\\
\mu_{n;m_{\sigma^{-1}1},\ldots
,m_{\sigma^{-1}n}}&\bigl(p\otimes q_{\sigma^{-1}1}\otimes\ldots
\otimes q_{\sigma^{-1}n}\bigr) .\sigma_{(m_1,\ldots,m_n)}, \\ &
\text{ all } p\in{\mathcal P(n)}, q_i\in {\mathcal P(m_i)},
\sigma\in\Sigma_n\\
\end{split}
\end{equation*}
 are the original axioms of the
definition given by May\ \cite{limay}.

\bigskip

Morphisms of operads are morphisms of $\Sigma$-spaces which
respect unit and composition maps.

\end{remark}

\begin{example}

Given a vector space $V$, there is an operad ${\mathcal End}_V$
defined by
\begin{equation*}
{\mathcal End}_V(n)={\rm Hom}(V^{\otimes n},V)
\end{equation*}
with unit given by $\id_V$ and composition maps induced by the
composition of maps. The right $\Sigma_n$-action is given by
$p.\sigma = p \circ\tilde\sigma$, where $\tilde\sigma$ is the
 left $\Sigma_n$-action of place permutation on $V^{\otimes n}$
 induced by $\sigma$.

\medskip

There also is an operad ${\mathcal{C}o\mathcal{E}nd}_V$ defined by
\begin{equation*}
{\mathcal{C}o\mathcal{E}nd}_V(n)={\rm Hom}(V,V^{\otimes n}).
\end{equation*}
Now the right $\Sigma_n$-action is defined by composition with the
right $\Sigma_n$-action  on $V^{\otimes n}$ given by
$(v_1\otimes\ldots\otimes
v_n).\sigma={\sigma^{-1}}.(v_1\otimes\ldots\otimes
v_n)$.

\end{example}

\begin{remark}

For any operad ${\mathcal P}$, we can consider the given
$\Sigma$-space just as a sequence of vector spaces and forget all
 $\Sigma_n$-actions.
This leads to a "nonsymmetric" version of the concept of operads.
It appeared earlier than the concept of (symmetric) operads,
 cf.\ the notion of a comp algebra in \cite{ligers}.

\end{remark}

\begin{defn}\label{defnonsymoperad}

A $K$-linear non-$\Sigma$  operad (or nonsymmetric operad)
$\underline{\mathcal P}$ is a sequence $\underline{\mathcal
P}(n)_{n\in\N}$ of vector spaces together with maps
\begin{equation*}
\mu_{n;m_1,\ldots ,m_n}: \underline{\mathcal P}(n)\otimes
\underline{\mathcal P}(m_1)\otimes \ldots  \otimes
\underline{\mathcal P}(m_n)\to \underline{\mathcal P}(m_1+\ldots
+m_n) , \text{ all } n, m_1,\ldots,m_n
\end{equation*}
and a unit map $1: K\to\underline{\mathcal P}(1)$ fulfilling  the
associativity and unity conditions indicated above.

 \end{defn}

\begin{remark}

By convention, in order to notationally distinguish operads from
non-$\Sigma$ operads, we use underlining (cf.\cite{limss}) to
indicate that a non-$\Sigma$ operad is given.

\end{remark}

\begin{defn}
Let ${\mathcal P}$ be a $K$-linear operad (or non-$\Sigma$
operad), and let the composition maps be denoted by
\begin{equation*}
\mu_{n;m_1,\ldots,m_n}: {\mathcal P}(n)\otimes {\mathcal
P}(m_1)\otimes \ldots \otimes {\mathcal P}(m_n)\to {\mathcal
P}(m_1+\ldots+m_n) , \text{ all } n, m_1,\ldots,m_n.
\end{equation*}

Then there are defined $\circ_i$-operations
 \begin{equation*}
\circ_i : {\mathcal P}(n)\otimes {\mathcal P}(m)\to {\mathcal
P}(m+n-1) , \text{ all } n, m\geq 1,\ 1\leq i\leq n
\end{equation*}
mapping $p\otimes q$ onto
 \begin{equation*}
 \mu_{n;1,\ldots,m,\ldots,1}(p\ \otimes\stackrel{1}{\id}\otimes\ldots\otimes \stackrel{i}{q}
 \otimes\ldots\otimes
 \stackrel{n}{\id}).
\end{equation*}

\end{defn}

\begin{remark}
Similar to the $\circ_i$-operations of Gerstenhaber (cf.\
\cite{ligers}), the $\circ_i$-operations for operads fulfill,
given any iteration
\begin{equation*}
(h\circ_i p)\circ_j q:{\mathcal P}(r)\otimes{\mathcal P}(n)\otimes
{\mathcal P}(m)\to {\mathcal P}(m+n+r-2),
 \end{equation*}
 the condition
\begin{equation*}
(h\circ_i p)\circ_j q=
\begin{cases}
(h\circ_j q)\circ_{i+m-1}p &: 1\leq j \leq i -1
\\
 h\circ_i (p \circ_{j-i+1}q) &: i\leq j \leq i+ n -1  \\
(h\circ_{j-n+1} q)\circ_i p &: i+ n \leq j. \\
\end{cases}
\end{equation*}

The composition maps are determined by these operations. To define
a (non-$\Sigma$) operad it thus suffices to give
$\circ_i$-operations fulfilling the condition above, and to
specify an operad unit, see \cite{limss}, p.\ 45.

\end{remark}

\begin{remark}\label{remmonoidaloperad}

Although we consider $K$-linear operads and non-$\Sigma$ operads,
 it should be noted that non-$\Sigma$ operads can more generally
 be defined for arbitrary monoidal categories different from
$(\mbox{Vect}_K, \otimes)$, and that operads can be defined for
arbitrary symmetric monodial categories (cf.\ \cite{limss}, II.1).
The original definition of Boardman, Vogt, and May (see
\cite{limay}) made use of the category of topological spaces.

\end{remark}

\begin{remark}\label{remregularoperad}

For every  non-$\Sigma$ operad $\underline{\mathcal P}$ there is
an operad ${\mathcal P}$, called the symmetrization of
$\underline{\mathcal P}$, with ${\mathcal P}(n)=
\underline{\mathcal P}(n)\otimes_K K\Sigma_n$ (all $n$) and
composition maps induced by the maps of $\underline{\mathcal P}$
and the maps $\mu_{n;m_1,\ldots,m_n}: K\Sigma_n\otimes
K\Sigma_{m_1}\otimes\ldots\otimes K\Sigma_{m_n}\to K\Sigma_m$ (see
Section \ref{sectconventions}).

We will also say that a given operad ${\mathcal P}$ is regular
(often also called non-$\Sigma$ in the literature), if ${\mathcal
P}$ is of the form ${\mathcal P}(n)= \underline{\mathcal
P}(n)\otimes_K K\Sigma_n$ (all $n$) for some non-$\Sigma$ operad
$\underline{\mathcal P}$.

\end{remark}

\begin{defn}

The generating series of a non-$\Sigma$ operad
$\underline{\mathcal P}$ is
\begin{equation*}
f^{\underline{\mathcal P}}(t):=\sum_{n\geq 1}\bigl(\dim
\underline{\mathcal P}(n)\bigr)t^n
\end{equation*}

The generating series of an
 operad ${\mathcal P}$ is
\begin{equation*}
f^{\mathcal P}(t):=\sum_{n\geq 1}\frac{\dim {\mathcal
P}(n)}{n!}t^n
\end{equation*}

\end{defn}

\begin{example}

 The operad ${\mathcal As}$ given by ${\mathcal
As}(n)=K\Sigma_n$, $n\geq 1$  is the symmetrization of the operad
$\underline{\mathcal As}$ given by $\underline{\mathcal As}(n)=K$
(all $n\geq 1$). The generating series is
\begin{equation*}
\sum_{n\geq 1}t^n=\frac{t}{1-t}
\end{equation*}
and it will become clear in the next section, that ${\mathcal As}$
describes the operad of (non-unitary) associative algebras.

\end{example}

\end{section}

\begin{section}{Algebras and coalgebras over operads}\label{sectalgcoalg}

\begin{defn}\label{defoperadalg}

For any operad ${\mathcal P}$, a ${\mathcal P}$-algebra consists
of a vector space $A$ together with a morphism
$\gamma=(\gamma(n)): {\mathcal P}\to {\mathcal End}_A$ of operads.

\bigskip

Each $p\in {\mathcal P}(n)$ yields a multilinear operation of $n$
arguments on $A$.

Equivalently, a ${\mathcal P}$-algebra is a vector space $A$
together with a family of $\Sigma_n$-invariant morphisms
$\gamma_A(n): {\mathcal P}(n)\otimes A^{\otimes n}\to A$, called
structure maps.

Morphisms of ${\mathcal P}$-algebras $A\to B$ are $K$-linear maps
$\varphi:A\to B$ compatible with the corresponding structure maps
$\gamma_A(n), \gamma_B(n)$ , i.e.\ such that the following diagram
commutes:
\begin{equation*}
\xymatrix{
          {\mathcal P}(n)\otimes_{\Sigma_n} A^{\otimes n}
          \ar[d]_{\id\otimes_{\Sigma_n}\varphi^{\otimes n} }
          \ar[r]^{\ \ \ \ \ \ \ \ \gamma_A(n)}
          & A
           \ar[d]^{\varphi}\\
            {\mathcal P}(n)\otimes_{\Sigma_n} B^{\otimes n}
        \ar[r]_{\ \ \ \ \ \ \ \ \gamma_B(n)} & B
            }
\end{equation*}

The space $F_{\mathcal P}(V)=\bigoplus_{n=0}^{\infty} {\mathcal
P}(n)\otimes_{\Sigma_n}V^{\otimes n}$ is the (underlying space of
the) free ${\mathcal P}$-algebra generated by the space $V$. Its
structure maps are induced by the composition maps
$\mu=\mu_{n;m_1,\ldots,m_n}$.

\end{defn}

\begin{example}\label{exascom}

The free commutative associative algebra generated by a vector
space $V$ is the symmetric algebra ${\bigoplus_{n}} (V^{\otimes
n})_{S_n}$, where the sum starts with $n=0$ to get the free
unitary commutative algebra, and with $n=1$ for the non-unitary
case.

Thus, the operad ${\mathcal Com}$ whose algebras are commutative
associative algebras (not necessarily unitary), is given by
\begin{equation*}
{\mathcal Com}(n)=K
\end{equation*}
(trivial representation of $\Sigma_n$) for each $n\geq 1$, and
${\mathcal Com}(0)=0.$

The generating series is
\begin{equation*}
\sum_{n\geq 1}\frac{t^n}{n!}=\exp(t)-1.
\end{equation*}

\medskip

 Choosing
\begin{equation*}
{\mathcal As}(n)=K\Sigma_n
\end{equation*}
the (module of the) regular representation for each $n\geq 1$, one
gets the operad of associative (not necessarily unitary) algebras,
with free algebra functor
\begin{equation*}
F_{\mathcal As}(V)={\displaystyle \bigoplus_{n=1}^{\infty}}
V^{\otimes n}.
\end{equation*}

Lie algebras and Poisson algebras yield classical examples of
operads, too. The operad $\Lie$ is given by the
$(n-1)!$-dimensional $K$-spaces generated by multilinear bracket
monomials with respect to Jacobi identity and anti-symmetry, cf.\
\cite{ligk}. For $\Pois$ one needs a commutative associative
binary operation and a Lie bracket operation $[,]$. They are
related by the identity $[a,b\cdot c]=b\cdot[a,c]+[a,b]\cdot c$.
(See Section \ref{sectfreeoperad} for the construction of operads
by generators and relations.)

\end{example}

\begin{remark}
We note that the operads $\Com$ and  $\As$ and also $\Pois$ and
$\Mag$ (the operad of magma algebras) are defined in a way such
that their algebras are not necessarily unitary.

For example, the elements of the free $\Com$-algebra $F_{\mathcal
Com}(V_X)$, $V_X$ the
  vector space with basis $X$,
are polynomials in commutative variables from $X$ without constant
terms, $F_{\mathcal Com}(V_X)=\overline{K[X]}$ (called the
augmentation ideal of the free unitary commutative algebra).

\end{remark}

\begin{defn} (cf.\cite{lifra}, 1.4.)

A ${\mathcal P}$-algebra $A$ is called  nilpotent, if for
sufficiently large $n$,
\begin{equation*}
\gamma_A(n)(\mu\otimes a_1\otimes\ldots a_{n})=0
\end{equation*}
for all $a_1,\ldots,a_n\in A$, $\mu\in{\mathcal P}(n)$.

\medskip

 An ideal of  a ${\mathcal P}$-algebra $A$ is a subspace $I$ of
$A$, such that
\begin{equation*}
\gamma_A(n)(\mu\otimes a_1\otimes\ldots \otimes a_{n-1}\otimes
b)\in I
\end{equation*}
 for all $\mu\in{\mathcal P}(n), a_i\in A, b\in I$.

\end{defn}

\begin{remark}

Given an ideal $I$ in a ${\mathcal P}$-algebra $A$, the quotient
$A/I$ is again a ${\mathcal P}$-algebra.

If $A$ is a ${\mathcal P}$-algebra together with a sequence $I_n,
n\geq 1$ of ideals, such that the ${\mathcal P}$-algebras $A/I_n$
are nilpotent, one can construct the completion $\widehat
A=\underleftarrow{\lim}\ A/I_n$ of $A$ with respect to the
topology given by $(I_n)$.

\end{remark}

\begin{defn}

A complete ${\mathcal P}$-algebra is a ${\mathcal P}$-algebra $A$
together with a sequence $I_n, n\geq 1$ of ideals, such that
$A/I_n$ is nilpotent for all $n$ and such that
$A=\underleftarrow{\lim}\ A/I_n$.

The free complete ${\mathcal P}$-algebra generated by a vector
space $V$ is given by
\begin{equation*}
\begin{split}
&\widehat F_{\mathcal P}(V):= \prod_{n=1}^{\infty} {\mathcal
P}(n)\otimes_{\Sigma_n}V^{\otimes n}.
\\
\end{split}
\end{equation*}
If $V=V_X$, $X=\{x_1,\ldots, x_m\}$, the elements of $\widehat
F_{\mathcal P}(x_1,\ldots,x_m)$ are called ${\mathcal P}$-power
series in variables $x_1,\ldots,x_m$.

\end{defn}

\begin{example}

For example, $\Com$-power series in variables $x_1,\ldots,x_m$ are
power series (without constant terms) in commuting variables,
$\widehat F_{\mathcal Com}(V_X)=\overline{K[[X]]}$.

Similarly, $\As$-power series in variables $x_1,\ldots,x_m$ are
power series (without constant terms) in non-commuting variables.

\end{example}

\begin{remark}\label{remcooperad}

Dual to the notion of an operad, there is the notion of a
co-operad ${\mathcal Q}=({\mathcal Q}(n))_{n\in \N}$.
  The definition is analogous to the definition of
operads. Now the endofunctor on Vect$_K$ given by
\begin{equation*}
\begin{split}
&\bigoplus_{n=0}^{\infty} \bigl({\mathcal Q}(n)\otimes V^{\otimes
n}\bigr)^{\Sigma_n}
\\
\end{split}
\end{equation*}
has to be a comonad.  (The space of invariants is used.)

 We get (co-)composition maps ${\mathcal
Q}(i_1+\ldots+i_n)\to {\mathcal Q}(n)\otimes {\mathcal
Q}(i_1)\otimes\ldots\otimes{\mathcal Q}(i_n)$.

If ${\mathcal P}^*=({\mathcal P}(n))^*_{n\in \N}$ denotes the
$K$-linear dual of a $\Sigma$-space ${\mathcal P}$, and if we
assume all ${\mathcal P}(n)$ to be finite dimensional, the axioms
for ${\mathcal P}$ being an operad correspond to the axioms for
${\mathcal P}^*$ being a co-operad.

\end{remark}

\begin{defn}\label{defpcoalg}(cf.\cite{limss}, p.165)

Let ${\mathcal P}$ be an operad. We assume that all ${\mathcal
P}(n)$ are finite dimensional.

A ${\mathcal P}$-coalgebra consists of a vector space $C$ together
with a morphism $\lambda_C=(\lambda_C(n)):{\mathcal P}\to
{\mathcal Co\mathcal{E}nd}_C$ of operads.

\bigskip

Equivalently, a ${\mathcal P}$-coalgebra is a vector space $C$
together with a family of $\Sigma_n$-invariant morphisms
$\lambda_C(n): {\mathcal P}(n)\otimes C\to C^{\otimes n}$, called
structure maps.

Morphisms of ${\mathcal P}$-coalgebras $C\to D$ are $K$-linear
maps $\varphi:C\to D$ compatible with the corresponding structure
maps.

A ${\mathcal P}$-coalgebra $C$ is called connected (see\
\cite{liff}, \S 4) or (co-)nilpotent (see \cite{limss}, p.165), if
the following condition holds:
\begin{equation*}
\text{ for all } c\in C \text{ there is } N\in \N \text{ such that
for } n > N, \lambda_C(n)(\mu\otimes c)=0, \text{ all
}\mu\in{\mathcal P}(n).
\end{equation*}

\end{defn}

\begin{remark}\label{remcofree}

The space $F^c_{\mathcal P}(V)=\bigoplus_{n=0}^{\infty}
\bigl({\mathcal P}^*(n)\otimes V^{\otimes n}\bigr)^{\Sigma_n}$ is
the underlying space of the cofree (co-)nilpotent ${\mathcal
P}$-coalgebra (co-)generated by the space $V$.
 Its structure maps are
induced by the dual composition maps
\begin{equation*}
{\mathcal P}^*(i_1+\ldots+i_n)\to {\mathcal P}^*(n)\otimes
{\mathcal P}(i_1)^*\otimes\ldots\otimes{\mathcal P}^*(i_n)
\end{equation*}
 of
(\ref{remcooperad}).

If the characteristic of $K$ is 0, for any vector space $V$ over
$K$ there is a projection $\pi:V\to V^{\Sigma_n}$, $v\mapsto
\frac{1}{n!}\sum_{\sigma\in \Sigma_n}\sigma(v)$ onto the
 space of invariants.

  Hence $V_{\Sigma_n}\cong V/\ker\pi\cong
 V^{\Sigma_n}$.
Thus via (non-canonical) vector space isomorphisms between
${\mathcal P}^*(n)$ and ${\mathcal P}(n)$ it is possible to
provide the free ${\mathcal P}$-algebra $F_{\mathcal
P}(V)=\bigoplus_{n=0}^{\infty} {\mathcal
P}(n)\otimes_{\Sigma_n}V^{\otimes n}$ with the structure of a
cofree (co-)nilpotent coalgebra.

The construction of non-nilpotent cofree coalgebras is more
complicated (see \cite{lifox}).
\end{remark}

\begin{example}\label{excofree}

Let $K=\Q$. For ${\mathcal P}={\mathcal As}$, and $V_X$ vector
space with basis $X$, the words $w_1.w_2\ldots w_r$, $r\geq 1$,
$w_i\in X$ form a basis of $F_{\mathcal As}(V_X)$. Now the
(co-)nilpotent cofree coalgebra-structure is given by
deconcatenation, that is, for $\mu\in {\mathcal As}(n),
w=w_1\ldots w_r,$ the image of $\mu\otimes w$ under $\lambda(n)$
is
\begin{equation*}
\sum_{1\leq j_1<\ldots<j_{n-1}<r}(w_1\ldots
w_{j_1})\otimes\ldots\otimes(w_{j_{n-1}+1}\ldots w_{r}).
\end{equation*}

Since all operations of ${\mathcal As}$ can be built up from one
binary multiplication $\tilde\mu\in {\mathcal As}(2)$ (together
with $\tilde\mu.\tau_1$, where $\tau_1$ is the transposition), the
coalgebra-structure is determined by the images of
$\tilde\mu\otimes w$, given by
\begin{equation*}
\Delta(w):=\sum_{j=1}^{r-1}(w_1\ldots w_{j})\otimes(w_{j+1}\ldots
w_{r}).
\end{equation*}
It is possible to extend the coalgebra structure onto $K\langle
X\rangle=K 1\oplus F_{\mathcal As}(V_X)$, where 1 is the empty
word, by setting
\begin{equation*}
\Delta(w):=\sum_{j=0}^{r}(w_1\ldots w_{j})\otimes(w_{j+1}\ldots
w_{r}).
\end{equation*}
Usually (the vector space) $K\langle X\rangle$ together with
$\Delta$ is called the standard tensor coalgebra.

\end{example}

\begin{example}\label{excomcofree}
 For ${\mathcal P}={\mathcal Com}$,
 the cofree (co-)nilpotent coalgebra occurs as a subcoalgebra of
the standard tensor coalgebra, namely the subcoalgebra of
symmetric tensors.
\end{example}

\end{section}

\end{chapter}

\begin{chapter}{Some combinatorics of trees}

\begin{section}{Abstract and planar trees}

For computations in not necessarily associative algebras, and also
in operad theory, trees are very useful to symbolize the ways of
associating variables (or arguments of an operation).

We have to make a difference between several types of trees. First
we recall the notions of rooted trees and planar rooted trees. We
skip the definition of graphs. A naive notion of a graph will
suffice. For a more sophisticated notion, involving half-edges,
see \cite{limss}, \S 5.3.

\begin{defn}

 A finite connected graph $\emptyset \neq T=({\rm
Ve}(T), {\rm Ed}(T))$, with a distinguished vertex $\rho_T$, is
called an abstract rooted tree, if for every vertex
$\lambda\in{\rm Ve}(T)$ there is exactly one path connecting
$\lambda$ and $\rho_T$.

The vertex $\rho_T$ is called the root of $T$. Thinking of the
edges as
 oriented towards the root, at each vertex there are incoming
 edges and one outgoing edge.
Modifying the standard convention (but cf.\ also \cite{limss},
p.50), we add to the root an outgoing edge
 that is not connected to any further vertex. (If we want to exclude
 this edge, we speak of the other edges as inner edges.)

We denote by $\ATree$ the set of abstract rooted trees.

\end{defn}

\begin{remark}

In the literature, abstract rooted trees are often only called
rooted trees. Since we are also going to deal with planar rooted
trees, and since planar rooted trees are not special (abstract)
rooted trees, we stress the word 'abstract'. We may skip the word
'rooted', because we are only going to consider rooted trees.

\end{remark}

\begin{defn}
The height of a vertex $\lambda\in{\rm Ve}(T)$ is the number of
edges separating it from $\rho_T$. The height of a rooted tree $T$
is the maximum height of its vertices.

 At a given vertex $\lambda$,
the number $n$ of incoming edges is called the arity
ar$_{\lambda}$ of $\lambda$. We write the set ${\rm Ve}(T)$ of
vertices as a disjoint union $\bigcup_{n\in\N} {\rm Ve}^n(T)$.
The vertices of arity $0$ are called leaves, and we denote ${\rm
Ve}^0(T)$ by ${\rm Le}(T)$.
 The elements of ${\rm Ve}^*(T)={\rm
Ve}(T)-{\rm Le}(T)$ are called internal vertices of $T$.

 A tree $T$ is called binary, if ${\rm Ve}^*(T)={\rm Ve}^2(T)$, i.e.\ if ar$_{\lambda}=2$
for all internal vertices $\lambda$.

 An abstract rooted tree $T$ together with a
chosen order of incoming edges at each vertex is called a planar
rooted tree (or ordered rooted tree), see Example
(\ref{exsametree}). We denote by $\PTree$ the set of planar rooted
trees.

\end{defn}

\begin{remark}\label{remsequences}

It is well-known that
the number of planar rooted trees with $n$ vertices is the $n$-th
Catalan number
\begin{equation*}
 c_n=\frac{(2(n-1))!}{n!(n-1)!}=\sum_{l=1}^{n-1}c_lc_{n-l}.
\end{equation*}
The sequence of Catalan numbers is
\begin{equation*}
 c_1=1, c_2=1, c_3=2,c_4=5, c_5=14, c_6=42, c_7=132,
 c_8=429, c_9=1430, \ldots
\end{equation*}
with generating series $f(t)=\sum_{n=1}^{\infty}c_nt^n$ given by
\begin{equation*}
\frac{1-\sqrt{1-4t}}{2}.
\end{equation*}

The numbers $c_n$ also count the number of planar binary rooted
trees with $n$ leaves (or $2n-1$ vertices).

\medskip

The numbers $a_n$ of abstract rooted trees with $n$ vertices, for
$n=1,\ldots,11$, are
\begin{equation*}
1,1,2,4,9,20,48,115,286,719,1842.
\end{equation*}
The generating series $f(t)=\sum_{n=1}^{\infty}a_nt^n$ fulfills
the equation
\begin{equation*}
f(t)=\frac{t}{\prod_{n\geq 1}(1-t^n)^{a_n}}
\end{equation*}
or equivalently the equation
\begin{equation*}
f(t)=t\exp\Bigr(\sum_{k\geq 1}\frac{f(t^k)}{k}\Bigl)
\end{equation*}
(cf.\ \cite{lihar}).

\medskip

The numbers of abstract binary rooted trees with $n$ leaves (or
$2n-1$ vertices), for $n=1,\ldots,11$, are
\begin{equation*}
1,1,1,2,3,6,11,23,46,98,207,451
\end{equation*}
with generating series $f(t)$ given by (cf.\ \cite{lipe}) the
equation
\begin{equation*}
f(t)=t+\frac{1}{2}f(t)^2 +\frac{1}{2}f(t^2).
\end{equation*}
The series can thus be written in the form
\begin{equation*}
\begin{split}
f(t)&=1-\sqrt{1-f(t^2)-2t}\\
&=1-\sqrt{\sqrt{1-f(t^4)-2t^2}-2t}=\ldots\\
\end{split}
\end{equation*}

Information on these integer sequences can be found in the On-Line
Encyclopedia of Integer Sequences \cite{lisl}.

\end{remark}

\begin{example}\label{exsametree}
We call rooted trees just trees, for short.
 We draw the root at the bottom (or bottom left, more exactly).
For every vertex of a planar tree, the chosen order of incoming
edges corresponds to an ordering of edges from left to right.
Every drawing of a tree provides us with a planar structure, which
we have to forget when dealing with abstract trees.

The drawings
 $
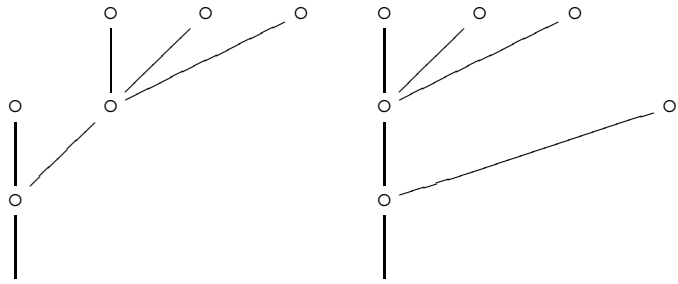
\xymatrix{  & \ar@{-}[d]\circ & \ar@{-}[dl]\circ & \ar@{-}[dll]\circ \\
 \circ \ar@{-}[d]& \circ \ar@{-}[dl]& \\ \circ\ar@{-}[d] & & \\ & \\}\ \ \ \ \
 \xymatrix{ \circ \ar@{-}[d] & \ar@{-}[dl]\circ & \ar@{-}[dll]\circ&  \\
  \circ \ar@{-}[d]& & &\ar@{-}[dlll]\circ \\ \circ\ar@{-}[d] & & \\
  &\\}$

represent the same abstract tree $T$ (of height 2), but different
planar trees $T^1, T^2$.

\medskip

In our drawing, we can put labels
 at
the vertices:
 $
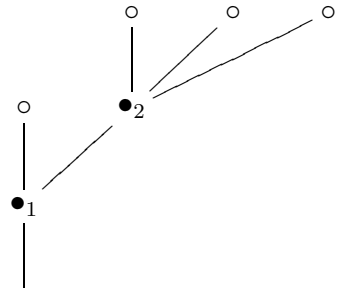
\xymatrix{  & \ar@{-}[d]\circ & \ar@{-}[dl]\circ & \ar@{-}[dll]\circ \\
 \circ \ar@{-}[d]& \bullet_2 \ar@{-}[dl]& \\ \bullet_1\ar@{-}[d] & & \\ &
 \\}$

Here we have used different sets of labels for leaves and internal
vertices.

The following definition of labeled trees and admissibly labeled
trees is useful to include labeled trees in operad theory as well
as other types of labelings.

\end{example}

\begin{defn}\label{defadmiss}

Let $M$ be a set and $T$ a planar (or abstract) tree.

Then a labeling of $T$ is a map $\nu:{\rm Ve}(T)\to M$. The tree
$T$ together with such a labeling is called a labeled tree.

\medskip

Let a collection $M_0, M_1, M_2,\ldots$ of
 sets be given, and let $M=\bigcup_{k\in\N}
M_k$.

A labeling  $\nu:{\rm Ve}(T)\to M$ of a planar (or abstract) tree
$T$ is called admissible, if the restrictions $\nu\vert{\rm
Ve^k}(T)$ are maps ${\rm Ve^k}(T)\to M_k$, i.e.\ it holds that:
\begin{equation*}
\nu(\lambda)\in M_k \text{ if \ ar}_{\lambda}=k.
\end{equation*}

The set  of planar rooted trees $T\in\PTree$ with admissible
labeling from $(M_k)_{k\in\N}$ is denoted by $\PTreeM$.

Given a fixed tree $T \in\PTree$, we define
\begin{equation*}
\TMseq:=\{ (T,\nu):\ \nu \text{\ admissible\ }\} \subseteq\PTreeM.
\end{equation*}

We make the same definitions for abstract rooted trees with
admissible labelings. The set of these trees is denoted by
$\ATreeM=\bigcup_{T\in\ATree} \TMseq.$

\end{defn}

\begin{remark}

We can identify the set $\PTree$ with the set
$\PTree\{M_k=\{\circ\}, k\in\N\}$, i.e.\ we consider non-labeled
trees as trivially labeled trees.

\medskip

In our drawings, we will often use labels $\bullet$ and $\circ$ to
distinguish between vertices that count and vertices that do not
count for a degree function.

More generally, one can work with weighted labels. Apart from the
$\circ$-label (weight 0), we are only going to use labels of
weight 1. The degree of a tree is then the number of all vertices
that are not labeled by a $\circ$.

\medskip

We mention that one can also associate levels to all vertices of
a given tree, to distinguish for example between the trees

 $
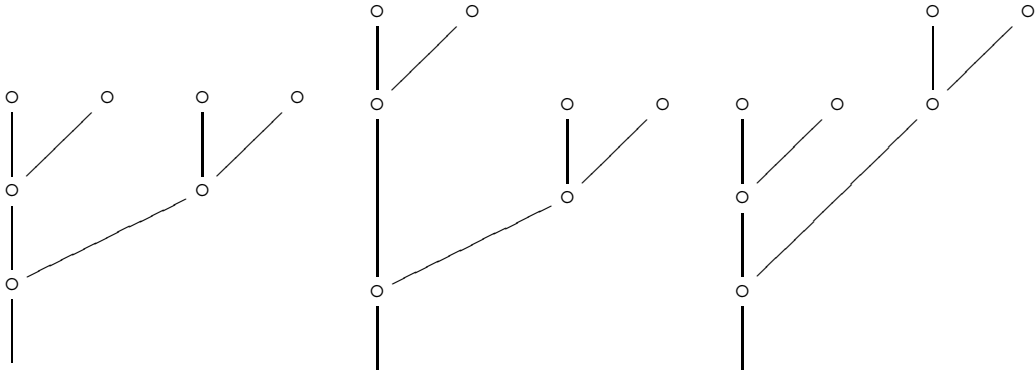
\xymatrix{ & \\
  \ar@{-}[d]\circ & \ar@{-}[dl]\circ & \ar@{-}[d]\circ
 & \ar@{-}[dl]\circ \\
 \circ \ar@{-}[d]& & \circ \ar@{-}[dll]& \\ \circ\ar@{-}[d] & & \\ & \\}\ \ \ \ \
 \xymatrix{
  \ar@{-}[d]\circ & \ar@{-}[dl]\circ  \\
 \circ \ar@{-}[dd]& & \ar@{-}[d]\circ
 & \ar@{-}[dl]\circ\\ & & \circ \ar@{-}[dll]& \\ \circ\ar@{-}[d] & & \\ & \\}\ \ \ \ \
  \xymatrix{ & & \ar@{-}[d]\circ
 & \ar@{-}[dl]\circ \\
  \ar@{-}[d]\circ & \ar@{-}[dl]\circ & \ar@{-}[ddll]\circ
 \\
 \circ \ar@{-}[d]& \\ \circ\ar@{-}[d] & & \\ & \\}\ \ \ \ \
 $
We are going to mention this type of trees only one time (in Section
\ref{sectfreeoperad}).

\end{remark}

\begin{defn}\label{deflogarithmic}
Let $a_n, n\geq 1,$ be a sequence of integers with generating
series $f(t)=\sum_{n=1}^{\infty} a_n t^n$.

The logarithmic derivative of $f(t)$ is the series
\begin{equation*}
g(t):=\frac{\partial}{\partial t}\log\bigl(1+f(t)\bigr),
\end{equation*}
and we say that the sequence $a'_n, n\geq 1,$ with
$\sum_{n=1}^{\infty} a'_n t^n=t\cdot g(t)$ is obtained from $a_n,
n\geq 1,$ by logarithmic derivation.

\end{defn}

\begin{example}\label{exlogcatalan}

The sequence of $\log$-Catalan numbers $c'_n$, starting with
\begin{equation*}
1, 1, 4, 13, 46, 166, 610, 2269, 8518, 32206,\ldots
\end{equation*}
has the generating series
\begin{equation*}
\frac{2t}{3\sqrt{1-4t}-1+4t\ .}
\end{equation*}
Since
\begin{equation*}
\frac{2}{3\sqrt{1-4t}-1+4t}=\frac{\partial}{\partial
t}\log\bigl(\frac{3-\sqrt{1-4t}}{2}\bigr),
\end{equation*}
it is obtained by logarithmic derivation from the Catalan numbers
$c_n, n\geq 1$.

In the set of all planar (rooted) trees with $n$ vertices,
 the number of vertices with even arity is given
by $c'_n$. The corresponding numbers of vertices with odd arity
have the generating series
\begin{equation*}
\begin{split}
\sum_{n=1}^{\infty} n c_n t^n &- \sum_{n=1}^{\infty} c'_n t^n =
\frac{t}{\sqrt{1-4t}}-\frac{2t}{(3-\sqrt{1-4t})\sqrt{1-4t}}
\\
 &=\frac{t(1-\sqrt{1-4t})}{(3-\sqrt{1-4t})\sqrt{1-4t}}
 =t^2 + 2t^3 + 7t^4 + 24t^5 + 86 t^6 + 314 t^7 + \ldots\\
\end{split}
\end{equation*}
(see \cite{lisl} A026641, \cite{lids} p.\ 258).

For example, in the set

$
 
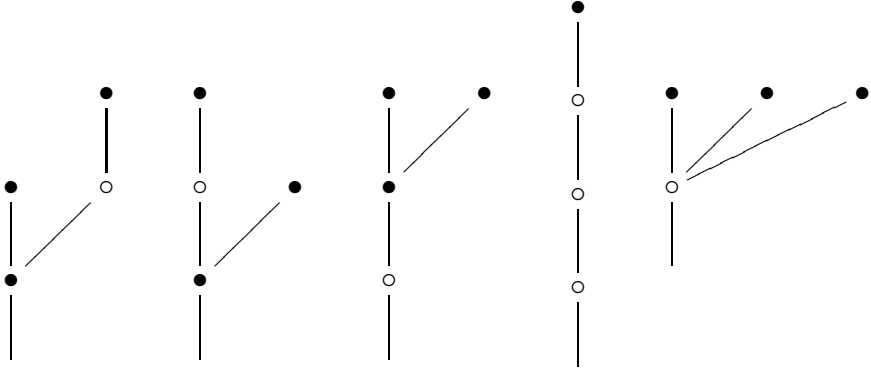
\xymatrix{ \\ & {\bullet}\ar@{-}[d] \\ {\bullet}\ar@{-}[d]& {\circ}\ar@{-}[dl]\\
{\bullet}\ar@{-}[d]\\ \\}\ \ \ \ \ \
 \xymatrix{ \\ {\bullet}\ar@{-}[d] \\ {\circ}\ar@{-}[d]& {\bullet}\ar@{-}[dl]\\
{\bullet}\ar@{-}[d]\\ \\}\ \ \ \ \ \
 \xymatrix{ \\ {\bullet}\ar@{-}[d] & {\bullet}\ar@{-}[dl]\\ {\bullet}\ar@{-}[d]\\
{\circ}\ar@{-}[d]\\ \\}\ \ \ \ \ \
 \xymatrix{{\bullet}\ar@{-}[d] \\ { \circ}\ar@{-}[d] \\
{\circ}\ar@{-}[d]\\ {\circ}\ar@{-}[d]\\
\\}\ \ \ \ \ \
 \xymatrix{\\ {
\bullet}\ar@{-}[d] & { \bullet}\ar@{-}[dl] &
 { \bullet}\ar@{-}[dll]\\ {\circ}\ar@{-}[d]\\ \\}$

of planar trees with $n=4$ vertices, we count seven vertices with
odd arity and $c'_4=13$ vertices with even arity.

\end{example}

\bigskip

\begin{defn}\label{defsubstitution}

Let $T^1, T^2$ be planar (or abstract) trees and let $b$ be a leaf
of $T^1$. Then the substitution of $T^2$ in $T^1$ at $b$, denoted
by $T^1 \circ_b T^2$, is obtained by replacing the leaf $b$ of
$T^1$ by the root of $T^2$.

A word $T^1.T^2\ldots T^r$ (or an ordered tuple $(T^1,T^2,\ldots
T^r)$, not necessarily non-empty) of planar trees is called a
planar forest. A disjoint union (or unordered tuple) of abstract
trees is called an abstract forest.

In both cases, given a forest $T^1\ldots T^n$ of $n\geq 0$ trees,
together with a label $\rho\in M$, there is a tree
$T=\vee_{\rho}(T^1\ldots T^n)$ defined by introducing a new root
of arity $n$ and grafting the trees $T^1,\ldots,T^n$ onto this new
root. The new root gets the label $\rho$. In the planar case, the
specified order determines the order of incoming edges at
$\rho_T$. The tree $T$ is called the grafting of $T^1\ldots T^n$
over $\rho$.

If there is no choice for a label $\rho$ (i.e.\ there is only one
label available) we simply write $\vee(T^1\ldots T^n)$. In the
literature, $\vee(T^1\ldots T^n)$ is often denoted by
$B_+(T^1\ldots T^n)$, e.g.\ in \cite{likrdiv}.

\end{defn}

\begin{defn}

 A (planar or abstract) tree $T$ is called reduced, if
ar$_{\lambda}\neq 1$ for all $\lambda\in {\rm Ve}(T)$. The set
$\PRTree$ of planar reduced trees is identified with the set
$\PTreeM$ given by $M_1=\emptyset$, and $M_k$ (for $k=0$ or $k\geq
2)$ a one element set, $\{\circ\}$ say.

\end{defn}

\goodbreak

\begin{example}

 The grafting $\vee(T^1, T^2)$ of $T^1$, $T^2$ from (\ref{exsametree}) is just

 $
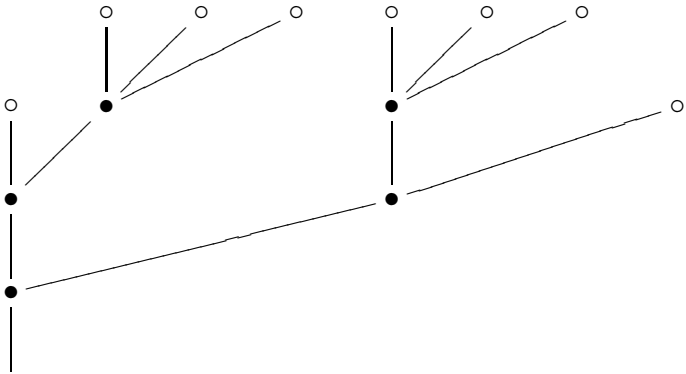
\xymatrix{ & \ar@{-}[d]\circ & \ar@{-}[dl]\circ &
 \ar@{-}[dll]\circ
&\circ \ar@{-}[d] & \ar@{-}[dl]\circ & \ar@{-}[dll]\circ&
 \\
 \circ \ar@{-}[d] & \bullet \ar@{-}[dl]& & & \bullet \ar@{-}[d]& & &\ar@{-}[dlll]\circ\\
  \bullet\ar@{-}[d] & & & &
  \bullet\ar@{-}[dllll] & & \\
  \bullet\ar@{-}[d]& & & & \\
  \\}$

It is a reduced tree.

The following tree is not reduced: $\xymatrix{\\ {
\bullet}\ar@{-}[d]\\ {\bullet}\ar@{-}[d]\\
\\}$

The tree $\vee(\underbrace{\bullet\ldots\bullet}_n)$ is called
$n$-corolla:

$\overbrace
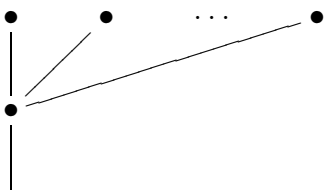
{\xymatrix{\\ \bullet\ar@{-}[d] & \bullet\ar@{-}[dl] &
\ldots & \bullet\ar@{-}[dlll]\\ {\bullet}\ar@{-}[d]\\
\\}}^n$

Examples for binary trees can be drawn as follows:

 $
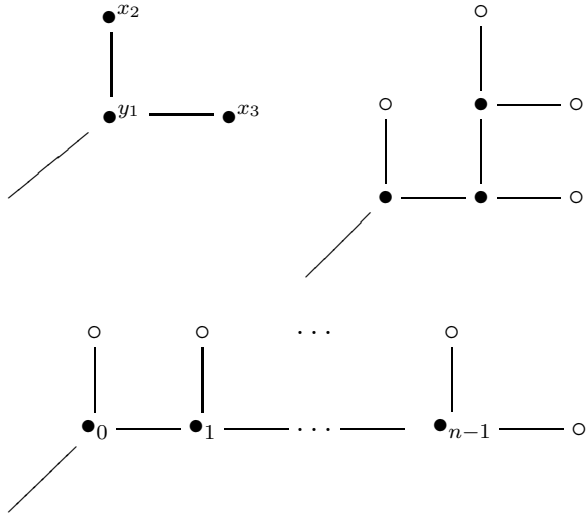
\xymatrix{ \\ & {\ \
\bullet^{x_2}}\ar@{-}[d] & \\ & {\ \ \bullet^{y_1}}\ar@{-}[dl]&{
\bullet^{x_3}}\ar@{-}[l]\\ \\} \ \
 \xymatrix{ \\&  & {\circ}\ar@{-}[d] \\ & {\circ}\ar@{-}[d] &
{\bullet}\ar@{-}[d]\ar@{-}[r]  & {\circ} \\
 & {\bullet}\ar@{-}[dl] &{\bullet}\ar@{-}[l]  &{\circ}\ar@{-}[l]  \\ \\} \
$

\medskip

$
\xymatrix{ & {\circ}\ar@{-}[d] &  {\circ}\ar@{-}[d] & \ldots &
 {\circ}\ar@{-}[d]\\
 & {\bullet_0}\ar@{-}[dl] &{\bullet_1}\ar@{-}[l]\ar@{-}[r]&\ldots &{\ \ \bullet_{n-1}}\ar@{-}[l]
  &{\circ}\ar@{-}[l]  \\ \\} $

If no labels are given, the last binary tree in the picture above
is called right comb (of height $n$). Left combs are defined
analogously.

\end{example}

\begin{defn}

Let $T$ be a planar (or abstract) tree, and let $v\in {\rm Ve}(T)$
be a vertex. The vertex $v$ determines a subgraph $S$ of $T$,
called the full subtree of $T$ with root $v$, such that: $v\in{\rm
Ve}(S)$, and for every vertex $w\in{\rm Ve}(S)$ all incoming edges
(vertices included) of $w$ in $T$ belong also to $S$.

\end{defn}

\begin{example}

$\newline$
 The tree
 $
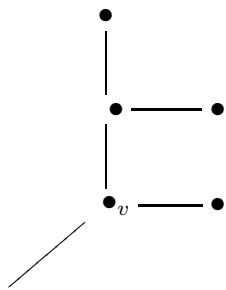
\xymatrix{ & \bullet\ar@{-}[d]\\ & {\ \
\bullet}\ar@{-}[d] &{ \bullet}\ar@{-}[l] \\ & {\ \
\bullet_v}\ar@{-}[dl]&{ \bullet}\ar@{-}[l]\\ \\} \ \ $ is a
 full subtree of the tree
$
 \xymatrix{  & {\bullet}\ar@{-}[d] &{\bullet}\ar@{-}[d] \\ & {\bullet}\ar@{-}[d] &
{\bullet}\ar@{-}[d]\ar@{-}[r]  & {\bullet} \\
 & {\bullet}\ar@{-}[dl] &{\ \ \bullet_v}\ar@{-}[l]  &{\bullet}\ar@{-}[l]  \\ \\} \
$
\end{example}

\begin{remark}\label{remdegraftmirror}

Every (abstract or planar) tree $T$ is the grafting of a forest
$\neg T$ uniquely determined by $T$. If the root $\rho$ of $T$ has
arity $n$, the forest $\neg T$ consists of the $n$ subtrees of $T$
given by the $n$ vertices connected to $\rho$ by one edge.
Especially, all trees in the forest $\neg T$ have heights less
than the height of $T$.

If $T$ is a planar binary tree with at least two leaves, then $T$
can be uniquely written as the grafting
$T=T^l\vee_{\rho}T^r:=\vee_{\rho}(T^l,T^r)$ of its left tree $T^l$
and its right tree $T^r$.

It is sometimes useful to call the empty set $\emptyset$ a tree
(and to define its height to be $-1$); obviously $\emptyset$
cannot be written as a grafting. Let ${\rm PRTree}':={\rm
PRTree}\cup\{\emp\}$.

There is a canonical de-grafting map $\neg$ from labeled non-empty
trees to forests of labeled trees (given by deleting the root
together with its label).
 The operator $\neg$ is often denoted by $B_{-}$.

Given a planar tree $T$, there is a unique tree $\bar T$,
recursively defined

\begin{equation*}
\overline{\vee(T^1\ldots T^n)}=\vee(\bar T^n\ldots \bar T^1),
\end{equation*}
 where $\overline{\vee(\emptyset)}= \vee(\emptyset)$ and
$\bar\emptyset=\emptyset$. In other words, $\bar T$ is obtained by
mirroring $T$ along the root axis.

It holds that $\overline{\bigl(\bar{T}\bigr)}=T$. The trees $T^1$
and $T^2$ from (\ref{exsametree}) are in correspondence via
$T\mapsto\bar T.$

\end{remark}

\end{section}

\begin{section}{Strings, reductions, and cuts}\label{sectstring}

 There is a correspondence
between planar rooted trees and (irreducible) parenthesized
strings. We sketch this correspondence in the general setting of
 admissibly labeled planar trees, where $M_0=\{x_1,x_2,\ldots \},
 M_k=\{\mu_1,\mu_2,\ldots \} (k\geq 1)$.

Kreimer's definition of irreducible parenthesized words in
\cite{likr} is completely analogous, but for the special case
where $x_i=\mu_i$ (all $i$).

\begin{remark}

Given an admissibly labeled planar tree $T$, we recursively
construct the corresponding parenthesized string.

If $T$ consists of its root $\rho_T$, then $\rho_T$ is a leaf
labeled by some $x_{i_1}$. The corresponding parenthesized string
is $(x_{i_1})$, i.e.\ an opening bracket followed by the letter
$x_{i_1}$ followed by a closing bracket.

Else, let $\rho_T$ be labeled by $\mu_{i_1}$, and let $T^1\ldots
T^n$ be the forest of labeled trees which remains after removing
the root with its incoming edges. Assume that $T_j$ has got the
corresponding  parenthesized string $w_j$, all $j$. Then
$(\mu_{i_1} w_1 \ldots w_n )$ is the parenthesized string
associated to $T$.

We get a string of letters and balanced brackets such that the
leftmost opening bracket is matched by the rightmost closing
bracket (irreducibility) and such that each letter has exactly one
opening bracket on its lefthandside.

It is easy to see that the tree can be reconstructed from its
string.

Reducible words are defined by concatenation of irreducible ones,
thus they correspond to forests.

For abstract trees, there is a completely similar construction.
The only difference is that some words have to be identified due
to the missing order of incoming edges.

\end{remark}

\begin{example}

The parenthesized strings
 $(x_1), (\mu_1(x_2)), (\mu_1(x_2)(x_3))$,
 $(\mu_1(\mu_2(x_3)))$, and $(\mu_1(x_2)(x_3)(x_4))$
 represent the trees

$
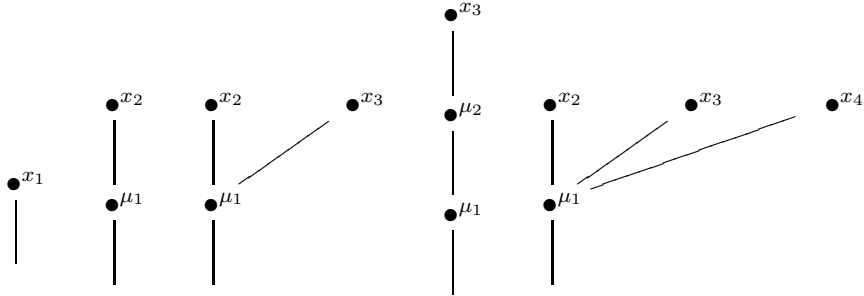
\xymatrix{\\ \\ {\ \ \bullet^{x_1}}\ar@{-}[d]\\ \\} \ \
\xymatrix{\\ {\ \ \bullet^{x_2}}\ar@{-}[d]\\ {\ \
\bullet^{\mu_1}}\ar@{-}[d]\\ \\} \ \
 \xymatrix{ \\ {\ \
\bullet^{x_2}}\ar@{-}[d] & {\ \ \bullet^{x_3}}\ar@{-}[dl]\\ {\ \
\bullet^{\mu_1}}\ar@{-}[d]\\ \\} \ \
 \xymatrix{ {\ \
\bullet^{x_3}}\ar@{-}[d] \\ {\ \ \bullet^{\mu_2}}\ar@{-}[d]\\ {\ \
\bullet^{\mu_1}}\ar@{-}[d]\\ \\} \ \
 \xymatrix{\\ {\ \
\bullet^{x_2}}\ar@{-}[d] & {\ \ \bullet^{x_3}}\ar@{-}[dl] &
 {\ \ \bullet^{x_4}}\ar@{-}[dll]\\ {\
\ \bullet^{\mu_1}}\ar@{-}[d]\\ \\}$

 An empty pair of brackets
without label is also allowed. It represents the empty tree.

\end{example}

\begin{example}\label{exmalcev}

In the following  we consider binary trees. Let
$M_0=\{x_1,x_2,\ldots\}$, and let $M_2$ be a one element set.

For every pair of brackets, the position of the closing bracket is
forced once the position of the opening bracket is given. Thus one
can omit the brackets and just use a letter $c$
 to mark an opening bracket given by an internal vertex.

For example the binary trees

 $
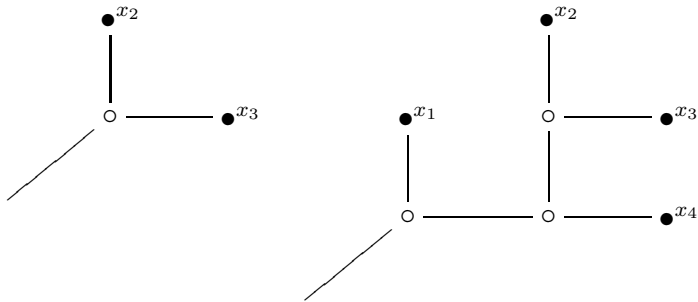
\xymatrix{ \\ & {\ \
\bullet^{x_2}}\ar@{-}[d] & \\ & {\circ}\ar@{-}[dl]&{
\bullet^{x_3}}\ar@{-}[l]\\ \\} \ \
 \xymatrix{ \\&  & {\ \
\bullet^{x_2}}\ar@{-}[d] \\ & {\ \ \bullet^{x_1}}\ar@{-}[d] &
{\circ}\ar@{-}[d]\ar@{-}[r]  & {\bullet^{x_3}} \\
 & {\circ}\ar@{-}[dl] &{\circ}\ar@{-}[l]  &{\bullet^{x_4}}\ar@{-}[l]  \\ \\} \
\
$

can be represented by the strings
 $c x_2 x_3$, $cx_1c^2x_2x_3x_4$.

\end{example}

\begin{remark}

Since the free magma generated by a set of variables $X$ consists
of parenthesized strings given by planar binary trees, we can
 call the set of planar binary trees with leaves labeled by $X$
 the free magma generated by $X$.

If a field $K$ is given, we can pass from the free magma generated
by $X$ to the free magma algebra (similarly to passing from
semi-groups or groups to semi-group algebras or group-algebras).

The representation given in Example (\ref{exmalcev}) is the Malcev
representation of the free magma algebra over $X=\{x_1,x_2,\ldots
\}$ in the free associative algebra generated by
$\{c,x_1,x_2,\ldots\}$. The free magma multiplication $\cdot$
corresponds to the operation $(v,w)\mapsto cvw$ in the free
associative algebra.

\end{remark}

\bigskip
\goodbreak

\begin{lemma}\label{lemleafrestrict}

Let $M_0=\{x_1,x_2,\ldots\},
 M_k=\{\mu_1,\mu_2,\ldots\}$ (for $k\geq 1$), and let
 $T$ be a (planar or abstract) admissibly labeled tree, with
 corresponding parenthesized string $w$.
Let $I\subseteq {\rm Le}(T)$ be a subset of the set of leaves.

We can delete all pairs of brackets (together with its letter)
that have no pair of bracket corresponding to a vertex from $I$ in
between, to obtain a parenthesized string $w\vert I$.

The string $w\vert I$ corresponds to a (not necessarily non-empty)
admissibly labeled tree $T\vert I$. The set of vertices of $T\vert
I$ corresponds to a subset of ${\rm Ve}(T)$, and the arity of $v$
in ${\rm Ve}(T\vert I)$ is less or equal than the arity of the
corresponding vertex in $T$.

\end{lemma}

\begin{proof}

The arity of a vertex $v$ in $T$ corresponds to the number of
irreducible strings in between the pair of brackets corresponding
to $v$. Removing strings as indicated does not increase this
number, and we get an admissibly labeled tree $T\vert I$.

\end{proof}

\begin{defn}

The tree $T\vert I$ is called the (non-reduced) leaf-restriction
of $T$ on $I\subseteq {\rm Le}(T)$.

\end{defn}

\begin{example}
Consider the tree $
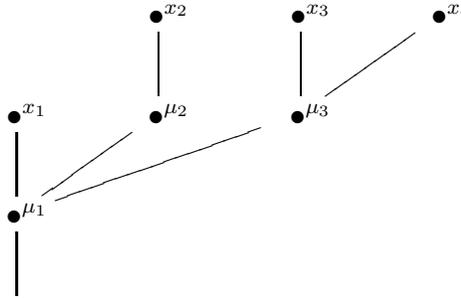
\xymatrix{ & {\ \ \bullet^{x_2}}\ar@{-}[d] & {\
\ \bullet^{x_3}}\ar@{-}[d]& {\ \ \bullet^{x_4}}\ar@{-}[dl]\\
 {\ \ \bullet^{x_1}}\ar@{-}[d] & {\ \ \bullet^{\mu_2}}\ar@{-}[dl] &
 {\ \ \bullet^{\mu_3}}\ar@{-}[dll]& \\
 {\ \ \bullet^{\mu_1}}\ar@{-}[d]& \\ &\\}$

 with parenthesized string
$(\mu_1(x_1)(\mu_2(x_2))(\mu_3(x_3)(x_4)))$.

The leaf-restriction on the first and third leaf is $
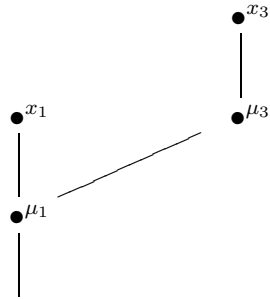
\xymatrix{ &&
{\ \ \bullet^{x_3}}\ar@{-}[d]\\ {\ \ \bullet^{x_1}}\ar@{-}[d] &  &
 {\ \ \bullet^{\mu_3}}\ar@{-}[dll]\\ {\
\ \bullet^{\mu_1}}\ar@{-}[d]\\ \\}$

  with parenthesized string
$(\mu_1(x_1)(\mu_3(x_3)))$.

\end{example}

\begin{remark}

Leaf-restriction is an example for a process induced by the
removal of some vertices of a tree. A similar process is induced
by removing all vertices of arity 1:
 There is a (canonical) map ${\rm
red}: \PTree\to\PRTree$ (and a similar map for abstract trees),
leaving the tree structure intact as much as possible. Admissible
labelings (of all vertices of arity $\neq 1$) are preserved.

\end{remark}

\begin{defn}

For $T$ a (planar or abstract) tree, the tree $ {\rm red}(T)$ is
called the reduction of $T$. Its set of vertices is ${\rm
Ve}(T)-{\rm Ve}^1(T)$, and for any pair $v,v'$ of vertices of
${\rm red}(T)$ there is an oriented path (or, equivalently, a path
not passing the root) from $v$ to $v'$ in  ${\rm red}(T)$ if and
only there is such a path in $T$ (cf.\ also \cite{ligepsp})

\end{defn}

\begin{example}
Consider the tree $\xymatrix{ {\ \ \bullet^{x_1}}\ar@{-}[d] & {\ \
\bullet^{x_2}}\ar@{-}[d] & {\ \ \bullet^{x_3}}\ar@{-}[d]& {\ \
\bullet^{x_4}}\ar@{-}[dl]\\
 {\ \ \bullet^{\mu_2}}\ar@{-}[d] & {\ \ \bullet^{\mu_3}}\ar@{-}[dl] &
 {\ \ \bullet^{\mu_4}}\ar@{-}[dll]& \\
 {\ \ \bullet^{\mu_1}}\ar@{-}[d]& \\ {\ \ \bullet^{\mu_0}}\ar@{-}[d]\\ & \\}$.

Its reduction is $\xymatrix{  &  & {\ \ \bullet^{x_3}}\ar@{-}[d]&
{\ \ \bullet^{x_4}}\ar@{-}[dl]\\
 {\ \ \bullet^{x_1}}\ar@{-}[d] & {\ \ \bullet^{x_2}}\ar@{-}[dl] &
 {\ \ \bullet^{\mu_4}}\ar@{-}[dll]& \\
 {\ \ \bullet^{\mu_1}}\ar@{-}[d]& \\ &\\}$.

\end{example}

\begin{defn}\label{defadmcut}

Let $T$ be a (planar or abstract) tree with root $\rho$, and
$C\subseteq {\rm Ve}(T)$. We call $C$ an admissible cut of $T$, if
for every vertex $v\in C$ all vertices of the full subtree given
by $v$ are also in $C$. The case $C=\emptyset$ is called the empty
cut. The case $C={\rm Ve}(T)$ is called the full cut.

Given such an admissible cut, let $R^C(T)$ be the not necessarily
non-empty tree (with root $\rho$, if $R^C(T)\neq\emp$), obtained
by removing all vertices of $C$ (together with their outgoing
edges).

From $T$ we can remove (the subgraph) $R^C(T)$ to get a (planar or
abstract) forest $C(T)$ with set of vertices $C$.

The pair $(C(T),R^C(T))$ is called result of the cut $C$.

\end{defn}

\begin{remark}

An admissible cut of $T$ can also be defined as a non-empty subset
of the set of (inner) edges of $T$ such that for every vertex
$v\in{\rm Ve}(T)$ on the path to the root there is at most one
edge selected, cf.\ \cite{lick}. This definition leads to the same
pair $(C(T),R^C(T))$ and is in fact equivalent, once we add the
full and empty cut.

\end{remark}

\begin{example}
Let $T$ be the following planar tree: $
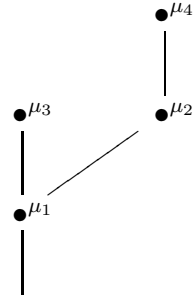
\xymatrix{   & {\ \
\bullet^{\mu_4}}\ar@{-}[d]\\
 {\ \ \bullet^{\mu_3}}\ar@{-}[d] & {\ \ \bullet^{\mu_2}}\ar@{-}[dl] & \\
 {\ \ \bullet^{\mu_1}}\ar@{-}[d]& \\ &\\}$

We are going to indicate (by $\circ$) which vertices are selected.
 Some admissible cuts of $T$ are:

We get $R^C(T)=\emp$ and $C(T)=T$ as the result of the cut
  $
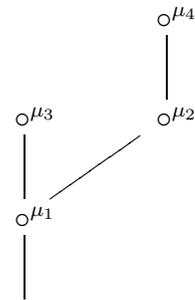
\xymatrix{   & {\ \
\circ^{\mu_4}}\ar@{-}[d]\\
 {\ \ \circ^{\mu_3}}\ar@{-}[d] & {\ \ \circ^{\mu_2}}\ar@{-}[dl] & \\
 {\ \ \circ^{\mu_1}}\ar@{-}[d]& \\ &\\}$

and $R^C(T)=
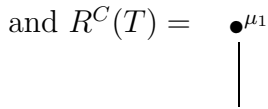
\xymatrix{
 {\ \ \bullet^{\mu_1}}\ar@{-}[d]& \\ &\\}$ and $C(T)=
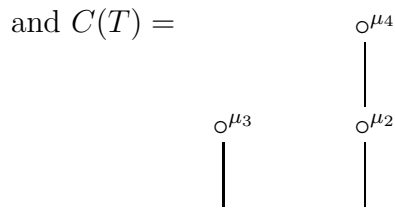
\xymatrix{   & {\ \ \circ^{\mu_4}}\ar@{-}[d]\\
 {\ \ \circ^{\mu_3}}\ar@{-}[d] & {\ \ \circ^{\mu_2}}\ar@{-}[d] \\
& & &\\ } $

as the result of the cut $
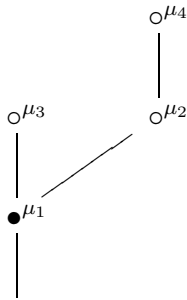
\xymatrix{   & {\ \
\circ^{\mu_4}}\ar@{-}[d]\\
 {\ \ \circ^{\mu_3}}\ar@{-}[d] & {\ \ \circ^{\mu_2}}\ar@{-}[dl] & \\
 {\ \ \bullet^{\mu_1}}\ar@{-}[d]& \\ &\\}
 $

Not an admissible cut is:
 $\xymatrix{   & {\ \
\circ^{\mu_4}}\ar@{-}[d]\\
 {\ \ \circ^{\mu_3}}\ar@{-}[d] & {\ \ \bullet^{\mu_2}}\ar@{-}[dl] & \\
 {\ \ \circ^{\mu_1}}\ar@{-}[d]& \\ &\\}
 $

\end{example}

\begin{remark}

While admissible cuts split trees into branches and rooted trunk,
leaf-restriction onto subsets $I_1, I_2$, where the set of leaves
is a disjoint union of $I_1$ and $I_2$, leads to a different type
of splitting. This splitting is best adapted to trees which allow
different labels only for leaves (and not for internal vertices),
and is described as follows.

\end{remark}

\begin{defn}\label{defleafsplit}
Let $M_0=\{x_1,x_2,\ldots\}$, $M_1=\emptyset$,
 $M_k=\{\bullet\} ($for $k\geq 2)$.
Let $T$ be a planar admissibly labeled tree, especially $T$ is
reduced.

Given a split ${\rm Le}(T)=I_1\uplus I_2$ of the set of leaves of
$T$ into two disjoint subsets $I_1, I_2$, we call the pair
\begin{equation*}
\bigl({\rm red}(T\vert I_1), {\rm red}(T\vert I_2)\bigr)
\end{equation*}
the induced leaf-split of $T$.

Given two planar admissibly labeled trees $T^1, T^2$ with $n_1,
n_2$ leaves, and a planar tree $T$ with $n_1+n_2$ leaves, we say
that $T$ is a shuffle of $T^1$ and $T^2$, if
\begin{equation*}
{\rm red}(T\vert I)=T^1, {\rm red}(T\vert I^c)=T^2
\end{equation*}
for some subset $I\subseteq {\rm Le}(T)$.

If ${\rm red}(T\vert I)=T^1, {\rm red}(T\vert I^c)=T^2$ (for some
$I$) we call $T^1$ a complement of $T^2$ in $T$ (and also $T^2$ a
complement of $T^1$ in $T$).

\end{defn}

\begin{example}

Consider the following planar tree $T$:

$
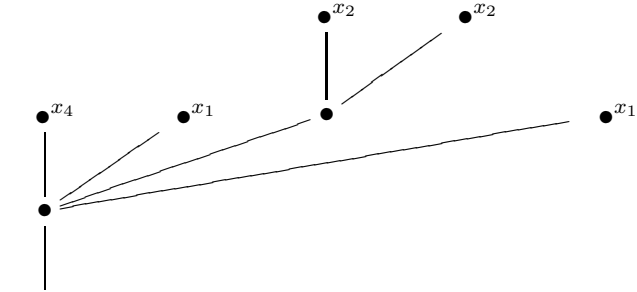
\xymatrix{ &  & {\ \ \bullet^{x_2}}\ar@{-}[d]& {\ \
\bullet^{x_2}}\ar@{-}[dl]\\
 {\ \ \bullet^{x_4}}\ar@{-}[d] & {\ \ \bullet^{x_1}}\ar@{-}[dl] &
 \bullet\ar@{-}[dll]& &{\ \ \bullet^{x_1}}\ar@{-}[dllll]  \\
 \bullet\ar@{-}[d]& \\ &\\}$

It is a shuffle of

$
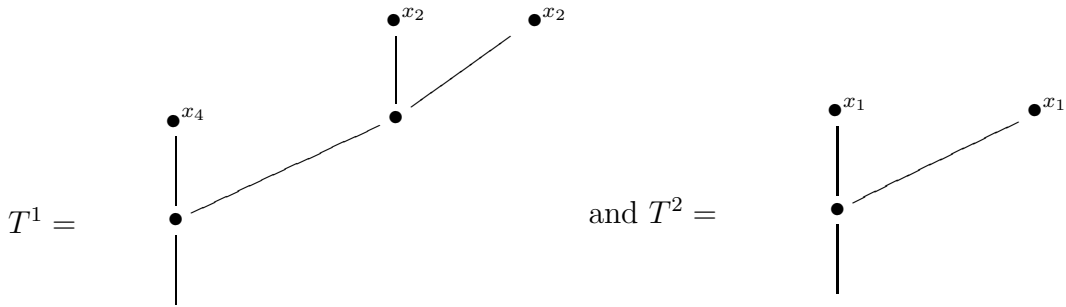
\xymatrix{ & &  & {\ \ \bullet^{x_2}}\ar@{-}[d]& {\ \
\bullet^{x_2}}\ar@{-}[dl]\\ & {\ \ \bullet^{x_4}}\ar@{-}[d] & &
 \bullet\ar@{-}[dll]& \\
 T^1=\!\! &
 \bullet\ar@{-}[d]& \\ & &\\}
\xymatrix{ && \\  &{\ \ \bullet^{x_1}}\ar@{-}[d] &  &
  {\bullet^{x_1}}\ar@{-}[dll]\\
  \text{\! and } T^2=\ &\bullet\ar@{-}[d]\\& \\}$

for the splitting given by the subsets $I=\{1,3,4\}$ and
$I^c=\{2,5\}$ of $\underline 5={\rm Le}(T)$.

The tree $T$ is also a shuffle of the trees

$
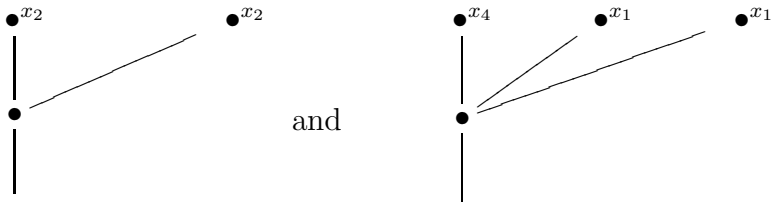
\xymatrix{  &{\ \ \bullet^{x_2}}\ar@{-}[d] &  &
  {\ \ \bullet^{x_2}}\ar@{-}[dll]\\
 \ &\bullet\ar@{-}[d]\\& \\}
\xymatrix{
 &{\ \ \bullet^{x_4}}\ar@{-}[d] & {\ \ \bullet^{x_1}}\ar@{-}[dl] &
 {\ \ \bullet^{x_1}}\ar@{-}[dll] \\
\text{ and } &\bullet\ar@{-}[d]& \\ &\\}$

\end{example}

\end{section}

\begin{section}{Planar binary trees and Stasheff polytopes}

We refine our notation for the various sets of trees, taking into
account the numbers of leaves and internal vertices. Stasheff
polytopes are helpful to link these sets.

We also have to introduce several operations for planar binary
trees.

\begin{defn}

For $n,p\in\N$, let
\begin{equation*}
\PTree_n:=\{ T\in \PTree : \#{\rm Ve}(T)=n\}
\end{equation*}
be the set of planar trees with $n$ vertices, and let
\begin{equation*}
\PTree^p:=\{ T\in \PTree : \#{\rm Le}(T)=p\}
\end{equation*}
be the set of planar trees with $p$ leaves.

Furthermore, let
\begin{equation*}
\PTree_n^p:=\PTree_n\cap\PTree^p.
\end{equation*}

We make the analogous definitions for abstract trees, planar
reduced trees, and the corresponding labeled sets.

Let
\begin{equation*}
\YTree^p:=\YTree\cap\PTree^p.
\end{equation*}

\end{defn}

\begin{remark}

If we consider planar binary trees $T$ and exclude the empty tree,
the number $p$ of leaves is one higher than the number of internal
vertices. Thus $\#{\rm Ve}(T)=2p-1$.

\medskip

Let $(M_k)_{k\in\N}$ be a collection of the form $M_0=\{\circ\}$,
$M_1=\emptyset$, $M_2\neq\emptyset$, $M_k=\emptyset$ (all $k\geq
3$).

Then we consider the labels from $M_2$ as labels of weight 1, and
the labels from $M_0$ as labels of weight 0, and make the
following definition.

\end{remark}

\begin{defn}

Let $\YTree\{M_2\}$ be the set of non-empty admissibly labeled
planar binary trees (i.e.\ the set of planar binary trees that are
labeled at their internal vertices).

We define the degree of an element of  $\YTree\{M_2\}$ to be the
number of internal vertices, and set
\begin{equation*}
\YTree^{(n)}:=\YTree^{n+1}, \text{ all } n\geq 0,
\end{equation*}
and similarly
\begin{equation*}
\YTree^{(n)}\{M_2\}:=\YTree^{n+1}\{M_2\}, \text{ all } n\geq 0,
\end{equation*}
i.e.\ the image of $\YTree^{(n)}$ under $\deg$ is $\{n\}$.

We also define the degree of a planar binary forest to be the sum
of the degrees of all its trees.

We identify $\YTree^{(n)}$ with the set ${\rm
YTree}^{(n)}\!\{\bullet\!\}$, and we also use the notation
\begin{equation*}
\YTree^{\infty}:=\bigcup_{n\geq 0}\YTree^{(n)}.
\end{equation*}

\end{defn}

\begin{remark}

The number $C_n=\#\PRTree^n$ of planar reduced trees with $n$
leaves is called the $n$-th
 super-Catalan number.

The generating series for the super-Catalan numbers is
$\frac{1}{4}(1+x-\sqrt{1-6x+x^2})$ (cf.\ \cite{lisl} A001003).

The first 10 super-Catalan numbers are
\begin{equation*}
1, 1, 3, 11, 45, 197, 903, 4279, 20793, 103049.
 \end{equation*}

\end{remark}

\bigskip

\begin{defn}

Given a planar binary tree $T$ in $\YTree^{\infty}$ or in
$\YTree^{\infty}\{M_2\}$, let $\alpha=\alpha(T)$ denote the first
leaf of $T$ (i.e.\ the leftmost leaf in a drawing which puts all
leaves on one line).

 Let
 $\omega=\omega(T)$ denote the last (i.e.\
rightmost) leaf of $T$.

For $1\leq i < \#{\rm Le}(T)$, the $i$-th internal vertex of $T$
is the highest internal vertex which belongs to both the paths
from the $i$-th and the $(i+1)$-th leaf to the root.

Given a second planar binary tree $S$, we define
\begin{equation*}
 T \backslash S:= T \circ_{\omega(T)} S.
\end{equation*}
\end{defn}

\begin{remark}\label{remcircalpha}

The operation $\backslash$ is called under-operation and was
introduced in \cite{lilrhopftree}. Clearly $\backslash$ is
associative and $S\backslash T\backslash Z$ is well-defined.

The analogous operation $\circ_{\alpha}$ given by $T
\circ_{\alpha}S=T \circ_{\alpha(T)}S$ plays the role of an
associative multiplication in a Hopf algebra defined by C.\
Brouder and A.\ Frabetti, see Section  \ref{sectbf}.

The opposite multiplication $\circ_{\alpha}^{op}$ is the
over-operation $S / T:= T \circ_{\alpha(T)} S$ of
\cite{lilrhopftree}.

Using the mirror-operation $T\mapsto \bar T$, one can express
$S/T$ as $\overline{\bigl(\bar T\bigr) \backslash \bigl(\bar
S\bigr)}$.

The tree $\vert$ consisting of the root serves as a unit for all
these operations, e.g.\ $S\backslash \vert=S=\vert\backslash S$.

\end{remark}

\begin{example}

The (non-labeled) right comb of height $n$ can be expressed as
\begin{equation*}
\underbrace{Y\backslash Y\backslash Y \ldots \backslash Y}_n
\end{equation*}
if $Y$ denotes the planar binary tree
 $\xymatrix{ \\ & {\circ}\ar@{-}[d] & \\ & {\bullet}\ar@{-}[dl]&{
\circ}\ar@{-}[l]\\ \\}$

\end{example}

\begin{defn}
 Given $n\geq 1$ planar binary trees $T^1,\ldots,T^n$, and a
sequence $w$ of $n$ labels from $M_2$,  we define
\begin{equation*}
\vup_w (T^1\ldots T^n)
\end{equation*}

to be the planar binary tree which can be obtained from the right
comb $C$ of height $n$ as follows: We replace the first leaf
$\alpha(C)$ by $T^1$, the second by $T^2$ and so on, leaving the
$(n+1)$-th leaf (i.e.\ $\omega(C)$) unaltered. We just write $\vup
(T^1\ldots T^n)$ if there is no choice of labels.

We denote the
 tree consisting of the root
(which has got the standard label from $M_0=\{\circ\}$) by
$\ladderonopen$ or simply by $\vert$.

We define $\vup(\emptyset)=\vert$.

\end{defn}

\begin{remark}\label{remrightcombpres}

The right comb of height $n$  can then be written as
\begin{equation*}
\vup_w (n):=\vup_w (\underbrace{\vert\ldots\vert}_n).
\end{equation*}

It is easy to see that the smallest set which contains $\vert$ and
is closed under $\vup_w$ operations contains all planar binary
trees.

For every planar binary tree in $\YTree^{\infty}\{M_2\}$ this
right comb presentation is unique. The left comb presentation is
similarly defined. The right (or left) comb presentation induces a
map $\varphi_r$ ($\varphi_l$ respectively) from planar binary
trees to (not necessarily binary) planar forests such that
\begin{equation*}
\varphi_r(\vert)=\emptyset,\ \ \varphi_r\Bigl(\vup_{w_1\ldots w_n}
(T^1\ldots T^n)\Bigr)
=\vee_{w_1}\bigl(\varphi_r(T^1)\bigr)\ldots\vee_{w_n}\bigl(\varphi_r(T^n)\bigr)
\end{equation*}

\end{remark}

\begin{example}

The following admissibly labeled planar binary tree

in $\YTree^{(10)}\{\underline{10}\}$

 $
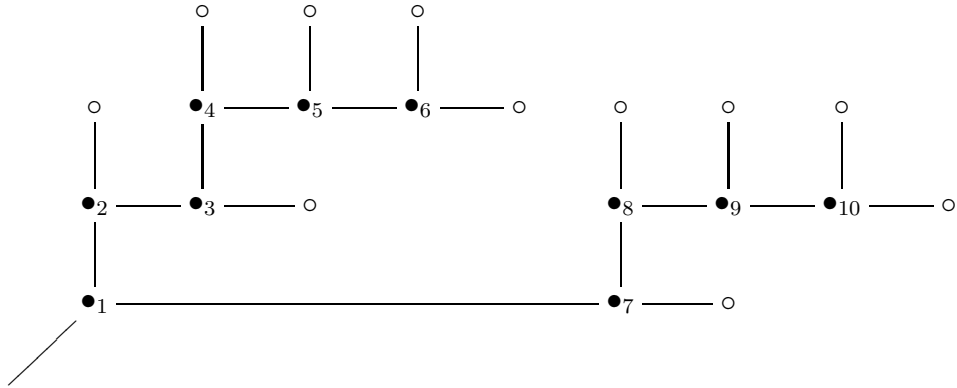
\xymatrix{ & &\circ\ar@{-}[d] &\circ\ar@{-}[d] &\circ\ar@{-}[d] & & & \\
 &\circ\ar@{-}[d] & \ar@{-}[d]\bullet_4 \ar@{-}[r]& \bullet_5\ar@{-}[r] &
\bullet_6 \ar@{-}[r] &\circ & \circ\ar@{-}[d] &
\circ\ar@{-}[d]&\circ\ar@{-}[d]
 \\
 & \bullet_2 \ar@{-}[d]\ar@{-}[r] & \bullet_3 \ar@{-}[r]&\circ & &
 & \ar@{-}[d]\bullet_8 \ar@{-}[r]& \ar@{-}[r]\bullet_9
 \ar@{-}[r] & \bullet_{10}\ar@{-}[r] &\circ \\
 & \bullet_1\ar@{-}[dl]\ar@{-}[rrrrr] & & & & &
  \bullet_7\ar@{-}[r] &\circ & \\
  \\}$

can be written in right comb presentation as
\begin{equation*}
\vup_{1,7}\Bigl(\vup_{2,3}\bigl(\vert,\vup_{4,5,6}(\vert\vert\vert)\bigr)
\ ,\ \vup_{8,9,10}(\vert\vert\vert)\Bigr)
\end{equation*}

Its image under $\varphi_r$ is this forest:

 $
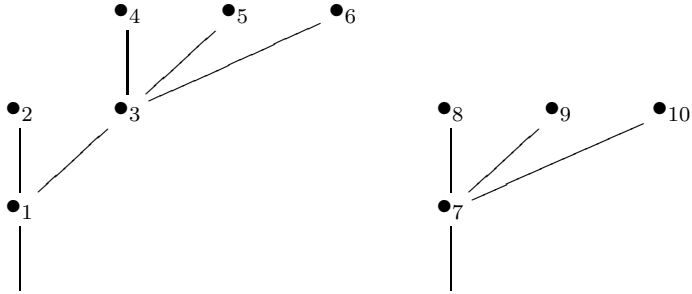
\xymatrix{& & \ar@{-}[d]\bullet_4 & \ar@{-}[dl]\bullet_5 &
 \ar@{-}[dll]\bullet_6
& &  & &
 \\
& \bullet_2 \ar@{-}[d] & \bullet_3 \ar@{-}[dl]& & & \bullet_8
\ar@{-}[d]& \ar@{-}[dl]\bullet_9
  & \ar@{-}[dll]\bullet_{10} \\
 & \bullet_1  \ar@{-}[d]& & & &
  \bullet_7 \ar@{-}[d]& & \\ & & & & &
  \\}$

\end{example}

\begin{lemma}\label{lemtreebij}

Let the degree of a planar forest be given by the total number of
vertices.

Then the map $\varphi_r$ (or $\varphi_l$) provides a bijection
between the set $\YTree^{(n)}\{M_2\}$ of degree $n$
 admissibly labeled planar binary trees  and
degree $n$ planar forests labeled by the set $M_2$.

The number of planar forests, labeled by the set $M_2$, of degree
$n-1$, as well as the number of admissibly ($M_2$-)labeled planar
binary trees of degree $n-1$ is given by
\begin{equation*}
 c_{n}\cdot {(\# M_2)}^{n-1}
 \end{equation*}

where $c_n$ is the $n$-th Catalan number.

The bijection between non-labeled planar trees with $n$ vertices
and planar binary
 trees
 with $n$ leaves occurs as a special case.

\end{lemma}

\begin{proof}
The Lemma is an easy consequence of our construction.
\end{proof}

\begin{remark}\label{remassociahedron}

There is a convex polytope $K_{n+1}$ of dimension $n-1$, $n\geq
1$, with one vertex for each planar binary tree with $n+1$ leaves,
which is called the Stasheff polytope or associahedron in
dimension $n-1$.

More exactly, $K_{n+1}$ is a cell complex in dimension $n-1$ with
the elements of $\YTree^{(n)}=\PRTree^{n+1}_{2n+1}$
 as 0-cells.

The associahedra were created by J.\ Stasheff \cite{listafirst} to
study higher homotopies for associativity. If we consider the
parenthesized strings (of 3 letters) given by the two planar
binary trees with 3 leaves, we can get from one to the other by
shifting a bracket, in other words (cf.\ \cite{limss}, I.1),
applying an associating homotopy $h(x,y,z)$ from $x(yz)$ to
$(xy)z$.

In $K_4$ for example, the five planar binary trees with 4 leaves
have to be arranged in a pentagon

 $
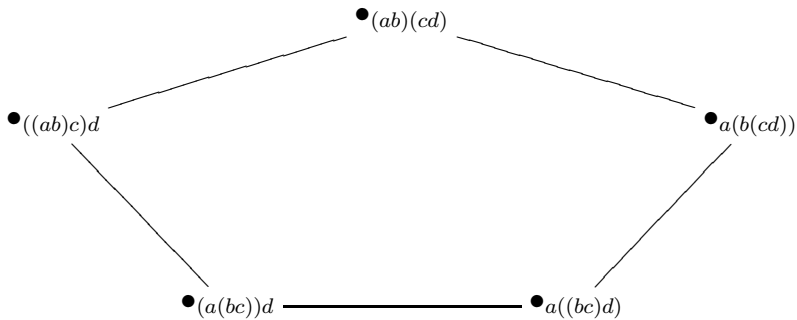
\xymatrix{& & \ar@{-}[dll]\bullet_{(ab)(cd)} \ar@{-}[drr]
 \\
 \bullet_{((ab)c)d} \ar@{-}[ddr] & & & & \bullet_{a(b(cd))} \ar@{-}[ddl]
 \\
 & & & \\
 & \bullet_{(a(bc))d}  \ar@{-}[rr]& &\bullet_{a((bc)d)}
    \\}$

such that each side corresponds to an application of $h(x,y,z)$.
These 5 sides can be labeled by the 5 reduced planar trees (with
$4$ leaves and 2 internal vertices) indicating the associating
homotopy. For example, the edge between $(a(bc))d$ and $((ab)c)d$
corresponds to the tree $
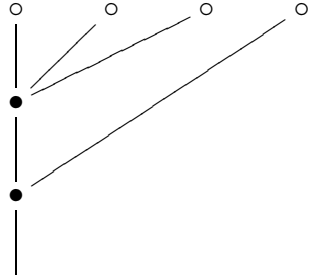
\xymatrix{\\ \circ\ar@{-}[d] &
\circ\ar@{-}[dl] & \circ\ar@{-}[dll] & \circ\ar@{-}[ddlll]\\
{\bullet}\ar@{-}[d]\\ \bullet\ar@{-}[d]\\
\\}$

\medskip

The $(n+1)$-corolla represents the top dimensional cell of the
polytope
 $K_{n+1}$.

By definition of the cell complex $K_{n+1}$ the
 cells of dimension $k$ are in bijection with
the elements of $\PRTree^{n+1}_{2n-k+1}$, for $k=0,\ldots,n-1$.

Thus the super-Catalan numbers count all cells of the Stasheff
polytope, and we can check that $K_4$, for example, has got 11
cells.

\medskip

The polytope $K_2$ is a point which corresponds to the unique
element $\Y$ of $\YTree^{(1)}=\PRTree^{2}_{3}$, and $K_3$ is an
interval.

The facets (i.e.\ codimension one cells) of $K_{n+1}$ are of the
form $K_{r+1}\times K_{s+1}$, $r,s\geq 1$, $r+s=n$, with label
obtained by grafting the $s$-corolla to the $i$-th leaf of the
$r$-corolla, $1\leq i\leq r$. One gets inclusion maps $\circ_i:
K_{r+1}\times K_{s+1}\to K_{r+s+1}$
 (cf.\ \cite{limss}, I.1.6).

Pentagons and squares are the facets of $K_5$.

\bigskip

In fact, the realization of $K_{n+1}$ as a convex polytope was an
open problem at first. Several solutions were given (cf.\
\cite{listaphys}). A simple realization, given in \cite{lilosta},
associates to each planar binary tree $T\in\YTree^{(n)}$ a
coordinate tuple $x(T)$ in (a hyperplane of) $\R^{n}$ as follows:

For $1\leq i < n+1=\#{\rm Le}(T)$, consider the subtree with root
given by the $i$-th internal vertex of $T$, and let $a_i$ be the
number of leaves on the left side, $b_i$ the number of leaves on
the right side (of the subtree's root).
 Then the $i$-th entry of $x(T)$
is given by $a_i b_i$.

For example, the coordinate tuples we get for $K_4$ are:

 $
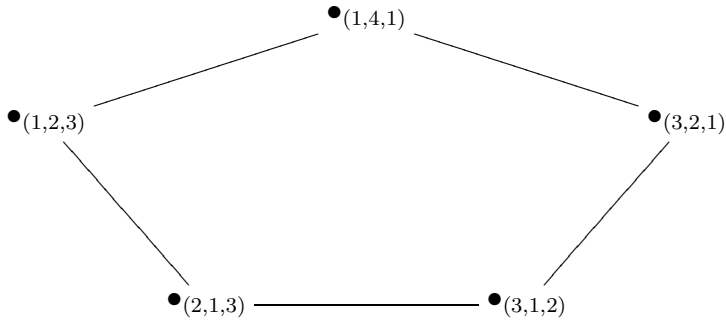
\xymatrix{& & \ar@{-}[dll]\bullet_{(1,4,1)} \ar@{-}[drr]
 \\
 \bullet_{(1,2,3)} \ar@{-}[ddr] & & & & \bullet_{(3,2,1)} \ar@{-}[ddl]
 \\
 & & & \\
 & \bullet_{(2,1,3)}  \ar@{-}[rr]& &\bullet_{(3,1,2)}
    \\}$

It is shown in \cite{lilosta}, Theorem 1.1, that the convex hull
of the points $x(T), T\in\YTree^{(n)},$ is a realization of the
Stasheff polytope of dimension $n-1$.

\medskip

It is possible to give an orientation to all the edges of the
Stasheff polytope, see \cite{lilrb}.
When 0-cells are represented as parenthesized words, arrows are directed such that
they correspond to shifting a bracket
 from left ((xx)x) to right (x(xx)).

 The induced partial
ordering on the set $\YTree^{(n)}$ (also called
Tamari order) and the weak Bruhat order on the Coxeter group $\Sigma_n$
are related by a projection $\Sigma_n\to\YTree^{(n)}$ given in \cite{lilrb}.

\end{remark}

\end{section}

\begin{section}{Free operads and quotients}\label{sectfreeoperad}

The forgetful functor from the category of operads to the category
$\Sigma$-Vect$_K$ has a left adjoint, the free operad functor
$\Gamma$ (cf.\ \cite{limss}, p.82, or \cite{lifrkos}).

Given a $\Sigma$-space $A$ with basis $\alpha_1\in A(n_1)$, $
\ldots$, $\alpha_l\in A(n_l)$ such that $n_i\geq 2$ for all $i$,
we present an explicit construction of $\Gamma(A)$ in terms of
reduced trees. We start with the case where we do not have to
consider a $\Sigma$-action (i.e.\ the case of non-$\Sigma$ operads
or the case of regular operads).

It is possible to use abstract trees or planar trees for the
description of free operads. Both approaches have advantages and
disadvantages
 (cf.\ \cite{limss}, p.83).
A common way to give the definition of a concrete operad is to
present it as a quotient of a free operad modulo an operad ideal.

\begin{defn}\label{defoperadideal}

Let $\underline{\mathcal P}$ (or $\mathcal P$, respectively) be a
non-$\Sigma$ operad (or an operad), and let $\mathcal R(n)$ be a
subspace (invariant subspace) of $\mathcal P(n)$ for each $n$.

Then $\mathcal R=(\mathcal R(n))_{n\in\N}$ is called ideal of
$\underline{\mathcal P}$ (or $\mathcal P$, respectively), if for
$p\circ_i q$ defined in $\underline{\mathcal P}$ (or $\mathcal P$)
it follows that $p\circ_i q \in\mathcal R$ whenever $p\in\mathcal
R$ or $q\in \mathcal R$.

\end{defn}

\begin{remark}

For every ideal $\mathcal R$ of an operad (or non-$\Sigma$ operad)
$\mathcal P$, the quotient $\mathcal P/\mathcal R$ given by
$(\mathcal P/\mathcal R)(n)={\mathcal P}(n)/{\mathcal R}(n)$ is
again a (non-$\Sigma$) operad.

Given a subset $\{r_i: i\in I\}$ of $\mathcal P$, the smallest
ideal containing $\{r_i: i\in I\}$ is also called the operad ideal
generated by the relations $r_i=0$.

\medskip

If we consider the free algebras over free operads $\mathcal P$,
we can get the free algebras over a quotient operad $\mathcal Q$
as quotient algebras. In this way,  the given operad ideals are
determined by multilinear relations $R$ in free algebras over free
operads. One often considers the class of all $\mathcal
Q$-algebras as a subvariety determined by the identities $R$ in
the variety of all $\mathcal P$-algebras, cf.\ \cite{lidre}.

\end{remark}

\begin{defn}
Let a collection $M=(M_k)_{k\geq 2},$ of sets be given, and set
 $M_0:=\{\circ\}, M_1:=\emptyset$.

Define  for $n\geq 2$:
\begin{equation*}
\Gamma(M)(n)=\bigoplus_{{\rm planar\ trees \  } T {\ \rm with \ }
n {\ \rm leaves} } K\ \TMseq
\end{equation*}
and let $\Gamma(M)(0)=0$, $\Gamma(M)(1)=K\cdot\vert$, where
$\vert$ is the tree consisting of the root.

The elements of $\Gamma(M)(n)$ are linear combinations of
admissibly labeled (necessarily reduced) planar trees with $n$
leaves.

Let $T^1$ be such a tree with $n$ leaves, and let $T^2$ be a
further tree having $m$ leaves. Then, for $i=1,\ldots,n$,
\begin{equation*}
T^1 \circ_i T^2
\end{equation*}
is given by the substitution of $T^2$ in $T^1$ at the $i$-th leaf,
see Definition (\ref{defsubstitution}).

By $K$-linear extension, we define $\circ_i$-operations
 \begin{equation*}
\circ_i : {\Gamma(M)}(n)\otimes {\Gamma(M)}(m)\to
{\Gamma(M)}(m+n-1) , \text{ all } n, m\geq 1,\ 1\leq i\leq n.
\end{equation*}

\end{defn}

\begin{lemma}\label{lemfreeoperad}

The sequence $\Gamma(M)$ of vector spaces, together with the unit
$\vert$ and composition maps determined by the given
$\circ_i$-operations, is the free non-$\Sigma$ operad generated by
the collection $M=(M_k)_{k\geq 2}$.
\end{lemma}
\begin{proof}
This is a special case of the free operad construction (given
for example in \cite{limss}); we give only a sketch of the proof.
Free composition products of elements $\alpha=\alpha(x_1.x_2.\ldots.x_k)\in M_k$
($k\geq 2$) can be written as admissibly labeled
leveled trees. We get vector spaces $\overline{\Gamma(M)}(n)$
generated by leveled trees with $n$ leaves.
The sequence $\Gamma(M)(n)$ is obtained by forgetting all levels, this
corresponds exactly to the relations imposed by associativity
of operad-composition.
\end{proof}

\begin{example}

 Let
$M_2$ consist of one generator $\alpha$, and let $M_k=\emptyset\
(k\geq 3)$.

Then all elements of $\Gamma(M)$ are linear combinations of planar
binary trees. For $n\geq 1$, we can identify a basis of
$\Gamma(M)(n)$ with the set of (non-labeled) planar binary trees
with $n$ leaves. Especially, $\dim\Gamma(M)(n)=c_n$.

The tree $\Y$ corresponds to the binary operation $\alpha$, and we
get ternary operations $\alpha\circ_1 \alpha$ and $\alpha\circ_2
\alpha$ as compositions:

 $ 
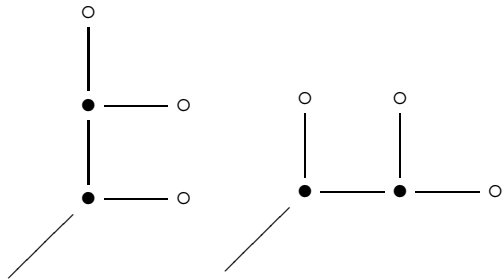
\xymatrix{ &  {\circ}\ar@{-}[d] \\  &
{\bullet}\ar@{-}[d]\ar@{-}[r]  & {\circ} \\
 & {\bullet}\ar@{-}[dl]  &{\circ}\ar@{-}[l]  \\ \\} \
  \xymatrix{ &  & \\ & {\circ}\ar@{-}[d] &
{\circ}\ar@{-}[d] \\
 & {\bullet}\ar@{-}[dl] &{\bullet}\ar@{-}[l]  &{\circ}\ar@{-}[l]  \\ \\} \
$

This is the non-$\Sigma$ operad $\underline{\Mag}$ of
(non-unitary) magma algebras.

\end{example}

\begin{example}\label{exminf}

Let each $M_k$,\ $k\geq 2$, consist of one generator $\vee^k$.

Then $\Gamma(M)(n)$ can be identified with the set of linear
combinations of (non-labeled) planar trees with $n$ leaves. The
$\circ_i$-operations are
 exactly the $\circ_i$-operations defined for the
the associahedron, see (\ref{remassociahedron}),
 given by grafting of trees.

Especially: the generator $\vee^k$ corresponds to the $k$-corolla,
which can be considered as the grafting operation $\vee$
restricted to planar forests consisting of $k$ trees (and their
$K$-linear combinations).

 This free non-$\Sigma$ operad $\underline{\mathcal K}$, with
 $\underline{\mathcal K}^n=\Gamma(M)(n)$ (all $n$)
 is called the
non-$\Sigma$ operad of Stasheff polytopes. The cells of the
associahedron in dimension $n-1$
 form a basis of $\underline{\mathcal K}^{n+1}$.
There is an obvious inclusion
$\underline{\Mag}\to\underline{\mathcal K}$ defined, and we will
denote $\underline{\mathcal K}$ also by $\underline\Minf$.

\end{example}

\begin{remark}

We sketch the construction of the free operad $\Gamma(A)$ for a
given $\Sigma$-space $A$ with basis $\alpha_1\in A(n_1)$, $
\ldots$, $\alpha_l\in A(n_l)$ ($n_i\geq 2$ for all $i$). It is
similar to the construction of the free non-$\Sigma$ operad.

The sum
\begin{equation*}
\bigoplus_{{\rm planar\ trees \  } T {\ \rm with \ } n {\ \rm
leaves} } K\ \TMseq
\end{equation*}
has to be replaced by a sum over trees $T$ with $n$ leaves,
labeled as follows:

The labeling map sends the set of leaves bijectively to the set
$\underline n$, and internal vertices of arity $k$ are labeled by
basis elements of $A(k)$.

The symmetric group $\Sigma_n$ acts by permuting the labels of the
leaves. For each internal vertex, we have to quotient out
relations
 
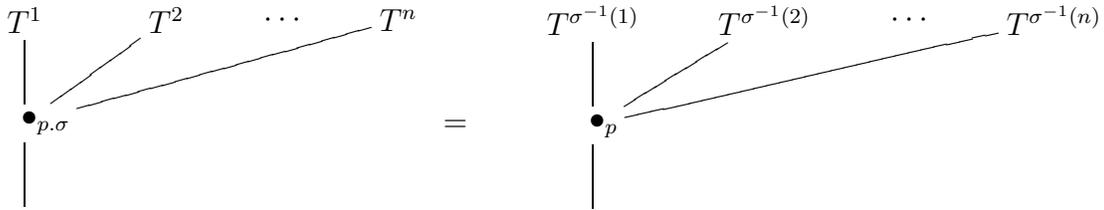
\begin{equation*}
 \xymatrix{
  T^1 \ar@{-}[d]& T^2 \ar@{-}[dl]& \cdots  &\ar@{-}[dlll] T^n \\ {\ \ \ \ \bullet_{p.\sigma}}\ar@{-}[d] & & \\
  &\\}
\ \xymatrix{
  & T^{\sigma^{-1}(1)} \ar@{-}[d]& T^{\sigma^{-1}(2)} \ar@{-}[dl]& \cdots  &\ar@{-}[dlll] T^{\sigma^{-1}(n)} \\
   =& {\ \ \bullet_{p}}\ar@{-}[d] & & \\
  & \\}
\end{equation*}

In fact, this is a colimit-construction. (One can also obtain a
quotient-free description using fixed representatives of abstract
trees, cf.\ \cite{lifrkos}, 3.6.1).

\end{remark}

\begin{example}

Let $A=A(2)$ be the $\Sigma_2$-module given by two generators
$\alpha$, $\beta$ with $\alpha.\tau=\beta$, where $\tau=(1,2)$ is
the transposition. Then
\begin{equation*}
{\phantom \Y}^1\Y^2_{\alpha}={\phantom \Y}^2\Y^1_{\beta}
\end{equation*}
and the free operad generated by $A$ is the symmetrization $\Mag$
of $\underline{\Mag}$, with ${\Mag}(n)=
\underline{\Mag}(n)\otimes_K K\Sigma_n$ (all $n$).

The symmetrization $\Minf$ of $\underline\Minf$ from
(\ref{exminf}) will be called the Stasheff operad, for short.

\medskip

If $A=A(2)$ is the $\Sigma_2$-module given by one generator
$\alpha$ with $\alpha.\tau=\alpha$, then
\begin{equation*}
{\phantom \Y}^1\Y^2_{\alpha}={\phantom \Y}^2\Y^1_{\alpha}
\end{equation*}
and
\begin{equation*}
\Gamma(M)(n)=\bigoplus_{{\rm abstract\ binary\ trees \  } T {\ \rm
with \ } n {\ \rm leaves} } K\cdot T
\end{equation*}
This operad is called the operad $\Cmg$ of commutative magma
algebras. It is not a regular operad.

Similarly, for the free operad $\Gamma$ generated by
$A(n)=K\{\alpha_n\}$\ ($n\geq 2$) with trivial action of
$\Sigma_n$, we have
\begin{equation*}
\Gamma(M)(n)=\bigoplus_{{\rm abstract\ reduced\ trees \  } T {\
\rm with \ } n {\ \rm leaves} } K\cdot T
\end{equation*}
The operad $\Cminf:=\Gamma$ might be called the operad of
commutative tree algebras.

\end{example}

\begin{example}

The non-$\Sigma$ operad $\underline{\As}$ is a quotient of
$\underline{\Mag}$ with respect to the associativity relation
 $ 
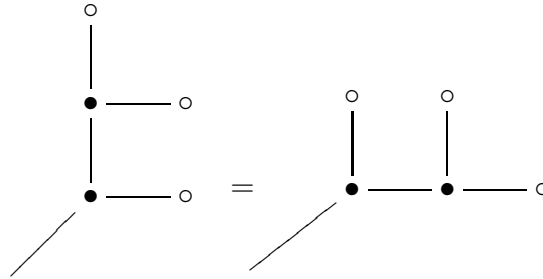
\xymatrix{ &  {\circ}\ar@{-}[d] \\  &
{\bullet}\ar@{-}[d]\ar@{-}[r]  & {\circ} \\
 & {\bullet}\ar@{-}[dl]  &{\circ}\ar@{-}[l]  \\ \\} \
  \xymatrix{ &  & \\ & {\circ}\ar@{-}[d] &
{\circ}\ar@{-}[d] \\
 \ =\ & {\bullet}\ar@{-}[dl] &{\bullet}\ar@{-}[l]  &{\circ}\ar@{-}[l]  \\ \\} \
$

given by $\alpha\circ_1\alpha=\alpha\circ_2\alpha$.

For each $n$ we get an identification of all operations with $n$
arguments in the quotient.

Consequently, the operad $\As$ (with $n!$ operations of $n$
arguments) is the analogous quotient of $\Mag$.

\end{example}

\begin{remark}

Whenever there is a presentation by binary generators which are
subject to quadratic relations, one speaks of a binary quadratic
operad. More generally a quadratic operad is generated by (not necessarily
binary) generators subject to relations that can be written as
linear combinations of admissibly labeled trees with two internal vertices.

Then (using orthogonal relations) there is a quadratic dual operad $\mathcal P^{!}$,
and there is a concept of Koszul duality, see \cite{ligk}.

We are not going to use the notion of Koszul operads,
 but it is an immediate observation that free operads are also Koszul operads.

\end{remark}

\begin{defn}\label{defdend} (see \cite{lilodia})

A dendriform algebra over $K$ is a $K$-vector space $A$ together
with two binary operations $\prec, \succ: A\otimes A\to A$ called
left and right, such that, if $x*y:=x\succ y + x\prec y$ (for
$x,y\in A)$:

\begin{equation*}
\begin{split}
(x\prec y)\prec z &= x\prec(y*z)\\
 (x\succ y)\prec z &= x\succ
(y\prec z)\\
 (x*y)\succ z &= x\succ (y\succ z)
\end{split}
\end{equation*}

The non-$\Sigma$ operad $\underline{\Dend}$ of dendriform algebras
is defined as the quotient of the free operad on two generators
$\prec, \succ$ with respect to the quadratic relations  induced by
the three equations above.

The (regular) operad $\Dend$ is the symmetrization of
$\underline{\Dend}$.

\end{defn}

\begin{remark}

It is easy to check that the sum $*$ of the operations $\prec$ and
$\succ$ is an associative operation (one sums up the lefthandsides
and righthandsides of the three equations and gets
$(x*y)*z=x*(y*z)$).

Free dendriform algebras will be described in Section
\ref{sectdend}.

A first example of dendriform algebras can be obtained  as
follows:

We equip the standard tensor coalgebra with the shuffle
multiplication $\sh$.

Then, for every pair of words $u=w_1.w_2\ldots w_p$, $v=
w_{p+1}\ldots w_r$, their product $u\sh v=w_1.w_2\ldots w_p \sh
w_{p+1}\ldots w_r$ is a finite sum, which can be 'splitted' in two
parts:
 One defines
$\prec$ by all shuffles with first letter $w_1$, and $\succ$ by
all shuffles with first letter $w_{p+1}$, see \cite{lilodia}.

\end{remark}

\begin{remark}\label{remprop}

While operads describe collections of operations with many inputs
but only one output,  the concept of PROPs, cf.\ \cite{limlbull}, handles
structures with many inputs and many outputs.

Since we are going to mention them sometimes in short remarks,
we sketch how PROPs arise as generalizations
of operads.

A PROP (compare \cite{limapropped}) is a sequence $\mathcal P$ of
$(\Sigma_m,\Sigma_n)$ bimodules together with two types of
compositions, horizontal
\begin{equation*}
{\mathcal P}(m_1,n_1)\otimes\ldots\otimes {\mathcal P}(m_s,n_s)\to
{\mathcal P}(m_1+\ldots+m_s,n_1+\ldots+n_s),
\end{equation*}
for all $m_1,\ldots,m_s,n_1,\ldots,n_s\in \N^*$, and vertical
\begin{equation*}
{\mathcal P}(m,n)\otimes{\mathcal P}(n,k)\to {\mathcal P}(m,k),
\end{equation*}
for all $m,n,k \in \N^*$.

 A PROP is supposed to satisfy axioms
which could be read off from the example of the endomorphism PROP
given by:
\begin{equation*}
{\mathcal End}_V(m,n)={\rm Hom}(V^{\otimes m},V^{\otimes n})
\end{equation*}
with horizontal composition given by the tensor product of linear maps,
and vertical composition by the ordinary composition of linear
maps.

Similarly to algebras over operads, there are defined (al)gebras
over PROPs. Many techniques for operads also hold for PROPs,
e.g.\ Vallette (see \cite{lival}) makes use of Koszul duality
for PROPs. In general, there is no free (al)gebra over a PROP
defined, though.

We already noted in Remark (\ref{remmonoidaloperad}) that operads
are also defined for monoidal categories different from
$(\mbox{Vect}_K, \otimes)$ (cf.\ \cite{limss}, II.1). This also
holds for the concepts of PROPs.

\end{remark}

\end{section}

\end{chapter}

\begin{chapter}{${\mathcal P}$-Hopf algebras}

\begin{section}{Hopf algebra theory}\label{secthopftheory}

Given two categories ${\bf C}, {\bf D}$, the study of
representable functors from ${\bf C}$ to ${\bf D}$ leads to the
study of co-${\bf D}$ objects in ${\bf C}$, see \cite{libh},
\cite{liabe}\S A5, as their representing objects. Here we assume
that the category ${\bf C}$ has finite categorical coproducts $A
\sqcup_{\bf C}\ldots \sqcup_{\bf C}A$.

 When
${\bf D}$ is the category of (semi-)groups, one speaks of co-group
objects (or co-semigroup objects) $A$ in ${\bf C}$. Being the
representing object, $A$ is provided with a morphism $\Delta: A
\to A \sqcup_{\bf C} A$. The associativity of (semi-)group
multiplication translates to coassociativity of $\Delta$, i.e.\
the following diagram commutes:

\begin{equation*}
\xymatrix{
          A \ \ar[d]_{\Delta}
          \ar[rr]^{\Delta}
          & & A\sqcup_{\bf C} A
           \ar[d]^{\Delta\sqcup_{\bf C}\id}\\
             A\otimes A \ar[rr]_{\id\sqcup_{\bf C}\Delta\ \ } & & A\sqcup_{\bf C} A\sqcup_{\bf C} A\\
            }
\end{equation*}

If ${\bf C}$ is the category of $\Com$-algebras, a representable
group-valued functor is known as an affine group scheme, and a
co-group (co-semigroup) object is known as a commutative Hopf
algebra (or bialgebra).

In the category of $\As$-algebras, the categorical coproduct is
given by the free product of $\As$-algebras.
 In Hopf algebra
theory, we consider the category of  $\As$-algebras together with
the tensor product $\otimes$ replacing the categorical coproduct.

Consequently, coassociativity of a $K$-linear map $\Delta: V \to V
\otimes V$, for any vector space $V$, means that the diagram
\begin{equation*}
\xymatrix{
          V \ \ar[d]_{\Delta}
          \ar[rr]^{\Delta}
          & & V\otimes V
           \ar[d]^{\Delta\otimes\id}\\
             V\otimes V \ar[rr]_{\id\otimes \Delta\ \ }& & V\otimes V\otimes V\\
            }
\end{equation*}
commutes.

A coalgebra $(C,{\Delta})$ in Hopf algebra theory just means a
vector space $C$
 together with a coassociative $K$-linear map ${\Delta}: C\to
C\otimes C$.

 This agrees with the definition of an $\As$-coalgebra
in Definition (\ref{defpcoalg}). For each $n\geq 0$ there is (up
to a constant factor) only one $n$-ary co-operation, called
${\Delta}^n$. The map ${\Delta}^n: C\to C^{\otimes n}$ is given by
\begin{equation*}
\begin{split}
&{\Delta}^0=\id, {\Delta}^1={\Delta}\\
&{\Delta}^n=({\Delta}\otimes
\id\otimes\ldots\otimes\id)\circ{\Delta}^{n-1}\\
\end{split}
\end{equation*}

Usually one looks at counitary coalgebras, dually defined to
unitary algebras:

The counit is a $K$-linear map $\varepsilon: C\to K$ with
commutative diagram
\begin{equation*}
\xymatrix{
           K\otimes C\ar[dr]_=
            & \ar[l]_{\varepsilon\otimes\id} C\otimes C
            \ar[r]^{\id\otimes\varepsilon}& C\otimes K\ar[dl]^=
           \\
            & C\ar[u]^{\Delta}
           \\
            }
\end{equation*}

If $g\in C$ is an element such that $\Delta(g)=g\otimes g$ and
$\varepsilon(g)=1$, then $g$ is said to be group-like.

For $g,h$ group-like elements, one defines (cf.\ \cite{limohopf},
p.4)
\begin{equation*}
 \Prim_{g,h}(C):=\{c\in C: \Delta(c)=c\otimes g+h\otimes c\}
\end{equation*}
and calls its elements $(g,h)$-primitives.

\medskip

Let $C=\cup_{i\in \N}C_i$ with finite-dimensional $C_i\subseteq
C_{i+1}$ such that
\begin{equation*}
\Delta(C_n)\subseteq \sum_{i=0}^{n} C_i\otimes C_{n-i}  \text{
(all } n).
\end{equation*}
Then $C$ is called a filtered coalgebra, cf.\cite{liabe}, p.91.

Setting $C^{(i)}=C_i/C_{i-1}$, one obtains a coalgebra
$\oplus_{i\in \N}C^{(i)}$, which is called the graded coalgebra
associated to $C=\cup_{i\in \N}C_i$.

A coalgebra is called graded, if $C=\oplus_{i\in \N}C^{(i)}$ and
\begin{equation*}
\Delta(C^{(n)})\subseteq \sum_{i=0}^{n} C^{(i)}\otimes C^{(n-i)}
\text{ (all } n).
\end{equation*}
 If furthermore there is only one group-like element $g$ in $C$,
$C^{(0)}=K g\cong K$, and $C^{(1)}=\Prim_{g,g}(C)$ then $C$ is
said to be strictly graded.

A unitary associative algebra $A$ together with the structure of a
(counitary) coalgebra on $A$ is called a bialgebra, if $\Delta$
and $\varepsilon$ are algebra homomorphisms. The unit $1=1_A$ is
group-like, and the elements of $\Prim_{1,1}(A)$ are called the
primitive elements of the bialgebra $A$.

Considering bialgebras as co-semigroup objects with categorical
coproduct replaced by $\otimes$, Hopf algebras
$A=(A,\mu,\Delta,\eta,\varepsilon,\sigma)$ are $\otimes$-cogroup
objects (in the category of unitary associative algebras). The
$\otimes$-coinverse map, called antipode, is a $K$-linear map
$\sigma:A\to A$ such that
\begin{equation*}
\mu\circ(\sigma\otimes\id)\circ\Delta=\eta\circ\varepsilon=\mu\circ(\id\otimes\sigma)\circ\Delta
\end{equation*}
(where $\eta:K\to A, 1\mapsto 1_A,$ and $\mu:A\otimes A\to A$
describe the unitary associative algebra structure).

Necessarily, $\sigma(1_A)=1_A$, and $\sigma(a)=-a$ for every
primitive element $a\in A$.

If $A_0=A^{(0)}\cong K$ in the filtered case, with $\varepsilon:
A\to K$ the augmentation map, an antipode can always be
constructed recursively by
\begin{equation*}
\sigma(a)=-a - {\sum\ }'\sigma\bigl(a_{(1)}\bigr)a_{(2)}
\end{equation*}
if $\sum' a_{(1)}\otimes a_{(2)}$ denotes
$\Delta'(a):=\Delta(a)-a\otimes 1 - 1\otimes a$.

Hence in the definition of filtered Hopf algebras it is not
necessary to require an antipode.

\medskip

In general, given a bialgebra $A$ one can construct its Hopf
envelope (see \cite{limanin}, cf.\ also \cite{lighquant}). Or one
can consider the completion $\hat A$ with respect to the topology
induced by the augmentation ideal $I=\ker\varepsilon$, i.e.\
$\widehat A=\underleftarrow{\lim} (A/I^n)$, which is a complete
Hopf algebra (cf.\ \cite{lih}).

For example, let $A=K\langle X_1, X_2, X_3\rangle$ be the free
unitary $\As$-algebra generated by variables $X_1, X_2, X_3$.
Together with the algebra homomorphism $\Delta$ given by
\begin{equation*}
\begin{split}
&\Delta(X_1)=X_1\otimes X_1,\ \  \Delta(X_2)=X_2\otimes X_2,\\
&\Delta(X_3)=X_1\otimes X_3+X_3\otimes X_2,\\
\end{split}
\end{equation*}
 and with the counit
defined by $\varepsilon(X_1)=1=\varepsilon(X_2)$,
$\varepsilon(X_3)=0$, $A$ is a bialgebra.

Then the Hopf envelope is the algebra $K\langle X^{(i)}_1,
X^{(i)}_2, X^{(i)}_3: i\geq 1\rangle/$ relations which assure that
the anti-algebra homomorphism given by
$\sigma(X^{(i)}_j)=X^{(i+1)}_j$ is an antipode.

The completion $\widehat A$ of $A$ is the power series algebra
$\widehat F_{\As}(V_{\{x_1, x_2,x_3\}})$,
 where $x_3$ corresponds to $X_3$, but $x_1$ and $x_2$
correspond to $X_1-1$ and $X_2-1$. The counit sends $x_i\mapsto
0$. The comultiplication $\Delta: \widehat A\to \widehat
A\hat\otimes\widehat A$ is given by
\begin{equation*}
\begin{split}
&\Delta(x_i)=x_i\otimes 1+1\otimes x_i + x_i\otimes x_i,
\text{ for }i=1,2\\ &\Delta(x_3)=x_3\otimes 1+ 1\otimes x_3+
x_1\otimes x_3+x_3\otimes x_2.\\
\end{split}
\end{equation*}
 One easily constructs the antipode. It is the anti-algebra
homomorphism given by
\begin{equation*}
\begin{split}
x_i&\mapsto \sum_{k=0}^{\infty}(-1)^k x_i^k  \text{ (for }
i=1,2)\\ x_3& \mapsto\sum_{k=0}^{\infty}(-1)^{k+1}\Bigl(
\sum_{j=0}^k x_1^jx_3 x_2^{k-j}\Bigr).\\
\end{split}
\end{equation*}

\end{section}

\begin{section}{Operads equipped with unit actions}

Before we can generalize the notion of "usual" bialgebra and Hopf
algebra structures on $\Com$-algebras or $\As$-algebras to the
setting of ${\mathcal P}$-algebras we need to require that the
tensor product of ${\mathcal P}$-algebras is provided with the
structure of a ${\mathcal P}$-algebra.

There are several approaches possible (cf.\cite{limoe},
\cite{lilosci}, \cite{lipseudo}). Here we follow \cite{lilosci}.

\begin{defn}\label{defcanunit}

Let ${\mathcal P}$ be an operad, ${\mathcal P}(0)=0, {\mathcal
P}(1)=K=K\id$.

Let a 0-ary element $\eta$ be adjoined to the $\Sigma$-space
${\mathcal P}$ by ${\mathcal P}'(i):=\begin{cases} {\mathcal P}(i)
&: i\geq 1\\ K \eta &: i=0.\\
\end{cases}$

A unit action on ${\mathcal P}$ is a partial extension of the
operad composition onto ${\mathcal P}'$ in the following sense:

We ask that
 composition maps $\mu_{n;m_1,\ldots,m_n}$ (fulfilling the
associativity, unitary, and invariance conditions)
 are defined on
\begin{equation*}
{\mathcal P}'(n)\otimes {\mathcal P}'(m_1)\otimes \ldots \otimes
{\mathcal P}'(m_n)\to {\mathcal P}'(m) ,\ m:=m_1+\ldots+m_n,
\end{equation*}
for all $m_j\geq 0$ ($j=1,\ldots,n$), for $n\geq 2,\ m>0$ (or
$n\leq 1, m\geq 0$).

\medskip

Especially defined are $K$-linear maps
\begin{equation*}
\underline{\ \ }\circ_{n,i}(\eta): {\mathcal P}'(n)\to {\mathcal
P}'(n-1), \ n\geq 2,\ 1\leq i\leq n.
\end{equation*}

\end{defn}

\begin{remark}\label{remcanunit}

Let a ${\mathcal P}$-algebra $\overline{A}$, with structure maps
$\gamma(n)$, be given. Let $p\in {\mathcal P}(n)$. We denote
$\gamma(0)(\eta)\in A:=K\oplus \overline{A}$ by $1$ (compare
Definition (\ref{defoperadalg})).

Thus the unit action gives sense to expressions
\begin{equation*}
p(a_{1},\ldots,a_n)=\gamma(n)(p)(a_1\otimes\ldots\otimes a_n)
\end{equation*}
in the case where some of the $a_i$ are 1.

For $p\in{\mathcal P}(n), n\geq 2,$
\begin{equation*}
p(1,\ldots,1)
\end{equation*}
is not defined. In Definition (\ref{defcanunit}), the compositions
of $p\otimes q_1\otimes\ldots\otimes q_n$ might more generally be
only required to exist if $p$ is in a set $M$ of fixed generators
of ${\mathcal P}.$

In Section (\ref{sectalgcoalg}) we considered operads and operad
algebras in the non-unitary case.

 For example, the operad $\Com$
leads to commutative associative algebras without unit. Since
there exists (up to a scalar factor) only one operation $\cdot_n$
in ${\mathcal P}(n)$ for each $n$, the canonical unit action is
 induced by $\cdot_n\circ_{n,i}(\eta)=\cdot_{n-1}$, $n\geq 2,\
1\leq i\leq n.$ This just means that in any $\Com$-algebra $A$,
occurrences of $1_A$ in products can be ignored.

Similarly, for regular operads that are generated by (single)
binary operations, like $\As$ and $\Mag$, there is a canonical
unit action induced by $1\cdot a=a\cdot 1=a$.

In these cases, the notion of a unitary ${\mathcal P}$-algebra
(defined as in definition (\ref{defoperadalg}), with ${\mathcal
P}'$ in place of ${\mathcal P}$) makes sense.

\end{remark}

\begin{defn}\label{defrespect}

If  there are operations
$\star_n\in\mathcal P(n)$, all $n\geq 2$, fulfilling
\begin{equation*}
\begin{split}
\star_n \circ_{n,i}(\eta) &=\star_{n-1}, \text{ all } i\\
 \star_2\circ_{2,1}(\eta) &=\star_2\circ_{2,2}(\eta)=\id.\\
\end{split}
\end{equation*}
then we say that the unit action respects the operations
$\star_n$, and we define
\begin{equation*}
\mu_{n;0,\ldots,0}(\star_n\otimes\eta\otimes\ldots\otimes\eta)=\eta.
\end{equation*}

\end{defn}

\begin{defn}\label{defcompatible}

The unit action of (\ref{defrespect}) is said to be compatible
with the relations of ${\mathcal P}$, if the relations still hold
on $A=K\ 1 \oplus\overline{A}$, for every (non-unitary)
 ${\mathcal P}$-algebra $\overline{A}$, whenever the expressions are defined.
We call $A=K\ 1 \oplus\overline{A}$ the unitary ${\mathcal
P}$-algebra associated to $\overline{A}$,  and we call the
projection $K\ 1 \oplus\overline{A}\to K$ with kernel
$\overline{A}$ the augmentation of $A$.

\end{defn}

\begin{proposition}\label{propunitpower}

Let $\mathcal R$ be an operad ideal of $\mathcal P=\As$ or
$\mathcal P=\Mag$, given by relations $r_i, i\in I$.
 If the corresponding quotient operad
${\mathcal P}/(r_i: i\in I)$ is equipped with a compatible unit
action, then also ${\mathcal P}/(r_i: i\in I)^n$ is equipped with
a compatible unit action, for all $n\geq 1$.

 Here,
$(r_i: i\in I)^n$ is the $n$-power of the ideal $(r_i: i\in I)$,
generated by products of $r_i$ ($n$ factors).

\end{proposition}
\begin{proof}

We consider free unitary ${\mathcal P}$-algebras $F^1V=K 1\oplus
F(V)$, and we factor out the relations induced by $(r_i: i\in
I)^n$. Since the unit action for ${\mathcal P}$ is compatible with
$(r_i: i\in I)\supseteq (r_i: i\in I)^n$, it is also compatible
with $(r_i: i\in I)^n$.
\end{proof}

\begin{defn}\label{defcoherent}

The unit action of (\ref{defrespect}) is called coherent, if for
all ${\mathcal P}$-algebras $\overline{A},\overline{B}$ it holds
that $(\overline{A}\otimes K 1)\oplus (K 1\otimes
\overline{B})\oplus (\overline{A}\otimes\overline{B})$
  is again a ${\mathcal P}$-algebra
 with, for all $p$ in some set of generators $M$ of ${\mathcal P}$:
\begin{equation*}
\begin{split}
&p(a_1\otimes b_1,a_2\otimes b_2,\ldots,a_n\otimes b_n):=
\star_n(a_1,a_2,\ldots,a_n)\otimes p(b_1,b_2,\ldots,b_n) \\
&\text{ for } a_i\in A, b_i\in B, \text{ in case that at least one
} b_j\in \overline B,\\
 \text{and } &p(a_1\otimes 1,a_2\otimes
1,\ldots,a_n\otimes 1):= p(a_1,a_2,\ldots,a_n)\otimes 1 \\
&\text{ (if
the righthandside is defined)}.
\\
\end{split}
\end{equation*}

\end{defn}
\begin{remark}\label{remtensorpower}

The definition given here generalizes the definition in
\cite{lilosci}, \S 3.2., where ${\mathcal P}$ is assumed to be a
binary quadratic operad (and where an associative operation $*$
plays the role of the operations $\star_n$).

The definition (\ref{defcoherent}) implies that, given a coherent
unit action, the tensor powers  $A^{\otimes n}$ are defined as
${\mathcal P}$-algebras, because

\begin{equation*}
\begin{split}
&p(a_1\otimes (b_1\otimes c_1),a_2\otimes(b_2\otimes
c_2),\ldots,a_n\otimes (b_n\otimes c_n))\\
 &=
\star_n(a_1,a_2,\ldots,a_n)\otimes \star_n(b_1,b_2,\ldots,b_n)
\otimes p(c_1,c_2,\ldots,c_n)\\ &=p((a_1\otimes b_1)\otimes
c_1,(a_2\otimes b_2)\otimes c_2,\ldots,(a_n\otimes b_n)\otimes
c_n)\\
\end{split}
\end{equation*}

The unit is $1\otimes 1\otimes\ldots\otimes 1$.

If $A=K\oplus F_{\mathcal P}(V)=K\oplus\bigoplus_{n=1}^{\infty} {\mathcal
P}(n)\otimes_{\Sigma_n}V_X^{\otimes n}$ is a free $\mathcal P$-algebra,
then
$A^{\otimes n}$ is as a ${\mathcal P}$-algebra (non-freely)
generated by the elements
\begin{equation*}
x_i^{\otimes n,p}:=\underbrace{1\otimes\ldots\otimes 1}_{p-1}\otimes
x_i\otimes\underbrace{1\otimes\ldots\otimes 1}_{n-p}\ \
(\text{for }1\leq p\leq n,\text{ all } x_i\in X).
\end{equation*}

\end{remark}

\begin{example}
\begin{itemize}
\item[1)]
 The operads $\Com, \As$ are equipped with canonical unit
actions, which are clearly coherent.
\item[2)]
In the same way, the operad $\Mag$ is equipped with a canonical
unit action, and we are going to define a similar unit action on
$\Minf$.
\item[3)]
The unit actions on $\Mag$ and $\Minf$ induce unit actions
on $\Cmg$ and $\Cminf$.
\item[4)]
The operad $\Dend$ is also equipped with a coherent unit action,
see Section (\ref{sectdend}).
\item[5)]
If $\mathcal R$ is an operad ideal of $\mathcal P=\As$ or
$\mathcal P=\Mag$, given by relations $r_i, i\in I$ as in
Proposition (\ref{propunitpower}), it is easy to see that
${\mathcal P}/[(r_i: i\in I),(r_i: i\in I)]$ is equipped with a
coherent unit action, too. (The commutator ideal may also be
replaced by the 'associator ideal' $(R,R,R)$ or combinations of
these ideals.)
\end{itemize}
\end{example}

\begin{remark}

The structure of a $\mathcal P$-algebra on $F_{\mathcal
P}(V)\otimes F_{\mathcal P}(V)$ leads to an operad morphism
${\mathcal P}\to {\mathcal P}\otimes {\mathcal P}$, which is
coassociative. Vice versa such an operad morphism determines the
$\mathcal P$-algebra structure on the tensor product
 (as noted by \cite{lifocopr}, \cite{livdl}).
 There is a
  corresponding notion of Hopf algebras over Hopf operads, cf.\ \cite{limoe}
  and \cite{ligejo},
 which does not cover the case of $\Dend$-algebras.

\end{remark}

\end{section}

\begin{section}{The definition of ${\mathcal P}$-Hopf algebras}

Since Hopf algebras combine operations and cooperations, there is
no operad whose algebras or coalgebras are Hopf algebras. To
describe them, and also generalizations with not necessarily
associative operations and not necessarily associative
cooperations, one would use PROPs, see Remark (\ref{remprop}).

Here we are only interested in the case, where the set of
cooperations is generated by one coassociative cooperation. Thus
we do not need the generality of PROPs and stay close to operad theory.

\begin{defn}\label{defpbialg}

Let ${\mathcal P}$ be equipped with a coherent unit action and let
$A=K\ 1 \oplus\overline{A}$ be a unitary ${\mathcal P}$-algebra.
Let $\Delta: A\to A\otimes A$ be a $K$-linear coassociative map
(called comultiplication map), such that $\Delta(1)=1\otimes 1$
and $\Delta(a)-a\otimes 1 - 1\otimes a \in {\overline A}\otimes
{\overline A}$ for all $a\in {\overline A}$.

Then $A$ together with $\Delta$ is called an (augmented)
${\mathcal P}$-bialgebra, if $\Delta$ is a morphism of unitary
${\mathcal P}$-algebras, i.e.\ if $\Delta\circ\mu_A=\mu_{A\otimes
A}\circ(\Delta\otimes\ldots\otimes\Delta),$ for all $\mu\in
{\mathcal P}$.

\end{defn}

\begin{remark}

In the general case,
the condition that $\Delta$ is a morphism of ${\mathcal
P}$-algebras does not imply
that the operations $\mu_A$ (for $\mu\in {\mathcal P}$) are
morphisms of ${\mathcal P}$-coalgebras.
In fact, we do not have provided $A\otimes A$ with
the structure of a ${\mathcal P}$-coalgebra in general.

\end{remark}

\begin{defn}
Elements $a$ of an (augmented) ${\mathcal P}$-bialgebra $A$ are
called primitive, if
\begin{equation*}
\Delta(a)-a\otimes 1 - 1\otimes a=0.
\end{equation*}

\end{defn}

\begin{defn}\label{defphopf}

Let $A=\cup_{n\in \N} A_n$ be an (augmented) ${\mathcal
P}$-bialgebra, with comultiplication $\Delta$. We assume that each
$A_i$ is a finite dimensional vector space, and that $A_i\subseteq
A_{i+1}$ (all $i$).

Then $A$ is called a filtered ${\mathcal P}$-Hopf algebra if
\begin{equation*}
\Delta(A_n)\subseteq \sum_{i=0}^{n} A_i\otimes A_{n-i}  \text{
(all } n).
\end{equation*}

Let $\Delta'(a):=\Delta(a)-a\otimes 1 - 1\otimes a.$

A filtered  ${\mathcal P}$-Hopf algebra $A=\cup_{n\in \N} A_n$ is
called connected, if $A_0=K\cdot 1$, and (for $n\geq 1$)
\begin{equation*}
A_n=\bigl\{a\in A: \Delta'(a)\in \sum_{i=1}^{n-1} A_i\otimes
A_{n-i}\bigr\}.
\end{equation*}

We call a filtered  ${\mathcal P}$-Hopf algebra $A=\cup_{n\in \N}
A_n$ a strictly graded ${\mathcal P}$-Hopf algebra, if
\begin{equation*}
\Delta'(A^{(n)})\subseteq \sum_{i=1}^{n-1} A^{(i)}\otimes
A^{(n-i)} \text{ (all } n\geq 1)
\end{equation*}
for $A^{(i)}\subseteq A_i$ such that $A=\oplus_{n\in \N} A^{(n)}$,
with $A^{(0)}=K$.

\medskip

Homomorphisms of ${\mathcal P}$-bialgebras and ${\mathcal P}$-Hopf
algebras  $(A,\Delta_A), (B, \Delta_B)$ are unitary ${\mathcal
P}$-algebra homomorphisms $\varphi:A\to B$ with $\varphi(\overline
A)\subseteq \overline B$ and
\begin{equation*}
\Delta_B\circ\varphi=(\varphi\otimes\varphi)\circ\Delta_A.
\end{equation*}

\end{defn}

\begin{remark}\label{remantipode}

In (\ref{defphopf}) we used the term filtered ${\mathcal P}$-Hopf
algebra instead of filtered ${\mathcal P}$-bialgebra.

If  $\mu$ is a binary operation of ${\mathcal P}$, $A$ is a
connected filtered ${\mathcal P}$-Hopf algebra, and $\mu\bigl(
\cup_{k=0}^n A_k\otimes A_{n-k}\bigr)\subseteq A_n$, then one can
construct recursively a $K$-linear map $\sigma_l$ by
\begin{equation*}
\sigma_l(a)=-a -
\mu\circ\bigl(\sigma_l\otimes\id\bigr)\circ\Delta'(a)
\end{equation*}
and a $K$-linear map $\sigma_r$ by
\begin{equation*}
\sigma_r(a)=-a -
\mu\circ\bigl(\id\otimes\sigma_r\bigr)\circ\Delta'(a).
\end{equation*}
The maps $\sigma_l, \sigma_r$ might be called left and right
antipodes. (They do not have to coincide, see Remark
(\ref{remcoadantipode})).

\end{remark}

\begin{example}\label{exdiagonal}

Let ${\mathcal P}$ be equipped with a coherent unit action and let
$A$ be the free unitary ${\mathcal P}$-algebra generated by a
vector space $V$. Thus there is a ${\mathcal P}$-algebra
homomorphism $\Delta_a$ defined by
\begin{equation*}
\Delta_a(v)=v\otimes 1 + 1\otimes v, \text{ all } v\in V
\end{equation*}
We call $\Delta_a$ the diagonal or co-addition.

Since $(\Delta_a\otimes\id)\bigl(\Delta_a(v)\bigr)=v\otimes
1\otimes 1 + 1\otimes v\otimes 1 + 1\otimes 1 \otimes
v=(\id\otimes \Delta_a)\bigl(\Delta_a(v)\bigr)$ for all $v\in V$,
the ${\mathcal P}$-algebra homomorphisms
$(\Delta_a\otimes\id)\circ\Delta_a$ and $(\id\otimes
\Delta_a)\circ\Delta_a$ coincide, and we get a strictly graded
${\mathcal P}$-Hopf algebra.

 The name co-addition is chosen in analogy to the
 usage in \cite{libh}, where any
abelian group valued representable functors leads to a co-addition
(on the representing object).

\end{example}

\begin{defn}\label{defcocom}

Let ${\mathcal P}$ be equipped with a coherent unit action such
that for every unitary ${\mathcal P}$-algebra $A$ the $K$-linear
map $\tau: A\otimes A\to A\otimes A, a_1\otimes a_2\mapsto
a_2\otimes a_1$ is a ${\mathcal P}$-algebra isomorphism.

Then, given a unitary ${\mathcal P}$-algebra $A$, we call a
coassociative ${\mathcal P}$-algebra homomorphism $\Delta: A\to
A\otimes A$ cocommutative, if $\tau\circ\Delta=\Delta$.

\end{defn}

\begin{remark}\label{remcocom}

By definition, whenever ${\mathcal P}$ is equipped with a coherent
unit action as in Definition (\ref{defcocom}), the co-addition
$\Delta_a$ is cocommutative.

If ${\mathcal P}$ is $\Com$ or $\As$, this is the classical
definition. We are going to see that $\Dend$ does not fulfill the
condition.

The cocommutative Hopf algebra given by the free unitary
$\As$-algebra $K\langle X\rangle$ together with the diagonal
$\Delta_a$ is well-known. Dually one can equip the standard tensor
coalgebra (i.e.\ $K\langle X\rangle$ with its deconcatenation
coalgebra-structure, see Example (\ref{excofree})) with the
shuffle multiplication $\sh=\Delta_a^{*}$, which is commutative
(cf.\ \cite{lireu}, Chapter 1).

\end{remark}

\begin{defn}

For any strictly graded ${\mathcal P}$-Hopf algebra
$A=\bigoplus_{n\in \N} A^{(n)}$, we define the vector space
\begin{equation*}
A^{*g}=\bigoplus_{n\in \N} (A^{*g})^{(n)}=\bigoplus_{n\in \N}
(A^{(n)})^*
\end{equation*}
where $V^*={\rm Hom}_K (V,K)$ for any vector space $V$.

We call $A^{*g}$ the graded dual of $A$.

We denote by $\Delta^*:A^{*g}\otimes A^{*g}\to A^{*g}$ the
$K$-linear map given by $\bigl(\Delta^*(f_1\otimes
f_2)\bigr)(a)=f_1(a_{(1)})f_2(a_{(2)})$, where $f_1, f_2\in
A^{*g}$ and $a\in A$ with $\Delta(a)=\sum a_{(1)}\otimes a_{(2)}$.

The maps $\Delta^*$ and $\Delta$ are adjoint with respect to the
canonical bilinear form $\langle, \rangle: A^{*g}\times A\to K$,
i.e.\
\begin{equation*}
\langle \Delta^*(f_1\otimes f_2),a\rangle=\langle f_1\otimes
f_2,\Delta(a)\rangle,
\end{equation*}
where $\langle f_1\otimes f_2, g_1\otimes g_2\rangle = \langle
f_1, g_1\rangle\otimes\langle f_2,g_2\rangle$.

\end{defn}

\begin{lemma}\label{lemdualhopf}

The space $A^{*g}$ together with the operations
$\mu^*$, $\mu\in{\mathcal P}$, is a ${\mathcal P}$-coalgebra in the sense of
Definition (\ref{defpcoalg}).

The space $A^{*g}$ together with $\Delta^*:A^{*g}\otimes A^{*g}\to
A^{*g}$ is a graded $\As$-algebra with unit $1\in
(A^{*g})^{(0)}=K$.

Together with its induced ${\mathcal P}$-Hopf algebra structure,
${(A^{*g})}^{*g}$ is isomorphic to $A$.

\end{lemma}

\begin{proof}

This is just the classical result, modified using the fact that a
${\mathcal P}$-algebra structure on $A$ induces a
 ${\mathcal P}$-coalgebra structure on $A^{*g}$.

\end{proof}

\begin{remark}\label{remdualhopf}

In the special case where ${\mathcal P}=\As$, the graded dual
$(A^{*g},\Delta^*)$ of a strictly graded Hopf algebra $(A,\Delta)$
is again a strictly graded Hopf algebra. Furthermore, a
cocommutative $\Delta$ leads to a $\Com$-Hopf algebra
$(A^{*g},\Delta^*)$; and a $\Com$-Hopf algebra $(A,\Delta)$ leads
to a cocommutative $\Delta^*$.

\end{remark}

\end{section}
\goodbreak

\begin{section}{$\Dend$-Hopf algebras of Loday and Ronco}\label{sectdend}

We are going to recall the construction of the free
$\Dend$-algebra $(K\YTree^{\infty}, \prec,\succ)$ from \cite{lilodia},
and we state some of its properties.

The $\Dend$-algebra $(K\YTree^{\infty}, \prec,\succ)$ can be
provided with the structure of a $\Dend$-Hopf algebra. Since the
unit $1$ plays the role of the planar binary tree which consists
of the root, $1$ is usually denoted as $\vert$.

\begin{lemma} (see \cite{lilosci})

For any $\Dend$-algebra $\overline{A}$, $x,y\in \overline{A}$, let

\begin{equation*}
\begin{split}
&\vert\prec y:=0,\ x\prec \vert:=x,\ x\succ \vert:=0,\ \vert\succ
y:=y
\\
&\vert\prec\vert\text{ and }\vert\succ\vert\text{ are not
defined.}\\
\end{split}
\end{equation*}

The induced unit action for $\Dend$ is coherent.

\end{lemma}

\begin{example} (cf.\ \cite{lilrhopftree}, \cite{lilosci})

Consider the vector space $K\YTree^{\infty}$ together with
 operations
$\prec$ and $\succ$, which are recursively defined as follows:

For $T, Z$ planar binary trees $\neq \vert$, we set
\begin{equation*}
T\prec Z=T^l\vee (T^r*Z),\  T\succ Z= (T*Z^l)\vee Z^r
\end{equation*}
where $\vert$ is neutral for the associative operation $*$ (the
sum of $\prec$ and $\succ$).

It is shown in \cite{lilodia} that $(K\YTree^{\infty}, \prec,\succ)$
is the free unitary $\Dend$-algebra on one generator $\Y$.

We consider $(K\YTree^{\infty}, *)$ as a graded associative
algebra, graded with respect to the canonical degree function
induced by
\begin{equation*}
\deg T = n \text{ for } T\in\YTree^{(n)}.
\end{equation*}

The associative algebra $(K\YTree^{\infty}, *)$ is a free
associative algebra with basis $\V(T), T\in \YTree^{\infty}$, cf.\
\cite{lilrhopftree}.

The free $\Dend$-algebra generated by $\Y$ can be provided with a
$\Dend$-Hopf algebra structure.

\end{example}

\begin{lemma}\label{lemdendhopf}(\cite{lilosci})

The free $\Dend$-algebra generated by $\Y$ is a $\Dend$-Hopf
algebra with $\Delta$ determined by

$\Delta(\Y)=\Y\otimes\vert + \vert\otimes\Y$.

\medskip

For $T, Z$ planar binary trees, it holds that
\begin{equation*}
\begin{split}
\Delta(T\vee Z)&= T\vee Z\otimes\vert
 + \sum T_{(1)}* Z_{(1)}\otimes
T_{(2)}\vee Z_{(2)}.\\
\end{split}
\end{equation*}

\end{lemma}

\begin{proof}

The first assertion is clear. For the second, we observe that
\begin{equation*}
\begin{split}
\Delta&(T\vee Z)=\Delta(T\succ\Y\prec
Z)=\Delta(T)\succ(\Y\otimes\vert+\vert\otimes\Y)\prec\Delta(Z)
\\
 & =
\sum(T_{(1)}\otimes
T_{(2)})\succ(\Y\otimes\vert)\prec(Z_{(1)}\otimes Z_{(2)})+ \sum
T_{(1)}* Z_{(1)}\otimes T_{(2)}\succ\!\!\Y\!\!\prec Z_{(2)}\\
 &=
 T\vee Z\otimes\vert + \sum T_{(1)}*
Z_{(1)}\otimes T_{(2)}\vee Z_{(2)},\\
\end{split}
\end{equation*}
because the second tensor component of $(T_{(1)}\otimes
T_{(2)})\succ(\Y\otimes\vert)\prec(Z_{(1)}\otimes Z_{(2)})$ is
zero if not $T_{(2)},Z_{(2)}\in K$, and $T_{(2)},Z_{(2)}\in K$
(inductively) implies $T_{(1)}=T$, $Z_{(1)}=Z$,

$(T_{(1)}\otimes T_{(2)})\succ(\Y\otimes\vert)\prec(Z_{(1)}\otimes
Z_{(2)})=T\succ\Y\prec Z\otimes\vert$.

\end{proof}

\begin{remark} In \ \cite{liroa},\cite{lirob} Ronco has
defined non-unital dendriform Hopf algebras as $\Dend$-algebras
$A$ together with a $K$-linear coassociative map $\Delta':
A\otimes A\to A$, such that:

\begin{equation*}
\begin{split}
 \Delta'(x\succ y) &= x\otimes y + \sum
x_{(1)}*y_{(1)}\otimes x_{(2)}\succ y_{(2)}+ \sum x*y_{(1)}\otimes
x_{(2)}\\
 & + \sum y_{(1)}\otimes x\succ
y_{(2)}+\sum x_{(1)}\otimes x_{(2)}\succ y
\\
\Delta'(x\prec y) &= y\otimes x + \sum x_{(1)}*y_{(1)}\otimes
x_{(2)}\prec y_{(2)} + \sum x_{(1)}*y\otimes x_{(2)} \\
 &+ \sum y_{(1)}\otimes x\prec
y_{(2)}+\sum x_{(1)}\otimes x_{(2)}\prec y  \\
\end{split}
\end{equation*}

After adjoining a unit $\vert$ to the non-unital dendriform Hopf
algebra, the map $\Delta$ given by
 $\Delta(\vert)=\vert\otimes \vert$,
 $\Delta(T)=T\otimes\vert+\vert\otimes T+\Delta'(T)$, $T\neq\vert$,
 is easily seen to be a morphism of $\Dend$-algebras (in the sense
 above).
For example, to check the formula for $\Delta'(x\prec y)$ one can
calculate $\Delta(x)\prec\Delta(y)-x\prec y\otimes \vert
-\vert\otimes x\prec y$ in the unitary case as an expression
consisting of seven terms: Two of these terms, namely the terms
$(x\otimes \vert)\prec(\vert\otimes y)$ and $(x\otimes
\vert)\prec\Delta'(y)$, are zero.

\bigskip

There is an obvious forgetful functor from $\Dend$-Hopf algebras
to $\As$-Hopf algebras.

\end{remark}

\begin{defn}

The $\As$-Hopf algebra structure  on
$\bigl(K\YTree^{\infty},*\bigr)$ given by the formula of
(\ref{lemdendhopf}), i.e.\ with $\Delta_{LR}(\vert)=\vert\otimes
\vert$ and
\begin{equation*}
\Delta_{LR}(T\vee Z)= T\vee Z\otimes\vert + \sum T_{(1)}*
Z_{(1)}\otimes T_{(2)}\vee Z_{(2)}
\end{equation*}
for $T, Z$ planar binary trees, is called the Hopf algebra of
Loday and Ronco, see \cite{lilrhopftree}.

We write $\Delta_{LR}(T)=T\otimes 1+ 1\otimes T +\Delta_{LR}'(T)$.

\end{defn}

\begin{lemma}\label{lemeasylr}
It holds that

\begin{itemize}
\item[(i)]
\begin{equation*}
\Delta_{\rm LR}(\V T)= \V T\otimes\vert + \sum T_{(1)} \otimes
\V(T_{(2)})
\end{equation*}
\item[(ii)]
\begin{equation*}
\Delta_{\rm LR}(T\vee\vert)= T\vee\vert\otimes\vert + \sum T_{(1)}
\otimes T_{(2)}\vee\vert
\end{equation*}
\item[(iii)]
\begin{equation*}
\begin{split}
\Delta_{\rm LR}(\vup&(T^1.T^2\ldots T^n))
\\
 &= \sum_{j=0}^n
T^1_{(1)}*\ldots *T^j_{(1)} * \vup (T^{j+1}\ldots T^n)\otimes\vup
(T^1_{(2)}.T^2_{(2)}\ldots T^j_{(2)})
\end{split}
\end{equation*}
\end{itemize}

\end{lemma}

\begin{proof}

Assertions (i) and (ii) follow directly from the formula of
(\ref{lemdendhopf}).

For (iii) we use induction, where the case $n=1$ is assertion
(ii).

For $j=0$,

$T^1_{(1)}*\ldots * T^j_{(1)} * \vup (T^{j+1}\ldots
T^n)\otimes\vup (T^1_{(2)} \ldots T^j_{(2)})=\vup (T^{1} \ldots
T^n)\otimes\vert.$

\medskip

Since $\vup (T^1 \ldots T^n)=T^1\vee\bigl(\vup (T^2 \ldots
T^n)\bigr)$,
\begin{equation*}
\begin{split}
&\Delta_{\rm LR}(\vup (T^1 \ldots T^n))-\vup (T^{1} \ldots
T^n)\otimes\vert= (*,\vee)\bigl( \Delta(T^1) \otimes \Delta(\vup
(T^2 \ldots T^n))\bigr)\\
 &= \sum_{j=1}^n
T^1_{(1)}*\ldots *T^j_{(1)} * \vup (T^{j+1} \ldots T^n)\otimes\vup
(T^1_{(2)} \ldots T^j_{(2)}).\\
\end{split}
\end{equation*}

\end{proof}

\begin{example}\label{exlrhopf}

For $f=\arbreup - \arbredown=2\V(\Y)-\Y*\Y$ we compute

\begin{equation*}
\Delta_{LR}(f)=V\bigl(\Y\bigr)\otimes\vert+
\vert\otimes\Bigl(\V\bigl(\Y\bigr)-\Y\vee\vert\Bigr) -
\Y\vee\vert\otimes\vert=f\otimes\vert + \vert\otimes f,
\end{equation*}
i.e.\ $f$ is primitive, $\Delta'_{LR}(f)=0$.

The image of $\Y\vee\Y$ under $\Delta'_{LR}$ is
\begin{equation*}
\Y*\Y\otimes\Y+\Y\otimes\Y*\Y.
\end{equation*}

It follows from formula (iii) that  $\Delta_{\rm LR}(\vup (n))=
\sum_{j=0}^n \vup ({n-j})\otimes\vup ({j})$. An analogous formula
holds for left combs.

Combining the combs of height 3 and the tree $\Y\vee\Y$, we get
that the element $g$ given by
\begin{equation*}
\arbrecinq+\arbreun-\arbretrois
\end{equation*}

is primitive.

\medskip

The Loday-Ronco Hopf algebra is not cocommutative.

 We can extend
$\Bigl\{\Y,f,g\Bigr\}$ to an algebra basis of
$\bigl(K\YTree^{\infty},*\bigr)$, by adding one non-primitive
element $h$ of degree $3$ and further elements of degree $\geq 4$.
For example, set

\begin{equation*}
h:=\arbrequatre-\arbredeux
\end{equation*}
with
\begin{equation*}
\begin{split}
\Delta'(h) &=  \arbredown\otimes \Y
-\Y\otimes\arbredown-\arbreup\otimes \Y +\Y\otimes\arbreup\\
 &=
\Y\otimes f - f\otimes \Y.\\
\end{split}
\end{equation*}

\end{example}

\end{section}

\bigskip
\goodbreak

\begin{section}{The Connes-Kreimer Hopf algebra of
renormalization}

\begin{defn}
Let $K\langle\PTree\rangle$ be graded with respect to the
canonical degree function induced by
\begin{equation*}
\deg T = n= \# {\rm Ve}(T), \text{ for } T\in\PTree_n.
\end{equation*}
The empty tree $\emp\in\PTree'=\PTree\cup\{\emp\}$ is identified
with the unit of $K\langle\PTree\rangle$.

\end{defn}

\begin{proposition}

 Let $\underline\Delta:K\langle\PTree\rangle\to
K\langle\PTree\rangle\otimes K\langle\PTree\rangle$ be the graded
unitary algebra homomorphism which maps any tree $T=\vee(T^1\ldots
T^n)\in \PTree$ to
\begin{equation*}
\underline\Delta(T)= T\otimes\emp + \sum T^1_{(1)}\ldots
T^n_{(1)}\otimes \vee\bigl(T^1_{(2)}\ldots T^n_{(2)}\bigr),
\end{equation*}
where we can assume that the images
 $\underline\Delta(T^k)=\sum T^k_{(1)}\otimes
T^k_{(2)}$ of $T^k\in\PTree'=\PTree\cup\{\emp\}$, $k=1,\ldots,n$,
are already defined.
\newline
Then $(K\langle\PTree\rangle,\underline\Delta)$, together with the
augmentation map as counit is a (graded) Hopf algebra.

\end{proposition}

\begin{proof} (Cf.\ also \cite{lifo}.)

 It is easy to check, that $\sum T^1_{(1)}\ldots
T^n_{(1)}\otimes \vee\bigl(T^1_{(2)}\ldots T^n_{(2)}\bigr)$ is
given by $\emp\otimes T$ plus terms $a_{(1)}\otimes a_{(2)}$ with
both $a_{(i)}$ of order $\geq 1$, thus the augmentation map is a
counit. We have to show that $\underline\Delta$ is coassociative.

 As $\underline\Delta(\emp)=\emp\otimes\emp$,
$\bigl( (\underline\Delta\otimes\id)\underline\Delta -
(\id\otimes\underline\Delta)\underline\Delta)\bigr)(\emp)=0$.

Let $T$ be a tree. We can assume that $\bigl(
(\underline\Delta\otimes\id)\underline\Delta -
(\id\otimes\underline\Delta)\underline\Delta)\bigr)(Z)=0$ for all
trees $Z$ of a degree smaller than $\deg(T)$.
\newline
Now $\underline\Delta\otimes\id$ maps
\newline
$T\otimes\emp + (\cdot,\vee)\bigl( \underline\Delta(T^1)
\otimes\ldots \otimes \underline\Delta(T^n)\bigr):=T\otimes\emp +
\sum T^1_{(1)}\ldots T^n_{(1)}\otimes \vee\bigl(T^1_{(2)}\ldots
T^n_{(2)}\bigr)$ onto
\begin{equation*}
\begin{split}
 &\Bigl(T\otimes\emp + (\cdot,\vee)\bigl(
\underline\Delta(T^1) \otimes\ldots \otimes
\underline\Delta(T^n)\bigr)\Bigr)\otimes\emp \\ &+
(\cdot,\cdot,\vee)\bigl((\underline\Delta\otimes\id)\underline\Delta(T^1),\ldots,
(\underline\Delta\otimes\id)\underline\Delta(T^n)\bigr).\\
\end{split}
\end{equation*}
On the other hand, $\id\otimes\underline\Delta$ maps the same
element $T\otimes\emp + (\cdot,\vee)\bigl( \underline\Delta(T^1)
\otimes\ldots \otimes \underline\Delta(T^n)\bigr)$ onto
\begin{equation*}
\begin{split}
T\otimes\emp\otimes\emp + \sum T^1_{(1)}\ldots T^n_{(1)}&\otimes
\underline\Delta\circ\vee\bigl(T^1_{(2)}\ldots T^n_{(2)}\bigr)
\\
=\Bigl(T\otimes\emp +(\cdot,\vee)&\bigl( \underline\Delta(T^1)
\otimes\ldots \otimes \underline\Delta(T^n)
\bigr)\Bigr)\otimes\emp\\
 &+
(\cdot,\cdot,\vee)\bigl((\id\otimes\underline\Delta)\underline\Delta(T^1),\ldots,
(\id\otimes\underline\Delta)\underline\Delta(T^n)\bigr).
\\
\end{split}
\end{equation*}
By induction, it follows that $\underline\Delta$ is coassociative.

\end{proof}

\begin{remark}

In \cite{likr}, Kreimer discovered a $\Com$-Hopf algebra
 for the use of
renormalization of quantum field theories. It was further studied
by  Connes and Kreimer as a $\Com$-Hopf algebra structure on
$K[\ATree]$, cf.\ \cite{lick}, \cite{lickrihi}, and
\cite{likrstruct}.

The Connes-Kreimer Hopf algebra is obtained from $( K\langle{\rm
PTree}\rangle,\underline\Delta)$ by considering the quotient Hopf
algebra with respect to the ideal generated by all commutators and
'forgetting' the planar structure of the trees.
 Its comultiplication $\Delta_{\rm CK}$ is the unitary
algebra homomorphism defined on $K[\ATree]$ by
\begin{equation*}
B_+(T^1\ldots T^n) \mapsto B_+(T^1\ldots T^n)\otimes\emp +
(\cdot,B_+)\bigl( \Delta_{\rm CK}(T^1) \otimes \ldots \otimes
\Delta_{\rm CK}(T^n)\bigr).
\end{equation*}
The graded dual of this $\Com$-Hopf algebra is isomorphic (via a
graded isomorphism) to a noncommutative cocommutative Hopf algebra
on trees introduced by Grossman and Larson \cite{ligla}, see
\cite{lipan},\cite{lihof}.

\end{remark}

\begin{remark}
There is an alternative description of the comultiplication
$\Delta_{\rm CK}$ and also of $\underline\Delta$ using the concept
of admissible cuts (see \ref{defadmcut}). Since any non-full
admissible cut of $T=\vee(T^1\ldots T^n)$ corresponds to $n$
admissible cuts (of $T^1,\ldots, T^n$), it is not hard to prove by
induction that the image of $T$ under $\underline\Delta$ (or
$\Delta_{\rm CK}$) is given by
\begin{equation*}
\displaystyle\sum_{C\text{admissible cut}}C(T)\otimes R^C(T),
\end{equation*}
see \cite{lick}.

\end{remark}

\begin{example}

The comultiplication $\underline\Delta$  maps

\begin{equation*}
\ladderone =\vee(\emp) \mapsto \ladderone\otimes\emp +
\emp\otimes\!\!\!\ladderone
\end{equation*}

For $f=2\laddertwo - \ladderone\ladderone$ we compute
\begin{equation*}
\begin{split}\underline\Delta(f)&=2\Bigl(\laddertwo\otimes \emp +
\emp\otimes\!\!\laddertwo+ \ladderone\otimes\!\!\!\ladderone\Bigr)
- \ladderone\ladderone\otimes\emp -
\emp\otimes\!\!\ladderone\ladderone -
2\ladderone\otimes\!\!\!\ladderone\\ &=f\otimes \emp + \emp\otimes
f\\
\end{split}
\end{equation*}

For $h=2\coratwo - \ladderthree- \laddertwo\!\ladderone$ we
compute
\begin{equation*}
\begin{split}\underline\Delta(h)&=2\Bigl(\coratwo\otimes \emp +
\emp\otimes\!\!\coratwo+
\ladderone\!\!\!\ladderone\otimes\!\!\!\ladderone+ \
2\ladderone\otimes\!\!\!\laddertwo\Bigr)\\ & -
\ladderthree\otimes\emp\ - \emp\otimes\!\!\ladderthree -
\ladderone\otimes\!\!\!\laddertwo-\laddertwo\otimes\!\!\!\ladderone\\
& -\laddertwo\!\ladderone\otimes\emp\ -
\emp\otimes\!\!\laddertwo\!\ladderone - \Bigl(\!\laddertwo +
\ladderone\!\!\!\ladderone\ \Bigr)\otimes\!\!\!\ladderone-
\ladderone\otimes\!\!\!\Bigl(\!\laddertwo +
\ladderone\!\!\!\ladderone\ \Bigr)\\
&=h\otimes \emp + \emp\otimes
h + \ladderone\otimes f - f\otimes\!\!\ladderone
\\
\end{split}
\end{equation*}

The elements $\ladderone$, $f$, and $h$ are in fact images of the
elements $\Y$, $f$, and $h$ of Example (\ref{exlrhopf}) under the
homomorphism $\Theta$ defined in the following.

\end{example}

\begin{defn}
 Let
 $\Theta:(K\YTree^{\infty},*)\to K\langle\PTree\rangle$ be the (graded)
  algebra homomorphism which is
 uniquely defined by $\Theta(\vert)=\emp$, $\Theta(\V(Z))=(\vee\circ
\Theta)(Z), Z\in \YTree^{\infty}$.

\end{defn}

\begin{proposition}\label{propcompa}(\cite{licompa}, \cite{lifo})

The graded algebra homomorphism $\Theta$ provides a Hopf algebra
isomorphism between the Hopf algebra
 $\bigl(K\YTree^{\infty},*,\Delta_{\rm LR}\bigr)$ of \cite{lilrhopftree}
and the noncommutative version $(K\langle{\rm
PTree}\rangle,\underline\Delta)$ of the Hopf algebra of
\cite{lick}.

\end{proposition}

\begin{proof}

 Via $\Theta$,
 $\YTree^{(1)}=\{\V(\vert)\}$ corresponds to $\PTree_{1}=\{\!\ladderone\}$.

 Next, $\YTree^{(2)}=\{\V(\Y)\}$ corresponds to ${\rm
PTree}_{2}=\{\!\laddertwo\}$.

Let $\xi:(K\langle\PTree\rangle,\cdot)\to (K\YTree^{\infty},*)$ be
the algebra homomorphism given by $\xi(\emp)=\vert$,
$\xi\bigl(T)=\V\bigl(\xi(\neg(T))\bigr)$, for planar trees $T$.
Then
$\Theta\bigl(\xi(T)\bigr)=\bigl(\vee\circ\Theta\circ\xi\circ\neg\bigr)(T)=
T$ by induction. Similarly
$(\xi\circ\Theta)\bigl(\V(Z)\bigr)=\V(Z)$, for all $Z$ in
$\YTree^{\infty}$.

It is left to show that $(\Theta\otimes
\Theta)\circ\Delta=\underline\Delta\circ\Theta$. Clearly
$(\Theta\otimes \Theta)(\Delta(\vert))=\emp\otimes
\emp=\underline\Delta(\emp)$. Again, we use induction over the
degree and assume that $(\Theta\otimes
\Theta)(\Delta(S))=\underline\Delta(S)$ for all $S\in
K\YTree^{\infty}$ of degree $< r$.

Let $\V(Z)\in \YTree^{\infty}$ be a tree of degree $r$. Thus
 $Z$ is a tree of degree $r-1$.

$\V(Z)$ is mapped by $(\Theta\otimes \Theta)\circ\Delta$ onto
 $\Theta(\V(Z))\otimes \emp +
\bigl(\Theta\otimes(\Theta\circ\V)\bigr)(\Delta(Z))$

 $=\vee(\Theta(Z))\otimes \emp
 +(id,\vee)\bigl((\Theta\otimes\Theta)\circ\Delta(Z) \bigr)$
 $=\vee(\Theta(Z))\otimes \emp+(id,\vee)\bigl(
 \underline\Delta(\Theta(Z))\bigr)$
 $=\underline\Delta\bigl(\vee(\Theta(Z))\bigr)
 =(\underline\Delta\circ\Theta)\bigl(\V(Z)\bigr)$.

The elements $\V(Z)$ of degree $r$ together with lower degree
elements generate all elements of degree $r$ in
$K\YTree^{\infty}$.
 Thus
 the equation holds, and $\Theta$ is a (graded) Hopf
algebra isomorphism.

\end{proof}

\end{section}

\begin{section}{The Brouder-Frabetti Hopf algebra}\label{sectbf}

The use of a noncommutative \As-Hopf algebra was proposed by C.\
Brouder and A.\ Frabetti \cite{libfb} for renormalization. It is
also a \Dend-Hopf algebra, via a nontrivial isomorphism.

 The Brouder-Frabetti Hopf algebra is
an \As-Hopf algebra structure on the vector space
$K\YTree^{\infty}$, with  multiplication map and comultiplication
$\Delta_{BF}$ both different from the corresponding maps of the
Loday-Ronco Hopf algebra. In \cite{libfa} the name
$\widetilde{H^{\alpha}}$ is given to the opposite algebra and
coalgebra structure.

\begin{defn}

Let the vector space $K\YTree^{\infty}$ be provided with the
operation $\circ_{\alpha}$ as an associative multiplication, see
Remark (\ref{remcircalpha}).

We denote the multiplication
$(\circ_{\alpha}\otimes\circ_{\alpha})\circ\tau_2$ on the
tensor product $K\YTree^{\infty}\otimes K\YTree^{\infty}$
 also by
$\circ_{\alpha}$.

The comultiplication $\Delta_{\rm BF}$ defined on
$\bigl(K\YTree^{\infty},\circ_{\alpha}\bigr)$ is the algebra
homomorphism given by $\Delta_{\rm BF}(\vert)=\vert\otimes \vert$,
$\Delta_{\rm BF}(\Y)=\Y\otimes\vert+\vert\otimes\Y$,
 and, for $T=T^l\vee T^r$ a planar binary tree $\neq \vert$,
\begin{equation*}
\Delta_{\rm BF}\bigl(\V(T)\bigr) =
\V\bigl(T\bigr)\otimes\vert+(\id \otimes \V)\bigl((\Delta_{\rm
BF}(\V(T^r))-\V(T^r)\otimes\vert)\circ_{\alpha}\Delta_{\rm
BF}(T^l)\bigr).
\end{equation*}

\end{defn}

\begin{remark}

We note that it suffices to define $\Delta_{\rm BF}$ on trees of
the form $\V(T)=\vert\vee T$, because these trees generate
$K\YTree^{\infty}$ with respect to $T\circ_{\alpha}S$ (recall that
this is the grafting of
$S$ onto the leftmost leaf of $T$, see (\ref{remcircalpha})).

For example, $\arbredown$ is the product $\Y\circ_{\alpha}\Y$.
Thus $\Delta_{\rm BF}$ maps it on

$\Delta_{\rm BF}(\Y)\circ_{\alpha}\Delta_{\rm BF}(\Y)=\vert\otimes
\arbredown+\arbredown\otimes\vert+ 2\Y\otimes\Y$.

To compute the image of the planar binary tree $\arbrequatre$
under ${\Delta_{\rm BF}}$, we apply the formula above to
$T=\arbredown$ with $T^l=\Y, T^r=\vert$ and get

\begin{equation*}
\begin{split}
&\V(T)\otimes\vert+(\id \otimes \V)\Bigl((\vert\otimes\Y)
\circ_{\alpha}(\vert\otimes\Y+\Y\otimes\vert)\Bigr)\\ &\\
 &=
\arbrequatre\otimes\vert+\vert\otimes\arbrequatre +
\Y\otimes\arbreup
\end{split}
\end{equation*}

Thus $\Delta_{\rm BF}$ is not cocommutative.

It is shown in \cite{libfb,libfa} that $\Delta_{\rm BF}$ is
coassociative.

\end{remark}

\goodbreak

\begin{proposition}\label{propbfformula}(cf.\ \cite{lipa})

For $T=\vup(T^n\ldots T^1)$, $n\geq 1$, $T^j\in \YTree^{\infty}$,
it holds that
\begin{equation*}
\Delta_{\rm BF}(\V(T))=\V(T)\otimes\vert+\sum
\bigl(T^1_{(1)}\circ_{\alpha}\ldots\circ_{\alpha}
T^n_{(1)}\bigr)\otimes \V\bigl(\vup(T^n_{(2)}\ldots
T^1_{(2)})\bigr),
\end{equation*}
where $\Delta_{\rm BF}(S)=\sum S_{(1)}\otimes S_{(2)}$.

\end{proposition}

\begin{proof} We adapt the argument of \cite{lipa} to our notation.

If $T=\vup(T^1)=T^1\vee\vert$, then by definition, $\Delta_{\rm
BF}\bigl(\V(T)\bigr)-\V(T)\otimes\vert$ is equal to
\begin{equation*}
\begin{split}
 & (\id \otimes \V)\bigl((\Delta_{\rm
BF}(\V(\vert))-\V(\vert)\otimes\vert)\circ_{\alpha}(\sum
T^1_{(1)}\otimes T^1_{(2)})\bigr)\\
 &= (\id
\otimes \V)\bigl((\vert\otimes\Y)\circ_{\alpha}(\sum
T^1_{(1)}\otimes T^1_{(2)})\bigr)\\
 &=\sum
T^1_{(1)}\otimes \V(T^1_{(2)}\vee\vert).\\
\end{split}
\end{equation*}
For $T=\vup(T^n\ldots T^1)$, since $T^l=T^n$, it holds that
$\Delta_{\rm BF}\bigl(\V(T)\bigr)-\V(T)\otimes\vert$ is equal to
\begin{equation*}
\begin{split}
  & (\id\otimes \V)\Bigl(\bigl(\Delta_{\rm BF}\bigl(\V(\vup(T^{n-1}\ldots
T^1)\bigr)-\V(\vup(T^{n-1}\ldots
T^1))\otimes\vert\bigr)\circ_{\alpha}\Delta_{\rm BF}(T^n)\Bigr)\\
 &= (\id\otimes \V)\Bigl(\bigl(\sum
T^1_{(1)}\circ_{\alpha}\ldots\circ_{\alpha}T^{n-1}_{(1)}\otimes
\V(\vup(T^{n-1}_{(2)}\ldots
T^1_{(2)}))\bigr)\circ_{\alpha}\Delta_{\rm BF}(T^n)\Bigr)\\
 &=
(\id\otimes \V)\sum
T^1_{(1)}\circ_{\alpha}\ldots\circ_{\alpha}T^{n-1}_{(1)}\circ_{\alpha}T^n_{(1)}\otimes\bigl(
\V(\vup(T^{n-1}_{(2)}\ldots
T^1_{(2)}))\circ_{\alpha}T^n_{(2)}\bigr)\\
 &=
(\id\otimes \V)\sum
T^1_{(1)}\circ_{\alpha}\ldots\circ_{\alpha}T^n_{(1)}\otimes
\vup(T^n_{(2)}\ldots T^1_{(2)})\\
\end{split}
\end{equation*}
by induction on $n$.
\end{proof}

\begin{example}

It follows immediately that all right combs\ \  $\Y$, $\arbreup,
\ldots, \ \ $ are primitive. Therefore
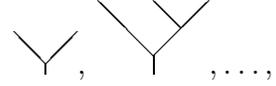
\begin{equation*}
\begin{split}
\Delta_{\rm BF}(\ \arbretrois)&=\Delta_{\rm
BF}(\arbreup)\circ_{\alpha}\Delta_{\rm BF}(\Y)
\\
&=T\otimes\vert+\vert\otimes T +
\Y\otimes\arbreup+\arbreup\otimes\Y.\\
\end{split}
\end{equation*}

\bigskip

Since $\Delta_{\rm BF}\Bigl(\arbredown\Bigr)= \vert\otimes
\arbredown+\arbredown\otimes\vert+ 2\Y\otimes\Y$,

the tree $T=\V\Bigl(\
\arbreun\!\!\!\Bigr)=\V\Bigl(\vup\bigl(\arbredown\bigr)\Bigr)$ is
mapped by $\Delta_{\rm BF}$ onto
\begin{equation*}
T\otimes\vert+\vert\otimes T +2\text{$\Y$ }\otimes\arbrequatre
+\arbredown\otimes\arbreup.
\end{equation*}

\end{example}

\bigskip

\begin{defn}
We define, following \cite{lipa}, for $T^1,\ldots,T^n$ planar
binary trees, the element $\Gamma(T^1,T^2,\ldots,T^n)$ of
$(K\YTree^{\infty},*)$ by
\begin{equation*}
\vup(T^1.T^2\ldots
T^n)-\sum_{j=1}^{n-1}\vup\Bigl(\underbrace{T^1.T^2\ldots T^{j-1}.
\bigl(T^{j}}_j*\vup(T^{j+1}.T^{j+2}\ldots T^n)\bigr)\Bigr),
\end{equation*}
where $\vup$ is applied $K$-linearly. For $n=1$, we get
$\Gamma(T^1)=\vup(T^1)=T^1\vee\vert$.

 Let
 $\Psi:(K\YTree^{\infty},\circ_{\alpha})\to (K\YTree^{\infty},*)$ be the (graded)
  algebra homomorphism which is
 uniquely defined by $\Psi(\vert)=\vert$, $\Psi\bigl(\Y\bigr)=\Y$,
 and
 \begin{equation*}
\Psi\bigl(\V(T)\bigr)
=(\V-\vup)\Bigl(\Gamma(\Psi(T^1),\ldots,\Psi(T^n))\Bigr)
=(\V-\Gamma)\Bigl(\Gamma(\Psi(T^1),\ldots,\Psi(T^n))\Bigr)
 \end{equation*}
if $T=\vup(T^n\ldots T^1)$ is the right comb presentation of
(\ref{remrightcombpres}).

\end{defn}

\begin{proposition}\label{propcompb}

The graded algebra homomorphism $\Psi$ provides a Hopf algebra
isomorphism between the Hopf algebra
$(K\YTree^{\infty},\circ_{\alpha},\Delta_{\rm BF})$ of
\cite{libfb} and the Hopf algebra
 $\bigl(K\YTree^{\infty},*,\Delta_{\rm LR}\bigr)$ of \cite{lilrhopftree}.

\end{proposition}

\begin{proof} We sketch the proof, following \cite{lipa}.

One checks, using Lemma (\ref{lemeasylr}),
\begin{equation*}
\Delta_{\rm LR}(\Gamma(T^1,\ldots,T^n))=
\Gamma(T^1,\ldots,T^n)\otimes\vert + \sum T^1_{(1)}*\ldots
*T^n_{(1)}\otimes\Gamma(T^1_{(2)},\ldots,T^n_{(2)}).
\end{equation*}

The application of $\Delta_{\rm LR}$ to
$\Psi\bigl(\V(T)\bigr)=(\V-\vup)\Bigl(\Gamma
(\Psi(T^1),\ldots,\Psi(T^n))\Bigr)$ and subtraction of the term
$\Psi\bigl(\V(T)\bigr)\otimes\vert$ yields, by Lemma
(\ref{lemeasylr}),
\begin{equation*}
\begin{split}
&\sum \Gamma(\Psi(T^1),\ldots,\Psi(T^n))_{(1)}\otimes
(\V-\vup)\bigl(\Gamma(\Psi(T^1),\ldots,\Psi(T^n))_{(2)}\bigr)
\\
&=\sum (\Psi(T^1))_{(1)}*\ldots*(\Psi(T^n))_{(1)}\otimes
(\V-\vup)\Gamma(\Psi(T^1)_{(2)},\ldots,\Psi(T^n)_{(2)}).\\
\end{split}
\end{equation*}

\bigskip

By induction, we can assume $(\Psi(T^j))_{(1)}=\Psi(T^j_{(1)})$
and $(\Psi(T^j))_{(2)}=\Psi(T^j_{(2)})$. Thus $\Delta_{\rm
LR}\bigl(\Psi\V(T)\bigr)$ is given by
\begin{equation*}
\Psi\V(T)\otimes\vert + \sum
(\Psi(T^1_{(1)})*\ldots*\Psi(T^n_{(1)}))\otimes
(\V-\vup)(\Gamma(\Psi(T^1_{(2)}),\ldots,\Psi(T^n_{(2)}))
\end{equation*}

On the other hand, in view of Proposition (\ref{propbfformula}),
we get that $\Psi\otimes\Psi\bigl(\Delta_{\rm BF}(\V(T))\bigr)$ is
given by
\begin{equation*}
\begin{split}
\Psi\V(T)\otimes\vert &+\sum
\Psi\bigl(T^1_{(1)}\circ_{\alpha}\ldots\circ_{\alpha}
T^n_{(1)}\bigr) \otimes \Psi\V\bigl(\vup(T^n_{(2)}\ldots
T^1_{(2)})\bigr)
\\
 = \Psi\V(T)\otimes\vert &+\sum \bigl(\Psi
T^1_{(1)}*\ldots * \Psi T^n_{(1)}\bigr)\otimes
(\V-\vup)(\Gamma(\Psi(T^1_{(2)}),\ldots,\Psi(T^n_{(2)})).\\
\end{split}
\end{equation*}
Thus $\Psi$ is a Hopf algebra homomorphism. One then checks that
$\Psi$ is a graded bijection.
\end{proof}

\begin{remark}

The construction of a Hopf algebra isomorphism from the
Brouder-Frabetti Hopf algebra $\bigl(K
\YTree^{\infty},\circ_{\alpha},\Delta_{\rm BF}\bigr)$ to the
graded dual of $(K\langle\PTree\rangle,\underline{\Delta})$ is
given in \cite{lifo} (see also \cite{lifoii}).
 There it is
also shown that the Hopf algebra $\bigl(K\langle{\rm
PTree}\rangle,\underline{\Delta}\bigr)$ is self-dual, and thus
also $\bigl(K \YTree^{\infty},\circ_{\alpha},\Delta_{\rm
BF}\bigr)$, $\bigl(K \YTree^{\infty},*,\Delta_{\rm LR})$.

For the Loday-Ronco Hopf algebra, Aguiar and Sottile
(\cite{liasloro}) construct a different basis $\{M_T:
T\in\YTree^{\infty}\}$, called the monomial basis.
They give a different construction for an isomorphism between
$\bigl(K\langle{\rm PTree}\rangle,\underline{\Delta}\bigr)$  and
$\bigl(\YTree^{\infty},\Delta_{\rm LR})$.

\end{remark}

\end{section}

\begin{section}{$\Prim{\mathcal P}$ and the Milnor-Moore Theorem}

Let ${\mathcal P}$ be equipped with a coherent unit action, and
let $V_X$ be the vector space with basis $X=\{x_1,x_2,\ldots\}$.
The co-addition $\Delta_a$ provides the free ${\mathcal
P}$-algebra $F_{\mathcal P}(V_X)$ with the structure of a strictly
graded ${\mathcal P}$-Hopf algebra, see Example
(\ref{exdiagonal}). There exists an operad $\Prim{\mathcal P}$, as
the following Lemma holds.

\begin{lemma}\label{lemprimoperad}(cf.\  \cite{lighmag},\cite{lilosci})

The vector space $\Prim(F_{\mathcal P}(V_X))$ is closed under
insertion: Every ${\mathcal P}$-algebra homomorphism
$\eta_{(g_1,g_2,\ldots)}$ given by $x_i\mapsto g_i$, for $g_i\in
\Prim(F_{\mathcal P}(V_X))$, maps primitive elements on primitive
elements.

\end{lemma}

\begin{proof}

Since $g_i\in \Prim(F_{\mathcal P}(V_X))$ is equivalent to the
compatibility of $\eta_{(g_1,g_2,\ldots)}$ with $\Delta_a$, the
map $\eta_{(g_1,g_2,\ldots)}$ is a ${\mathcal P}$-Hopf algebra
homomorphism, and the assertion follows.

\end{proof}

\begin{defn}(cf.\cite{lilosci})

Let ${\mathcal P}$ be equipped with a coherent unit action.
Let $\Delta_a$ be the co-addition on $F_{\mathcal P}(V_X)$,
$X=\{x_1,x_2,\ldots\}$.

 Then we define $\Prim {\mathcal P}$ to be
the operad with free algebra functor $\Prim F_{\mathcal P}$ and
composition maps induced by insertion.

\end{defn}

\begin{example}

It holds that $\Prim\Com=\mbox{Vect}_K$.
The operad $\Prim\As$ is $\Lie$.

This follows from the Theorem of Friedrichs (cf.\ \cite{lireu}),
which states that Lie polynomials are exactly the polynomials in
non-commuting associative variables which are primitive under
$\Delta_a$.

\end{example}

\begin{remark}

If $A$ is a cocommutative connected graded $\As$-Hopf algebra
(over a characteristic 0 field $K$), then
\begin{equation*}
A\cong U(\Prim(A))
\end{equation*}
is isomorphic to the universal enveloping bialgebra of its space
of primitive elements. This is the theorem of Milnor-Moore
(\cite{limimo}, see also \cite{liqui}).

Together with the theorem of Poincar\'e-Birkhoff-Witt, it follows
that the category of such Hopf algebras is equivalent to the category
of Lie algebras.

An analogon of the Milnor-Moore theorem has been proven by Ronco
\cite{lirob} for $\Dend$-Hopf algebras (compare also
\cite{lichamimo}).

The role of $\Lie$-algebras is played by $\Brace$-algebras
(defined below), special pre-Lie algebras equipped with
 $n$-ary operations $\langle \ldots \rangle: A^{\otimes n}\to A$ for each $n$.

The pre-Lie bracket $\langle\ ,\ \rangle: \Prim(F_{\Dend}
V)^{\otimes 2}\to \Prim(F_{\Dend} V)$ is given by
\begin{equation*}
\langle f, g\rangle=f\prec g - g\succ f.
\end{equation*}
Ronco's results show that the primitive elements of the
Loday-Ronco Hopf algebra, the noncommutative planar Connes-Kreimer
Hopf algebra, and the Hopf algebra of Brouder-Frabetti, are a free
brace algebra on one generator. It follows that
$\Prim\Dend=\Brace$ (cf.\ also \cite{lilosci}).

For different constructions of the space of primitive elements see
also \cite{lifo}, \cite{lias}, and \cite{lihnt}.

In \cite{lilr04} Loday and Ronco prove a Milnor-Moore and
Poincar\'e-Birkhoff-Witt theorem for Hopf algebras $H$ that are
equipped with two multiplications $*$ and $\cdot$, such that the
comultiplication $\Delta$ is a homomorphism with respect to $*$,
and such that $(H,\cdot,\Delta)$ is unital infinitesimal:
\begin{equation*}
\Delta(x\cdot y)=(x\otimes 1)\cdot \Delta(y) + \Delta(x)\cdot
(1\otimes y) - x\otimes y
\end{equation*}

They show that (non-dg) $B_{\infty}$-algebras occur as the primitive
elements.

We sketch the definition of (non-dg) $B_{\infty}$-algebras below, following
\cite{lilr04}.

\end{remark}

\begin{defn} (cf.\ \cite{lilr03},\cite{liroa},\cite{lirob})

Let $A$ be a $K$-vector space together with a family
\begin{equation*}
M_{pq}: A^{\otimes p}\otimes A^{\otimes q}\to A,\ \  p\geq 0,
q\geq 0
\end{equation*}
of $(p+q)$-ary operations be given, with
\begin{equation*}
M_{00}=0, \ M_{10}=\id_A=M_{01}, \text{ and } M_{n0}=0=M_{0n}
\text{ for } n\geq 2.
\end{equation*}

For all $k\geq 0$, and all $i=(i_1,i_2,\ldots,i_k),
j=(j_1,j_2,\ldots,j_k)$, one denotes by
\begin{equation*}
M_{i_1j_1}M_{i_2j_2}\ldots M_{i_kj_k}: A^{\otimes p}\otimes
A^{\otimes q}\to A^{\otimes k}, \text{ with } p:=i_1+\ldots+i_k,
q:=j_1+\ldots+j_k,
\end{equation*}
the map which sends $u_1.u_2\ldots u_p\otimes v_1.v_2\ldots v_q$
to
\begin{equation*}
M_{i_1j_1}(u_1\ldots u_{i_1}\otimes v_1\ldots
v_{j_1}).M_{i_2j_2}(\ldots u_{i_1+i_2}\otimes\ldots
v_{j_1+j_2})\ldots M_{i_kj_k}(\ldots u_p\otimes \ldots v_q).
\end{equation*}

Let $*:A^{\otimes p}\otimes A^{\otimes q}\to A$ be given by
\begin{equation*}
u_1.u_2\ldots u_p * v_1.v_2\ldots v_q=\sum_{k\geq 1}\sum_{i,j}
M_{i_1j_1}\ldots M_{i_kj_k}(u_1.u_2\ldots u_p\otimes v_1.v_2
\ldots v_q),
\end{equation*}
where the sum is over all $k$-tuples $i,j$ with $p=i_1+\ldots+i_k,
q=j_1+\ldots+j_k$.

Then
\begin{equation*}
\begin{split}
u*v&=M_{11}(u\otimes v)+u.v+v.u
\\
u.v*w&=M_{21}(u.v\otimes w)+(u*w).v+u.(v*w)-u.w.v
\\
u*v.w&=M_{12}(u\otimes v.w)+(u*v).w+v.(u*w)-v.u.w.
\\
\end{split}
\end{equation*}

The space $A$ together with the operations $M_{pq}$ is called a
$B_{\infty}$-algebra, iff for every triple $(i,j,k)$ of positive
integers the relation
\begin{equation*}
(u_1.u_2\ldots u_i*v_1.v_2\ldots v_j)*w_1.w_2\ldots w_k
=u_1.u_2\ldots u_i*(v_1.v_2\ldots v_j*w_1.w_2\ldots w_k)
\end{equation*}
holds.

For $(i,j,k)=(1,1,1)$ this is
\begin{equation*}
\begin{split}
M_{11}\bigl(u\otimes M_{11}(v\otimes w)\bigr) &+
M_{12}\bigl(u\otimes(v.w+w.v)\bigr)
\\
 &=M_{21}\bigl((u.v+v.u)\otimes w\bigr)+M_{11}\bigl(M_{11}(u\otimes
v)\otimes w\bigr).\\
\end{split}
\end{equation*}

If $M_{pq}=0$ for all $(p,q)$ with $p\geq 2$, then $A$ is called a
brace algebra, and one denotes the $n$-ary operations $M_{1,n-1}$
for $n\geq 2$ by
\begin{equation*}
\langle \ldots \rangle: A^{\otimes n}\to A.
\end{equation*}

Then, e.g.\ the relation for $(i,j,k)=(1,1,1)$, reads
\begin{equation*}
\langle\langle u.v\rangle.w\rangle=\langle u.\langle v.w
\rangle\rangle + \langle u.v.w\rangle +\langle u.w.v\rangle.
\end{equation*}

Interchanging $v$ and $w$ one gets the pre-Lie relation

\begin{equation*}
\langle\langle u.v\rangle.w\rangle-\langle u.\langle
v.w\rangle\rangle =\langle\langle u.w\rangle.v\rangle-\langle
u.\langle w.v\rangle\rangle,
\end{equation*}
cf. \cite{lichaliv}, \cite{ligers}.

\end{defn}

\end{section}

\end{chapter}

\begin{chapter}{Primitive elements of $\Mag$- and $\Minf$-algebras}

\begin{section}{The free $\Mag$- and $\Minf$-algebras}

Let $K$ be a field of characteristic 0, and let
$X=\{x_1,x_2,\ldots\}$ be a finite or countable set of variables.
Let $V_X$ be the vector space with basis $X$.

\begin{remark}

We recall, see Example (\ref{exminf}),  that the $n$-th space of
the Stasheff operad $\Minf$, i.e.\ the space $\Gamma(M)(n)$ of the
free operad $\Gamma(M)$ generated by operations $\vee^k\in M_k,
k\geq 2$, is given by
 linear
combinations of reduced (non-labeled) planar trees with $n$
leaves. The operad structure is given by insertion, see
(\ref{exminf}).

Super-Catalan numbers are the coefficients of the generating
series of $\Minf$:
\begin{equation*}
f^{\Minf}(t)=\sum_{n\geq 1}\frac{\dim {\Minf}(n)}{n!}t^n
=\sum_{n\geq 1}C_n t^n= \frac{1}{4}(1+t-\sqrt{1-6t+t^2}).
\end{equation*}

For any $k$, the operation $\vee^k$ is given by the grafting
operation $\vee$ restricted to planar forests consisting of $k$
trees.

Let $M^X=(M^X_k)_{k\geq 0}$ denote the sequence
\begin{equation*}
M^X_0=X,\ \  M^X_1=\emp,\ \ M^X_k=M_k=\{\vee^k\} \text{ (for }
k\geq 2).
\end{equation*}
Then the free $\Minf$-algebra
\begin{equation*}
F_{\Minf}(V_X)={\displaystyle \bigoplus_{n=1}^{\infty}}
F_{\Minf}(V_X)^{(n)}={\displaystyle \bigoplus_{n=1}^{\infty}}
\Minf(n)\otimes_{\Sigma_n} \bigl(V_X\bigr)^{\otimes n}
\end{equation*}
has the set $\PRTrx$ of reduced planar trees with leaves labeled
by $X$ as a vector space basis. It is naturally graded, such that
the planar trees with $n$ leaves form a basis of
$F_{\Minf}(V_X)^{(n)}$.

\end{remark}

\begin{lemma}

 The operad $\Minf$ is equipped with a (unique) unit action which
respects the operations $\vee^k$, i.e.\
\begin{equation*}
\begin{split}
\vee^k \circ_{k,i}(\eta) &=\vee^{k-1}, \text{ all } i\\
 \vee^2\circ_{2,1}(\eta) &=\vee^2\circ_{2,2}(\eta)=\id.\\
\end{split}
\end{equation*}
This  unit action is coherent.
\end{lemma}
\begin{proof}

Since $\Minf$ is freely generated by the operations $\vee^k$, we
can define the unit action by the formulas above, see Definition
(\ref{defrespect}). The unit action is compatible with the (empty)
set of relations of $\Minf$, and coherence also follows.

\end{proof}

\begin{defn}
We denote the free $\Minf$-algebra with unit $1$ by
\begin{equation*}
\begin{split}
&K\{X\}_{\infty}={\displaystyle \bigoplus_{n=0}^{\infty}}
K\{X\}_{\infty}^{(n)},\\ &K\{X\}_{\infty}^{(0)}= K 1, \ \
K\{X\}_{\infty}^{(n)}=F_{\Minf}(V_X)^{(n)}\ (n\geq 1).\\
\end{split}
\end{equation*}
We set $\vee^k(1,1,\ldots,1)=1$ (all $k$).

We identify $1$ with the empty tree $\emp$ and call
\begin{equation*}
\PRTrprimex=\PRTrx\cup\{\emp\}
\end{equation*}
 the set of monomials of $K\{X\}_{\infty}$.

\end{defn}

\begin{remark}

The binary quadratic operad $\Mag$ with generating series
\begin{equation*}
f^{\Mag}(t)=\sum_{n\geq 1}\frac{\dim {\Mag}(n)}{n!}t^n
=\sum_{n\geq 1}c_n t^n=\frac{1-\sqrt{1-4t}}{2}
\end{equation*}
is contained as a sub-operad in $\Minf$, i.e.\ there is an
inclusion of operads
\begin{equation*}
\Mag\to\Minf.
\end{equation*}
(The vector subspace given by planar binary trees with labeled
leaves is closed under the operation $\vee^2$.)

The unit action of $\Minf$ extends the canonical unit action on
$\Mag$ induced by $1\cdot a= a\cdot 1 = a$ (for every
$\Mag$-algebra $A$, every $a\in A$), see (\ref{remcanunit}).

\medskip

The free $\Mag$-algebra with unit 1 can be identified with the
space of labeled binary trees $K\{X\}\subset K\{X\}_{\infty}$,
equipped with the free binary operation $\cdot=\vee^2$.

We pass freely from planar binary trees to parenthesized strings,
see (\ref{exmalcev}). We call a planar binary tree with labeled
leaves a monomial of $K\{X\}$, and we call the elements of
$K\{X\}$ (and $K\{X\}_{\infty})$ polynomials.

We can embed  $K\{X\}_{\infty}$ into its completion
$K\{\{X\}\}_{\infty}=\prod_{n=0}^{\infty} K\{X\}_{\infty}^{(n)}$,
the free complete $\Minf$-algebra generated by $X$. Similarly
defined is the free $\Mag$-power series algebra $K\{\{X\}\}$.

\medskip

We have the canonical degree and order functions.

The degree of a tree $T\in K\{X\}_{\infty}^{(n)}$ is $n$, i.e.\
the degree of $T$ is the number of its leaves.

The multi-degree of a tree $T\in K\{X\}_{\infty}$ is ${\bf n}=
(n_1,n_2,\ldots)$, if $n_j$ is the number of leaves with label
$x_j$ for all $j\in\N$.

\end{remark}

\begin{remark}
While $\As(n)=K\Sigma_n$ corresponds to the regular
representation, the $\Sigma_n$-module $\Mag(n)$ is given by $c_n$
copies of $K\Sigma_n$, and $\Minf(n)$ is given by $C_n$ copies of
$K\Sigma_n$.

The associated $GL$-modules describe the
 spaces $K\langle X\rangle^{(n)}$, $K\{X\}^{(n)}$, and $K\{X\}_{\infty}^{(n)}$.

The first degree where $K\{X\}^{(n)}$ and $K\{X\}_{\infty}^{(n)}$
differ is degree 3.

Here the regular representation
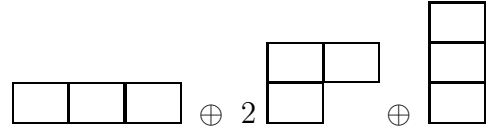
\begin{equation*}
 {\rm \vbox{\offinterlineskip
\hbox{\boxtext{\phantom{X }}\boxtext{\phantom{X
}}\boxtext{\phantom{X }}}
 } \rm}\ \oplus\ 2
 \ {\rm
 \vbox{\offinterlineskip \hbox{\boxtext{\phantom{X
}}\boxtext{\phantom{X }}} \hbox{\boxtext{\phantom{X }}}
 } \rm}
\oplus\  {\rm
 \vbox{\offinterlineskip \hbox{\boxtext{\phantom{X
}}} \hbox{\boxtext{\phantom{X }}} \hbox{\boxtext{\phantom{X }}} }
\rm}
\end{equation*}
occurs in two copies for $K\{X\}^{(3)}$, and in three copies for
$K\{X\}_{\infty}^{(3)}$. The extra copy is induced by the
3-corolla.

\medskip

For $X=\{x_1,x_2,\ldots,x_m\}$, the Hilbert series
$\Hilb(A,t_1,\ldots,t_m)$ of the multigraded vector spaces
$A=K\{X\}_{\infty}$ and $A=K\{X\}$ is defined by
\begin{equation*}
\Hilb(A,t_1,\ldots,t_m)=\sum_{{\bf n}\in \N^m} \dim
A^{(n_1,n_2,\ldots,n_m)}t_1^{n_1}t_2^{n_2}\cdots t_m^{n_m}
\end{equation*}
where the vector space $A^{(n_1,n_2,\ldots,n_m)}$ is the space of
homogeneous elements having multi-degree $(n_1,n_2,\ldots,n_m)$.

Again, the coefficients for $K\{X\}_{\infty}$ and $K\{X\}$ are
just given by the coefficients for $K\langle X\rangle$, multiplied
by $C_n$ or $c_n$.

\end{remark}

\end{section}

\bigskip
\goodbreak

\begin{section}{Partial derivatives on $\Minf$-algebras}

\begin{defn}\label{defderiv}
 A
derivation $D$ on a $\Mag$-algebra $A$ is a $K$-linear map which
fulfills the Leibniz rule
\begin{equation*}
D(u\cdot v)=D(u)\cdot v+u\cdot D(v) \text{ for all } u,v\in A.
\end{equation*}

A derivation $D$ on a $\Minf$-algebra $A$ is a $K$-linear map
which fulfills
\begin{equation*}
D(\vee^n(v_1,\ldots,v_n))=\sum_{i=1}^n
\vee^{n}(v_1,\ldots,v_{i-1},D(v_i),v_{i+1},\ldots,v_n) \text{ for
all } v_i\in A.
\end{equation*}

\end{defn}

\begin{lemma}

Every mapping $\delta: X\to K\{X\}$, or $X\to K\{X\}_{\infty}$
respectively, can be extended to a unique derivation of $K\{X\}$,
$K\{X\}_{\infty}$ respectively. Here $\delta\vert_K=0$.

\end{lemma}

\begin{proof}
As usual, one defines $\delta(u\cdot v)=\delta(u)\cdot v+u\cdot
\delta(v)$, or
$\delta\bigl(\vee^n(v_1,\ldots,v_n)\bigr)=\sum_{i=1}^n
\vee^{n}(v_1,\ldots,v_{i-1},\delta(v_i),v_{i+1},\ldots,v_n)$ for
$u,v,v_i\in X$.

\end{proof}

\begin{defn}

Let $\partial_k:K\{X\}_{\infty}\to K\{X\}_{\infty}$
 be the derivation given by

\begin{equation*}
\partial_k(x_l):=\begin{cases} 1 &:\ k=l \\ 0 &:\ k\neq l\\
\end{cases}
\end{equation*}

Let $\partial_{kj}:K\{X\}_{\infty}\to K\{X\}_{\infty}$
 be the derivation given by

\begin{equation*}
\partial_{kj}(x_l):=\begin{cases} x_j &:\ k=l \\ 0 &:\ k\neq l\\
\end{cases}
\end{equation*}

The restrictions on $K\{X\}$ are also denoted by
$\partial_k,\partial_{kj}$.

\end{defn}

\begin{example}

Let $f$ be the monomial $(x_1\cdot ((x_1\cdot x_2)\cdot x_2))$,
which corresponds to the planar binary tree
$
 
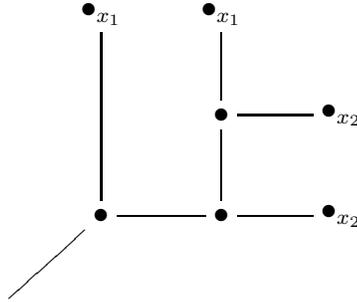
\xymatrix{  & {\bullet}_{x_1}\ar@{-}[dd] &{\bullet}_{x_1}\ar@{-}[d] \\ & &
{\bullet}\ar@{-}[d]\ar@{-}[r]  & {\bullet}_{x_2} \\
 & {\bullet}\ar@{-}[dl] &{\bullet}\ar@{-}[l]  &{\bullet}_{x_2}\ar@{-}[l]  \\ \\} \
$

Then $\partial_2(f)$ is given by $2(x_1\cdot (x_1\cdot x_2))$ or
$
 
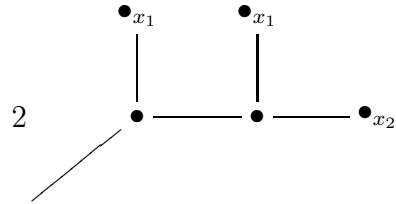
\xymatrix{  & {\bullet}_{x_1}\ar@{-}[d] &{\bullet}_{x_1}\ar@{-}[d] \\
2\ & {\bullet}\ar@{-}[dl] &{\bullet}\ar@{-}[l]
&{\bullet}_{x_2}\ar@{-}[l]  \\ \\} \
$

Then $\partial_1(f)$ is given by $(x_1\cdot (x_2\cdot
x_2))+((x_1\cdot x_2)\cdot x_2)$, i.e.\

$
 \xymatrix{  & {\bullet}_{x_1}\ar@{-}[d] &{\bullet}_{x_2}\ar@{-}[d] \\
\ & {\bullet}\ar@{-}[dl] &{\bullet}\ar@{-}[l]
&{\bullet}_{x_2}\ar@{-}[l]  \\ \\} \
 \xymatrix{  &  &{\bullet}_{x_1}\ar@{-}[d] \\ &+ &
{\bullet}\ar@{-}[d]\ar@{-}[r]  & {\bullet}_{x_2} \\
 & & {\bullet}\ar@{-}[dl]  &{\bullet}_{x_2}\ar@{-}[l]  \\ & \\} \
$
\end{example}

\begin{example}\label{exderiveight}
If $f$ is the monomial
\begin{equation*}
\vee^2\Bigl(\vee^2\bigl(x_1,
\vee^3(x_2,x_2,x_2))\bigr),\vee^2\bigl(\vee^3(x_2,x_2,x_1),x_2\bigr)\Bigr),
\end{equation*}
which corresponds to

 $
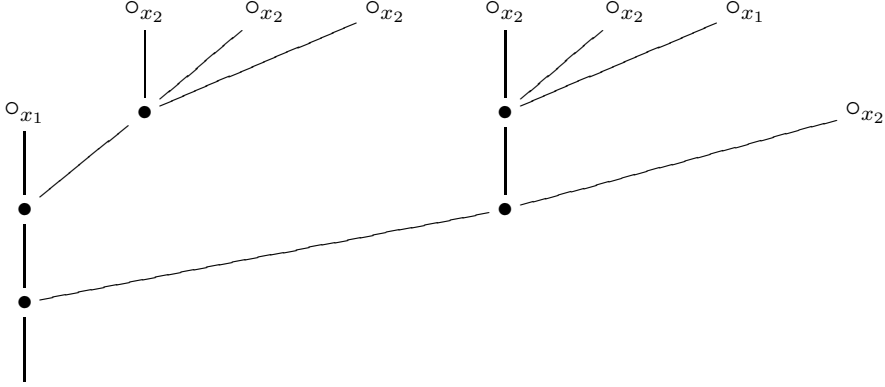
\xymatrix{ & \ar@{-}[d]\circ_{x_2} & \ar@{-}[dl]\circ_{x_2} &
 \ar@{-}[dll]\circ_{x_2}
&\circ_{x_2} \ar@{-}[d] & \ar@{-}[dl]\circ_{x_2} &
\ar@{-}[dll]\circ_{x_1}&
 \\
 \circ_{x_1} \ar@{-}[d] & \bullet \ar@{-}[dl]& & & \bullet \ar@{-}[d]& & &\ar@{-}[dlll]\circ_{x_2}\\
  \bullet\ar@{-}[d] & & & &
  \bullet\ar@{-}[dllll] & & \\
  \bullet\ar@{-}[d]& & & & \\
  \\}$

then $\partial_{12}(f)$ is given by
\begin{equation*}
\begin{split}
&\vee^2\Bigl(\vee^2\bigl(x_2,
\vee^3(x_2,x_2,x_2))\bigr),\vee^2\bigl(\vee^3(x_2,x_2,x_1),x_2\bigr)\Bigr)\\
&+\vee^2\Bigl(\vee^2\bigl(x_1,
\vee^3(x_2,x_2,x_2))\bigr),\vee^2\bigl(\vee^3(x_2,x_2,x_2),x_2\bigr)\Bigr).\\
\end{split}
\end{equation*}
and $\partial_{1}(f)$ is given by
\begin{equation*}
\begin{split}
&\vee^2\Bigl(\vee^3(x_2,x_2,x_2)),\vee^2\bigl(\vee^3(x_2,x_2,x_1),x_2\bigr)\Bigr)\\
&+\vee^2\Bigl(\vee^2\bigl(x_1,
\vee^3(x_2,x_2,x_2))\bigr),\vee^2\bigl(\vee^2(x_2,x_2),x_2\bigr)\Bigr).\\
\end{split}
\end{equation*}

\end{example}

\begin{proposition}\label{proppartial}

Let $T\in K\{X\}_{\infty}$ be a monomial, i.e.\ an element of the
set $\PRTrx$, $x_k\in X$.

Let $I_k\subseteq {\rm Le}(T)$ be the subset given by all leaves
with label $x_k$.

Then, for any $j$, $\partial_{kj}(T)$ is given by
\begin{equation*}
\partial_{kj}(T)=\sum_{\nu\in I_k} T_{\nu\to x_j}
\end{equation*}
where $T_{\nu\to x_j}$ denotes the tree which is obtained from $T$
by labeling the leaf $\nu$ by $x_j$ (instead by $x_k$).

If, for $\nu\in I_k$, $\{\nu\}^c$ denotes the complement ${\rm
Le}(T)-\{\nu\}$ of $\nu$ in the set of leaves of $T$, then
$\partial_k(T)$ is given by
\begin{equation*}
\partial_k(T)=\sum_{\nu\in I_k} {\rm red}(T\vert\{\nu\}^c)
\end{equation*}
where ${\rm red}(T\vert\{\nu\}^c)$ is the reduction of the
leaf-restriction of $T$ onto $\{\nu\}^c$, see
(\ref{lemleafrestrict}) and (\ref{defleafsplit}).

\end{proposition}

\end{section}

\begin{section}{Co-addition $\Mag$- and $\Minf$-Hopf algebras}

Since the unit actions for $\mathcal P\in \{\Mag,\Minf\}$ are
coherent, it makes sense to look at the $\mathcal P$-algebras
$K\{X\}_{\infty}\otimes K\{X\}_{\infty}$ and $K\{X\}\otimes
K\{X\}$.

Both $K\{X\}$ and $K\{X\}_{\infty}$ are naturally graded,
$K\{X\}_{\infty}=\bigoplus_{n\in \N}K\{X\}_{\infty}^{(n)}$,
$K\{X\}=\bigoplus_{n\in \N}K\{X\}^{(n)}$. Homogeneous elements of
degree $n$ are the $K$-linear combinations of trees with $n$
leaves. The operations $\vee^n$ clearly respect this grading.

\begin{defn}\label{defcoadd}

Let $\Delta_a: K\{X\}_{\infty}\to K\{X\}_{\infty}\otimes
K\{X\}_{\infty}$
 be the $\Minf$-algebra homomorphism defined by
\begin{equation*}
  x_i\mapsto x_i\otimes 1 + 1\otimes
 x_i, \text{ for all } x_i\in X.
\end{equation*}

The map $\Delta_a$ is called the co-addition.

\end{defn}

\begin{remark}

We note that $\Delta_a$ is cocommutative. The notion of
cocommutativity is defined for $\Minf$ and $\Mag$, see Definition
(\ref{defcocom}).

The restriction $\Delta_a\vert_{K\{X\}}$ is a $\Mag$-algebra
homomorphism
\begin{equation*}
\Delta_a: K\{X\}\to K\{X\}\otimes K\{X\},
\end{equation*}
which is
also cocommutative.

For example, to compute the image of $(x_1\cdot x_2)\cdot x_3$ under
$\Delta_a$,
one computes
\begin{equation*}
\begin{split}
\Delta_a(x_1\cdot x_2)&=(x_1\otimes 1 + 1\otimes x_1)
\cdot (x_2\otimes 1 + 1\otimes x_2)\\
&= x_1\cdot x_2\otimes 1
+ 1\otimes x_1\cdot x_2 + x_1\otimes x_2+x_2\otimes x_1\\
\Delta_a\bigl((x_1\cdot x_2)\cdot x_3\bigr)&=
(x_1\cdot x_2\otimes 1 + 1\otimes x_1\cdot x_2
+ x_1\otimes x_2+x_2\otimes x_1)\cdot(x_3\otimes 1+1\otimes x_3)\\
&= (x_1\cdot x_2)\cdot x_3\otimes 1 + 1\otimes (x_1\cdot x_2)\cdot x_3
+ x_1\cdot x_2\otimes x_3 + x_1\otimes x_2\cdot x_3 \\
&\ + x_2\otimes x_1\cdot x_3
+ x_3\otimes x_1\cdot x_2 + x_2\cdot x_3\otimes x_1 +x_1\cdot x_3\otimes x_2.
\\
\end{split}
\end{equation*}

\end{remark}

\begin{proposition}\label{propminfhopf}

The free unitary $\Minf$-algebra $K\{X\}_{\infty}$ on $X$ together
with $\Delta_a$ is a strictly graded $\Minf$-Hopf algebra.

The free unitary $\Mag$-algebra $K\{X\}\subseteq K\{X\}_{\infty}$
 is a strictly graded
$\Mag$-Hopf algebra.

\end{proposition}

\begin{proof}

We recall from (\ref{exdiagonal}) that
$(\Delta_a\otimes\id)\Delta_a=(\id\otimes\Delta_a)\Delta_a$,
because both homomorphisms map $x_i$ on
 $x_i\otimes 1\otimes 1 + 1\otimes
 x_i\otimes 1+ 1\otimes 1\otimes x_i$.
Here we note that the map $\Delta_a$ is (by definition) a
$\Minf$-algebra homomorphism.

It follows that
\begin{equation*}
\Delta'_a\bigl(K\{X\}_{\infty}^{(n)}\bigr) \subseteq
\sum_{i=1}^{n-1} K\{X\}_{\infty}^{(i)}\otimes
K\{X\}_{\infty}^{(n-i)}
\end{equation*}
using the $\Minf$-structure of the tensor product together with
the fact that the grading is respected by the operations $\vee^n$.

For the case of $\Mag$, the analogous statements hold.

\end{proof}

\begin{remark}\label{remcoadantipode}

Consider the $\Mag$-Hopf algebra  $(K\{X\}, \Delta_a)$. There is a
unique $K$-linear map $\sigma_l$, the left-antipode, such that
$\sigma_l(w) + w +\sum' \sigma_l(w_{(1)})w_{(2)}=0$, $\deg w \geq
1$, where $\Delta_a(w)=w\otimes 1+1\otimes w +\sum' w_{(1)}\otimes
w_{(2)}$, see Remark (\ref{remantipode}).

For primitive $f$, $\sigma_l(f)=-f$. But $\sigma_l$ is not the
anti-homomorphism induced by $x_i\mapsto -x_i$.

Let $\bar t$ denote the involution induced by
$\overline{t_1t_2}=\bar t_2 \bar t_1$, compare Remark
(\ref{remdegraftmirror}).

Then $\sigma_r:=\bar \sigma_l$ given by $\bar \sigma_l (\bar
t)=\overline{\sigma_l(t)}$ fulfills $\bar \sigma_l(w) + w +\sum
w_{(1)}\bar  \sigma_l(w_{(2)})=0$. Therefore $\sigma_r$ is the
right-antipode. Left and right-antipode do not coincide.

 For
example, let $X=\{x\}$ consist of one element. Then
$\sigma_l(x)=-x$, $\sigma_l(x\cdot x)=x\cdot x$, and we
recursively compute that $\sigma_l(x\cdot (x\cdot x))=2x\cdot
(x\cdot x)-3(x\cdot x)\cdot x$, $\sigma_l((x\cdot x)\cdot
x)=3x\cdot (x\cdot x)-4(x\cdot x)\cdot x$.

On the other hand $\sigma_r(x\cdot (x\cdot x))=3(x\cdot x)\cdot x
-4x\cdot (x\cdot x)$, $\sigma_r((x\cdot x)\cdot x)=2(x\cdot
x)\cdot x -3x\cdot (x\cdot x)$.

 The order of $\sigma_l$ (and also the order of $\sigma_r$)
is infinite.

\end{remark}

\begin{proposition}\label{propdelta}

Let $T\in K\{X\}_{\infty}$ be a monomial.

 For $I\subseteq {\rm Le}(T)$, let
$\bigl( {\rm red} (T\vert I), {\rm red}(T\vert I^c) \bigr)$ be the
leaf-split induced by $(I,I^c)$, see (\ref{defleafsplit}).

Then
\begin{equation*}
\Delta_a(T)=\sum_{I\subseteq {\rm Le}(T)}{\rm red} (T\vert
I)\otimes {\rm red}(T\vert I^c).
\end{equation*}

The same formula holds for the restriction of $\Delta_a$ to
$K\{X\}$, and $T\in K\{X\}$ a monomial.

\end{proposition}

\begin{proof}

For $T=1$, ${\rm Le}(T)=\emptyset=I$ and the formula says that
$\Delta_a(1)=1\otimes 1$. For $T=x_i\in X$,
\begin{equation*}
\Delta_a(T)=T\otimes 1+ 1\otimes T= {\rm red}(T\vert {\rm
Le}(T))\otimes {\rm red}(T\vert \emptyset)+ {\rm red}(T\vert
\emptyset)\otimes{\rm red}(T\vert {\rm Le}(T))
\end{equation*}
Else $T=\vee^{n}(T^1\ldots T^n)$, $n\geq 2$. Since $\Delta_a$ is
by definition a $\Minf$-algebra homomorphism,
\begin{equation*}
\Delta_a(T)=\vee^{n}\bigl( \Delta_a(T^1)\ldots
\Delta_a(T^n)\bigr).
\end{equation*}
After an iterated use of the $\Minf$-algebra homomorphism
property, we arrive at a sum that has one summand for every subset
$I\subseteq {\rm Le}(T)$. Here the choice between letting the
$j$-th leaf belong to $I$ or its complement $I^c$ corresponds to
the choice of a factor $x_{\nu(j)}\otimes 1$ instead of a factor
$1\otimes x_{\nu(j)}$, where $x_{\nu(j)}$ is the label of the
$j$-th leaf of $T$.

By construction, the reduced leaf-restriction ${\rm red} (T\vert
I)$ is the first tensor component of the summand given by $I$ in
$\Delta_a(T)$. Analogously, ${\rm red} (T\vert I^c)$ is the second
tensor component.

The same proof holds, when $K\{X\}_{\infty}$ is replaced by
$K\{X\}$ and
 $\Minf$ is replaced by $\Mag$.

\end{proof}

\begin{defn}\label{defdiffop}

 For any monomial $T$ in
$K\{X\}_{\infty}$ we define a $K$-linear map $\partial_T:
K\{X\}_{\infty}\to  K\{X\}_{\infty}$ by
\begin{equation*}
\Delta_a(f)=\sum_{T\in\PRTrx} T\otimes
\partial_T(f),
\end{equation*}
where  $f\in K\{X\}_{\infty}$. Especially, for
$T=1_{K\{X\}_{\infty}}=\emp$,  $\partial_{\emp}=\id$.

\goodbreak

More generally, for $g$ homogeneous, $g=\sum a_i T^i, a_i\in K$,
we also define $\partial_g(f):=\sum a_i
\partial_{T^i}(f)$.

 We analogously define $\partial_T(f)$ and $\partial_g(f)$ for
the binary case, where $f,g$ (or $T$) are elements of $K\{X\}$.

We call $\partial_T$ a generalized differential operator. It has
the following properties.

\end{defn}

\begin{proposition}\label{propdiffop}

\begin{itemize}
\item[(i)] For all monomials $S, T$ in $K\{X\}_{\infty}$,
\begin{equation*}
\partial_S\circ\partial_T=\partial_T\circ\partial_S
\end{equation*}
\item[(ii)]
For any $x_k\in X$, $\partial_{x_k}$ is the derivation
$\partial_k$.
\item[(iii)] For $T$ a monomial,
and $f_1,\ldots, f_p \in K\{X\}_{\infty}$,
\begin{equation*}
\partial_{T}(\vee^p(f_1\ldots f_p))=
\sum_{\vee^p(T^1.T^2\ldots
T^p)=T}\vee^p\bigl(\partial_{T^1}(f_1)\ldots\partial_{T^p}(f_p)\bigr),
\end{equation*}
where the sum is over all not necessarily non-empty trees $T^1,
T^2,\ldots, T^p$ such that $\vee^p(T^1.T^2\ldots T^p)=T$.
\newline
Especially
\begin{equation*}
\partial_{T}(f_1\cdot f_2)=
f_1\cdot\partial_{T}(f_2)+\partial_{T^1}(f_1)\cdot\partial_{T^2}(f_2)+
\partial_{T}(f_1)\cdot f_2
\end{equation*}
if $T=\vee^2(T^1.T^2)$ and $f_1,f_2 \in K\{X\}$ are binary.
\end{itemize}
\end{proposition}

\begin{proof}
\begin{itemize}
\item[1)]
Since $\Delta_a$ is cocommutative, we get the formula
\begin{equation*}
\begin{split}
{\Delta_a}^2(f)&=\sum_{S,T\in\PRTrx} S\otimes T\otimes
\partial_S\partial_T(f)\\
&=\sum_{S,T\in\PRTrx} S\otimes T\otimes
\partial_T\partial_S(f).\\
\end{split}
\end{equation*}
Hence (i) holds.
\item[2)]
In the sum
\begin{equation*}
\Delta_a(T)=\sum_{I\subseteq {\rm Le}(T)}{\rm red} (T\vert
I)\otimes {\rm red}(T\vert I^c).
\end{equation*}
of (\ref{propdelta}), ${\rm red} (T\vert I)=x_k$ exactly when
$I=\{\nu\}$, $\nu\in I_k$. Here $I_k\subseteq {\rm Le}(T)$ is the
subset given by all leaves with label $x_k$.
\newline
Since $\partial_k(T)$ is given by
\begin{equation*}
\sum_{\nu\in I_k} {\rm red}(T\vert\{\nu\}^c),
\end{equation*}
see Proposition (\ref{proppartial}), assertion (ii) follows.
\item[3)]
Since the co-addition $\Delta_a$ is an algebra homomorphism with
respect to the operations $\vee^{n}$,
 we get (iii).
\end{itemize}
\end{proof}

\begin{corollary}

Let $S, T$ in $K\{X\}_{\infty}$ be monomials, and let
 $\mu_S(T)$ be the number of subsets $I$ of ${\rm Le}(T)$ which
yield $S$ by reduced leaf-restriction onto $I$, i.e.\ for which
${\rm red}(T\vert I)=S$. Then the following formulas hold.
\begin{itemize}
\item[(i)] The generalized differential $\partial_S(T)$ is given by
a sum $\sum_{S'} \mu_{S'}(T)S'$ over trees $S'$ which are
complements, see (\ref{defleafsplit}), of $S$ in $T$.
\item[(ii)]
Furthermore, if $S$ is of the form $S=\vee^2(S^1.S^2)$, and $T$ is
of the form $T=\vee^2(T^1.T^2)$ (with non-trivial $S_i, T_i$) then
\begin{equation*}
\mu_{\vee^2(S^1.S^2)}(\vee^2(T^1.T^2))
 = \mu_{\vee^2(S^1.S^2)}(T^1)+
\mu_{\vee^2(S^1.S^2)}(T^2)+\mu_{S^1}(T^1) \mu_{S^2}(T^2).
\end{equation*}
\end{itemize}
\end{corollary}

\begin{proof}
Formula (i) is a reformulation of Definition (\ref{defdiffop}) in
the case where $f=T$ is a monomial. Formula (ii) is an easy
consequence of Proposition (\ref{propdiffop}) (iii).
\end{proof}

\begin{example}

Let $f$ be the monomial
\begin{equation*}
\vee^2\Bigl(\vee^2\bigl(x_1,
\vee^3(x_2,x_2,x_2))\bigr),\vee^2\bigl(\vee^3(x_2,x_2,x_1),x_2\bigr)\Bigr),
\end{equation*}
of Example (\ref{exderiveight}).

Then $\partial_{\vee^4(x_2,x_2,x_2,x_2)}(f)=0$. The generalized
derivative $\partial_{\vee^2(\vee^3(x_2,x_2,x_2),x_2)}(f)$ is
given by
$\vee^2\bigl(x_1,\partial_2(\vee^2(\vee^3(x_2,x_2,x_1),x_2)\bigr)$,
i.e.\ by

 $
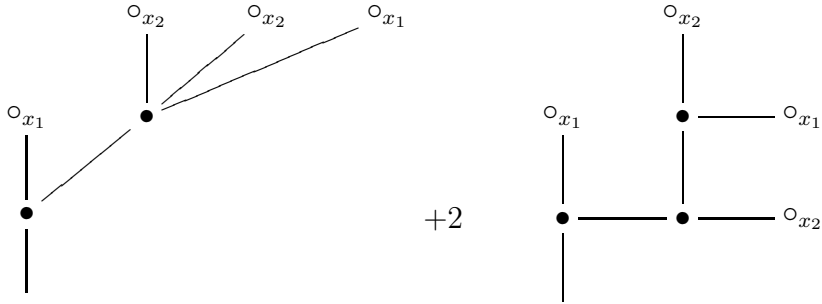
\xymatrix{ & \ar@{-}[d]\circ_{x_2} & \ar@{-}[dl]\circ_{x_2} &
 \ar@{-}[dll]\circ_{x_1} \\
 \circ_{x_1} \ar@{-}[d] & \bullet \ar@{-}[dl]\\
  \bullet\ar@{-}[d]  \\
   \\
  }
\xymatrix{ & & \ar@{-}[d]\circ_{x_2} \\
 & \circ_{x_1} \ar@{-}[d] & \bullet \ar@{-}[d] &\ar@{-}[l]\circ_{x_1}\\
  +2 &\bullet\ar@{-}[d] & \bullet\ar@{-}[l]& \ar@{-}[l]\circ_{x_2}\\
 &  \\
  }$

\bigskip

One notes that these terms also occur in
$\frac{1}{4!}\bigl(\partial_2\bigr)^4(f)$.

More generally, for $x_k\in X$, if $Z:=\{x_k\}$ and
$W^{(n)}:=\PRTrupnz$ is the set of reduced trees with $n$ leaves
that are all labeled by $x_k$, then:
\begin{equation*}
\sum_{S \in W^{(n)}}\partial_S=\frac{1}{n!}(\partial_k)^n.
\end{equation*}
This can be shown by a recursive argument, as
\begin{equation*}
\frac{1}{n!}\bigl(\partial_k\bigr)^n(\vee^p(T^1\ldots
T^p))=\frac{1}{n!}\sum_{i_1+\ldots+i_p=n} {n \choose
i_1,\ldots,i_p}\vee^p\bigl(\bigl(\partial_k\bigr)^{i_1}(T^1)\ldots
\bigl(\partial_k\bigr)^{i_p}(T^p)\bigr),
\end{equation*}
\begin{equation*}
\sum_{S \in W^{(n)}}\partial_S(\vee^p(T^1\ldots T^p))
=\sum_{i_1+\ldots+i_p=n}\ \ \sum_{S^1 \in W^{(i_1)},\ldots, S^p\in
W^{(i_p)}}\vee^p\bigl(\partial_{S^1}(T^1)\ldots\partial_{S^p}(T^p)\bigr).
\end{equation*}

\end{example}

\begin{lemma}\label{critpseudo}
For $f$ homogeneous of degree $n$, $f$ is primitive if and only if
$\partial_T(f)=0$ for all monomials $T\in K\{X\}_{\infty}$ with
$1\leq \deg T < \frac{n+1}{2}$.
\end{lemma}
\begin{proof}
By Definition (\ref{defdiffop}), $f$ is primitive if and only if
$\partial_T(f)=0$ for all monomials $T\in K\{X\}_{\infty}$.
 Using the cocommutativity of $\Delta_a$, the criterion follows.
\end{proof}

\end{section}

\begin{section}{Generalized Taylor expansions}

Given any element $f$ in a unitary $\Com$-algebra
$K[x_1,\ldots,x_m]/I$, the coefficients $a_{(n_1,\ldots,n_m)}\in
K$ in the formula
\begin{equation*}
f=\sum a_{(n_1,\ldots,n_m)} x_1^{n_1}\ldots x_m^{n_m},
\end{equation*}
can be successively determined by a process of calculating
 highest non-vanishing partial derivatives
$\partial_1^{n_1}\ldots \partial_m^{n_m}$ (and doing
subtractions).

A similar Taylor expansion also exists for noncommutative (cf.\
\cite{lidrecodim}, \cite{ligea}) and non-associative algebras
(cf.\ \cite{lighmag}, \cite{lidh}).

We are going to describe it for the $\Minf$-algebra
$K\{X\}_{\infty}$ or, more generally, for $\Minf$-algebras
$A=K\{X\}_{\infty}/I$. Here $X$ is the ordered set of variables
$x_1<x_2<\ldots,$ and $I$ is a multihomogeneous ideal which is
invariant under all partial derivatives $\partial_k$. We denote
the residue classes of $x_k\in X$ by the same symbols $x_k$, and
we denote the induced partial derivatives on $A$ also by
$\partial_k$.

\begin{lemma}

The subspace $A_0$ of $A$ given by all $f\in A$ such that
$\partial_k(f)=0$ is a $\Minf$-algebra.

\end{lemma}
\begin{proof}
The field $K$ is contained in $A_0$. For $f_1,\ldots,f_n\in A_0$,
$\partial_k\bigl(\vee^n(f_1.f_2\ldots f_n)\bigr) =\sum_{i=1}^n
\vee^n(\ldots.f_{i-1}.0.f_{i+1}\ldots)=0$, all $k$. Thus
$\vee^n(f_1.f_2\ldots f_n)\in A_0$.

\end{proof}

\begin{defn}

For $b\in A$, let $R_b: A\to A$, $a\mapsto a\cdot b:=\vee^2(a,b)$
be the binary right multiplication by $b$.

 For $f\in A$, $j\geq 0$, and
$x_k\in X$,
 let
\begin{equation*}
[f]_{\bullet} x_k^j=R_{x_k}^j(f)=(\ldots((f\underbrace{\cdot
x_{k})\cdot x_{k})\cdots x_{k})}_j.
\end{equation*}

If $x_1,x_2,\ldots, x_n\in X$, and ${\bf
j}=(j_1,j_2,\ldots,j_n)\in \N^n$, let
\begin{equation*}
[f]_{\bullet} x^{\bf j}=\bigl( R_{x_n}^{j_n}\circ
R_{x_{n-1}}^{j_{n-1}}\circ\cdots\circ R_{x_1}^{j_1}\bigr)(f).
\end{equation*}

\end{defn}

\begin{proposition}\label{proptaylor}
\begin{itemize}
\item[(i)]
There is a unique Taylor expansion of $f$ with respect to $x_k$,
\begin{equation*}
f=\sum_{j\geq 0} [a^{(k)}_{j}(f)]_{\bullet} x_k^{j},
\end{equation*}
 for every $f\in A=K\{X\}_{\infty}/I$,
with $\partial_k\bigl(a^{(k)}_{j}(f)\bigr)=0$ (for all $j$).
\newline
The elements $a^{(k)}_{j}$ are homogeneous of degree $n-j$ if $f$
is homogeneous of degree $n$.
\item[(ii)]
There is a unique Taylor expansion
\begin{equation*}
f=\sum_{\bf j} [a_{\bf j}(f)]_{\bullet} x^{\bf j}
\end{equation*}
 for every $f\in A$,
with $a_{\bf j}(f)\in A_0$.
\newline
The elements $a_{\bf j}(f)$ are homogeneous of degree $n-\sum_i
j_i$ if $f$ is homogeneous of degree $n$.
\newline
The map $\Phi: A\to A_0$ given by $f\mapsto a_{\bf 0}(f)$ is a
projector onto $A_0$.
\item[(iii)]
If $f$ depends only on the variables $x_1, x_2,\ldots, x_m$, and
$f=\sum_{\bf j} [a_{\bf j}]_{\bullet} x^{\bf j}$,
\newline
 then $f\cdot
x_m$ is given by
\begin{equation*}
\sum_{{\bf j}\in \N^{m-1}\times \N^*} [b_{\bf j}]_{\bullet} x^{\bf
j}, \text{ where, for all {\bf j}, } a_{\bf
j}=b_{(j_0,\ldots,j_{m-1},1+j_m)}.
\end{equation*}
Especially $f\cdot x_m$ is in the kernel of $\Phi$.
\end{itemize}
\end{proposition}

\begin{proof}
\begin{itemize}
\item[1)] For any $a^{(k)}_{j}\in A, j\geq 0,$ with $\partial_k(a^{(k)}_{j}(f))=0$,
\begin{equation*}
\partial_k \sum_{j\geq 0}
[a^{(k)}_{j}]_{\bullet} x_k^{j}=
 \sum_{j\geq 0}
[(j+1)a^{(k)}_{j+1}]_{\bullet} x_k^{j},
\end{equation*}
by the Leibniz rule.
\newline
For $f\in A$, let $n=n(f)$ be the maximal $n\in \N\cup\{-\infty\}$
such that $\partial_k(f)\neq 0$ (or $-\infty$ for $f=0$). Clearly
$n\leq \deg f$.
\newline
We prove the uniqueness and existence of the Taylor expansion (i)
by induction on $n$. In the case $n=0$, the Taylor expansion of
$f$ with respect to $x_k$ is given by $a^{(k)}_{0}(f)=f$. It is
unique, because $\partial_k(a\cdot x_k)\neq 0$ for all $a\neq 0$.
\newline
For $n\geq 1$, let $\partial_k f$ be given by the unique Taylor
expansion $f=\sum_{j\geq 0} [b^{(k)}_{j}]_{\bullet} x_k^{j}$,
where $b^{(k)}_{j}=a^{(k)}_{j}(\partial_k f)$.
\newline
 Let $\tilde f$ be given by
$\tilde f=\sum_{j\geq 1} [a^{(k)}_{j}]_{\bullet} x_k^{j}$, where
the elements $a^{(k)}_{j}$ are given by
$(j+1)a^{(k)}_{j+1}=b^{(k)}_{j}$.
Since $\partial_k(f-\tilde
f)=0$,
 $a^{(k)}_{0}:=f-\tilde f \in A_0$.
\newline
 Thus the desired Taylor expansion of $f$ is $
\sum_{j\geq 0} [a^{(k)}_{j}]_{\bullet} x_k^{j}$.
\newline
By construction, for homogeneous $f$, we get homogeneous elements
$a^{(k)}_{j}(f)$ of the asserted degree.
\item[2)]
Let $f$ depend only on the variables $x_1,\ldots, x_m$, and
let $f=\sum_{j\geq 0} [a^{(m)}_{j}(f)]_{\bullet} x_m^{j}$ be its
Taylor expansion with respect to $x_m$.
\newline
Continuing with the variables $x_{m-1},x_{m-2},\ldots, x_1$, we
get the desired Taylor expansion
\begin{equation*}
f=\sum_{\bf j} [\ldots [a_{\bf j}(f)]_{\bullet} x_1^{j_1}\ldots
]_{\bullet}x_m^{j_m}=\sum_{\bf j} [a_{\bf j}(f)]_{\bullet} x^{\bf
j}.
\end{equation*}
For $f\in A$, let now ${\bf n}={\bf n}(f)\in \N^m\cup\{-\infty\}$
be lexicographically maximal such that
$\partial_1^{n_1}\cdots\partial_m^{n_m}(f)\neq 0$, or $-\infty$
for $f=0$.
\newline
To show the uniqueness of the Taylor expansion (ii) we assume that
a nontrivial presentation $0=f=\sum_{{\bf j}\in J} [a_{\bf
j}(f)]_{\bullet} x^{\bf j}$ with nonzero $a_{\bf j}(f)\in A_0$ is
given, with $J\neq\emp$. On the one hand, ${\bf n}$ is the maximal
element in $J$. On the other hand, ${\bf n}=-\infty$, which is a
contradiction.
\newline
Again, the assertion about homogeneous elements follows by
construction.
\newline
For $a\in A_0$, $a=a_{\bf 0}(a)$ is the unique Taylor expansion,
thus $\Phi$ is surjective with $\Phi\circ\Phi=\Phi$.
\item[3)] Assertion (iii) follows easily from the construction of the Taylor expansion (ii).
\end{itemize}
\end{proof}

\begin{remark}
It should be noted that there are a lot of other possibilities to
define Taylor expansions for $A=K\{X\}_{\infty}/I$. The expansion
given above makes use of binary right multiplications, analogously
one can use left multiplications.

If we define
\begin{equation*}
[f]_{\infty} x^{{\bf
j}}:=\vee^{(1+j_1+\ldots+j_n)}(f,\underbrace{x_1,\ldots,x_1}_{j_1},\ldots,\underbrace{x_n,\ldots,x_n}_{j_n}),
\end{equation*}
we get a unique Taylor expansion
\begin{equation*}
f=\sum_{\bf j} [a_{\bf j}]_{\infty} x^{\bf j}
\end{equation*}
 for every $f\in A$,
with $a_{\bf j}=a_{\bf j}(f)\in A_0$.

The right Taylor expansion of (\ref{proptaylor}) and its
left-analogon have the advantage that they can be restricted to
$\Mag$-algebras as well. While, for $f\in K\{X\}\subseteq
K\{X\}_{\infty}$, the expression $f=\sum_{\bf j} [a_{\bf
j}(f)]_{\bullet} x^{\bf j}$ is defined in $K\{X\}$, this is not
the case for the non-binary expansion $f=\sum_{\bf j} [a_{\bf
j}]_{\infty} x^{\bf j}$.

\end{remark}

\begin{corollary}\label{corconstrep}
(compare \cite{lidh})

Let $X=\{x_1,x_2,\ldots,x_m\}$. Let $A$ be of the form
$K\{X\}_{\infty}/I$ or $A=K\{X\}/I$, where $I$ is a
multihomogeneous ideal invariant under all partial derivatives.
\begin{itemize}
\item[(i)]
The Hilbert series of the multigraded algebra $A$ and its
subalgebra $A_0$ of constants are related by
\begin{equation*}
\text{\Hilb}(A,t_1,t_2,\ldots,t_m)= \prod_{j=1}^m\frac{1}{1-t_j}
\text{\Hilb}(A_0,t_1,t_2,\ldots,t_m).
\end{equation*}
\item[(ii)]
Furthermore, if the ideal $I$ is invariant with respect to the
action of $GL_m(K)$, then
\begin{equation*}
A^{(n)}=\bigoplus_{j=0}^n A_0^{(j)}\otimes \underbrace{\rm
 \vbox{\offinterlineskip \hbox{\boxtext{\phantom{X
}}\boxtext{\phantom{X }}\boxtext{$\cdots$}\boxtext{\phantom{X
}}\boxtext{\phantom{X }}}
 } \rm}_{n-j}, \text{ for all } n,
\end{equation*}
as modules of the general linear group.
\end{itemize}
\end{corollary}

\begin{proof}

We consider the case where $A$ is of the form $K\{X\}_{\infty}/I$,
the binary case being completely analogous.

If $v_k, k\in J$ is a multihomogeneous basis
 of the vector space $A_0$,
then the Taylor expansion of Proposition (\ref{proptaylor}) shows
that $A$ has a basis
\begin{equation*}
[v_k]_{\bullet} x^{\bf j}=\bigl( R_{x_n}^{j_n}\circ
R_{x_{n-1}}^{j_{n-1}}\circ\cdots\circ R_{x_1}^{j_1}\bigr)(v_k),\
k\in J, {\bf j}=(j_1,j_2,\ldots,j_m)\in \N^m.
\end{equation*}

The vector space homomorphism $A\to A_0\otimes_K
K[x_1,\ldots,x_m]$ given by
\begin{equation*}
 [v_k]_{\bullet} x^{\bf j}
\end{equation*}
is an isomorphism of multigraded vector spaces.

Moreover it respects the action of the general linear group.

Since $K[x_1,\ldots,x_m]^{(n)}$ is the (module of the) trivial
representation with Young diagram
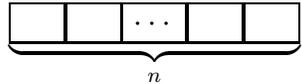
\begin{equation*}
\underbrace{\rm
 \vbox{\offinterlineskip \hbox{\boxtext{\phantom{X
}}\boxtext{\phantom{X }}\boxtext{$\cdots$}\boxtext{\phantom{X
}}\boxtext{\phantom{X }}}
 } \rm}_{n},
\end{equation*}
assertion (ii) follows.

For assertion (i), we note that the Hilbert series of
$K[x_1,\ldots,x_m]$ is
\begin{equation*}
\text{\Hilb}(K[x_1,\ldots,x_m],t_1,\ldots,t_m)
=\prod_{j=1}^m\frac{1}{1-t_j},
\end{equation*}
and that
 the Hilbert series of the tensor product is equal to the
product of the Hilbert series of the factors (cf.\ \cite{liufn}).

\end{proof}

\begin{remark}

To obtain the Taylor expansions (i) and (ii) of Proposition
(\ref{proptaylor}), we are going to use the following algorithm:

Let $f\in A$ depend only on the variables $x_1, x_2,\ldots, x_m$.

Let $n$ be maximal such that $(\partial_m)^n f\neq 0$. Then
$a^{(m)}_n(f)=\frac{1}{n!}(\partial_m)^n f$.

For $\tilde f := f - [a^{(m)}_n(f)]_{\bullet} x_m^j$, it holds
that $(\partial_m)^n \tilde f= 0$, and we can use $\tilde f$ to
obtain the coefficients $a^{(m)}_n(f), j<n$.

Analogously we continue with the variables
$x_{m-1},x_{m-2},\ldots, x_1$.

\end{remark}

\begin{remark}\label{rembasisonevar}

Since the free $\Mag$-algebra $K\{x_1\}$ in one variable $x_1$ has
the Hilbert series
\begin{equation*}
1+\frac{1-\sqrt{1-4t_1}}{2}=\sum_{n=0}^{\infty}c_n t_1^n, \text{
with } c_0:=1,
\end{equation*}
it follows from Corollary (\ref{corconstrep}) that the subalgebra
$K\{x_1\}_0$ has the Hilbert series
\begin{equation*}
(1-t_1)\Bigl(1+\frac{1-\sqrt{1-4t_1}}{2}\Bigr)=1+\sum_{n=1}^{\infty}(c_n-c_{n-1})
t_1^n=1+ t_1^3 + 3 t_1^4  + 9 t_1^5 + \ldots
\end{equation*}
(compare \cite{lidh}, \cite{lighmag}).

Hence the kernel of $\Phi: K\{x_1\}\to K\{x_1\}_0$ consists of all
$f\cdot x_1, f\in K\{x_1\}$, and a vector space basis of
$K\{x_1\}_0^{(n)}$, $n\geq 2$, is given by the set
\begin{equation*}
\{\Phi(T):\  T=T^1\cdot T^2\in W^{(n)}, T^1,T^2\neq 1, \text{ and
} T^2\neq x_1 \}.
\end{equation*}

\end{remark}

\end{section}

\begin{section}{Associators and the non-associative Jacobi relation}

We are going to use Taylor expansion to obtain the first primitive elements
in $K\{X\}_{\infty}$ and $K\{X\}$.

In degrees $n=0,1,2$, the homogeneous components of degree $n$
agree for the vector spaces $K\{X\}_{\infty}$, $K\{X\}$, and
$K\langle X \rangle$. We describe the constants, which are also
primitives, in the next example. Then, for $n\geq 3$, we treat the
two cases $\Mag$ and $\Minf$ separately.

\begin{example}
 Clearly all homogeneous elements of degree
$0$ are constants, all elements of degree $1$ are not constants.

The Taylor expansion of $x_{i_1}\cdot x_{i_2},
i_1\leq i_2$ is given by $[1]_{\bullet}
x_{i_1}x_{i_2}:=[1]_{\bullet} x^{\bf j}$, where $j_k=\begin{cases}
2& k=i_1=i_2\\ 1 &: i_1\neq i_2, k\in \{i_1,i_2\}\\0 &:
\text{else.}\\
\end{cases}$

The commutators $[x_{i_2},x_{i_1}]=x_{i_2}\cdot
x_{i_1}-x_{i_1}\cdot x_{i_2}$, $i_1< i_2$, are constants. Thus the
Taylor expansion of $f=x_{i_2}\cdot x_{i_1}$, $i_1<i_2$, is given
by
\begin{equation*}
[a_{\bf 0}(f)]_{\bullet} x^{\bf 0} + [1]_{\bullet}x_{i_1}x_{i_2}
\end{equation*}
with $a_{\bf 0}(f)=a_{\bf 0}(x_{i_2}\cdot
x_{i_1})=[x_{i_2},x_{i_1}]$.

In other words, looking at the representation of $\Sigma_2$ given
by multilinear homogeneous elements of degree 2, the constant part
corresponds to  the sign representation ${\rm
 \vbox{\offinterlineskip \hbox{\boxtext{\phantom{X
}}} \hbox{\boxtext{\phantom{X }}}
 } \rm}\ .$

The non-constant part (in degree 2) corresponds to the trivial
representation.

\end{example}

\begin{remark}

By Lemma (\ref{critpseudo}), the primitive elements of order $\geq
2$ form a subspace of the space $(K\{X\}_{\infty})_0$ of constants.

Furthermore it follows that all homogeneous elements of degree $\leq 3$
in $(K\{X\}_{\infty})_0$ are primitive.

\end{remark}

\begin{example}

Let $f=x_1\cdot(x_1\cdot x_1)\in K\{X\}\subset K\{X\}_{\infty}$.
 Then 3 is the maximal $n$ with
$(\partial_1)^n f\neq 0$. Thus we set
$a_3(f)=a_3^{(1)}(f)=\frac{1}{3!}(\partial_1)^3f=1$.

For $\tilde f:=f-[a_3(f)]_{\bullet}x_1^3=f-(x_1\cdot x_1)\cdot
x_1$ we repeat the step and find that $\tilde f\in A_0$. Thus the
Taylor expansion of $f$ is given by
\begin{equation*}
a_0(f)=x_1\cdot(x_1\cdot x_1)-(x_1\cdot x_1)\cdot x_1,\ \
a_3(f)=1, \text{ and } a_j(f)=0 \text{ else}.
\end{equation*}

Analogously the Taylor expansion of $f=\vee^3(x_1.x_1.x_1)\in K\{X\}_{\infty}$ is
given by
\begin{equation*}
a_0(f)=\vee^3(x_1.x_1.x_1)-(x_1\cdot x_1)\cdot x_1,\ \ a_3(f)=1,
\text{ and } a_j(f)=0 \text{ else}.
\end{equation*}

\end{example}

\begin{defn}

For $f,g,h\in K\{X\}_{\infty}$, we denote by
\begin{equation*}
(f,g,h)_b=(f\cdot g)\cdot h - f\cdot (g\cdot h)
\end{equation*}
the (binary) associator. The corresponding operation is denoted by
$(x_1,x_2,x_3)_b$.

We define as a ternary associator in $K\{X\}_{\infty}$ the element
\begin{equation*}
(f,g,h)_t=(f\cdot g)\cdot h - \vee^3(f,g,h).
\end{equation*}
\end{defn}

\begin{proposition}\label{propminfdegthree}

The spaces of binary and ternary operations of $\Prim\Minf$ can be
described as follows:

The space $\Prim\Minf(2)$ is given by $K\cdot
[x_1,x_2].$

The space $\Prim\Minf(3)$ has dimension 14 and is generated (as a
$\Sigma_3$-module) by the operations
\begin{equation*}
(x_1,x_2,x_3)_t, \ (x_1,x_2,x_3)_b, \text{ and } [[x_1,x_2],x_3]
\end{equation*}
which fulfill $[[x_2,x_1],x_3]=-[[x_1,x_2],x_3]$ and the equation
\begin{equation*}
\begin{split}
[[x_1,x_2],x_3]+[[x_3,x_1],x_2]
+[[x_2,x_3],x_1]
=\ \ &(x_1,x_2,x_3)_b-(x_2,x_1,x_3)_b\\
 +&(x_3,x_1,x_2)_b-(x_1,x_3,x_2)_b\\
 +&(x_2,x_3,x_1)_b-(x_3,x_2,x_1)_b
\\
\end{split}
\end{equation*}
(called the non-associative Jacobi relation, compare
\cite{lilosci}, and also \cite{lisu}).

\medskip

Moreover, in terms of representations of the symmetric groups,
\begin{equation*}
\Prim\Minf(3) =
 2\ {\rm \vbox{\offinterlineskip
\hbox{\boxtext{\phantom{X }}\boxtext{\phantom{X
}}\boxtext{\phantom{X }}}
 } \rm}\oplus\ 5
\ {\rm
 \vbox{\offinterlineskip \hbox{\boxtext{\phantom{X
}}\boxtext{\phantom{X }}} \hbox{\boxtext{\phantom{X }}}
 } \rm}
\oplus\ 2\  {\rm
 \vbox{\offinterlineskip \hbox{\boxtext{\phantom{X
}}} \hbox{\boxtext{\phantom{X }}} \hbox{\boxtext{\phantom{X }}} }
\rm}\ .
\end{equation*}

\medskip
\end{proposition}
\begin{proof}
\begin{itemize}
\item[1)]
We use that all constant homogeneous elements of degree 3 are
primitive, and apply Corollary (\ref{corconstrep}). Let
$A=K\{X\}_{\infty}$.
\newline
The subspace of
$A^{(3)}$ orthogonal to the space $A_0^{(3)}$ of primitive
elements is
\begin{equation*}
 {\rm \vbox{\offinterlineskip
\hbox{\boxtext{\phantom{X }}\boxtext{\phantom{X
}}\boxtext{\phantom{X }}}
 } \rm}\oplus\Biggl(\
{\rm
 \vbox{\offinterlineskip \hbox{\boxtext{\phantom{X
}}} \hbox{\boxtext{\phantom{X }}}
 } \rm}\otimes \boxtext{\phantom{X }}\ \Biggr),
\end{equation*}
where we use Young diagrams instead of $GL$-modules. The same
description applies to the corresponding $\Sigma_3$-modules (which
contain only multilinear elements).
\newline
 By Young's rule,
\begin{equation*}
{\rm
 \vbox{\offinterlineskip \hbox{\boxtext{\phantom{X
}}} \hbox{\boxtext{\phantom{X }}}
 } \rm}\otimes \boxtext{\phantom{X }}\
=\ {\rm
 \vbox{\offinterlineskip \hbox{\boxtext{\phantom{X
}}\boxtext{\phantom{X }}} \hbox{\boxtext{\phantom{X }}}
 } \rm}
\oplus\  {\rm
 \vbox{\offinterlineskip \hbox{\boxtext{\phantom{X
}}} \hbox{\boxtext{\phantom{X }}} \hbox{\boxtext{\phantom{X }}} }
\rm}\ .
\end{equation*}
\item[2)]
By comparison with three copies of the
regular representation, one immediately arrives at the description
of $\Prim\Minf(3)$ as a representation of
$\Sigma_3$ given above. Counting dimensions yields $14=2+5\cdot
2+2$.
\item[3)]
The commutators $[x_1,x_2]=-[x_2,x_1]$ and associators
$(x_1,x_2,x_3)_b$, $(x_1,x_2,x_3)_t$ are primitive.
\newline
If we define
$a(x_1,x_2,x_3):=[[x_1,x_2],x_3]-(x_1,x_2,x_3)_b+(x_2,x_1,x_3)_b$,
then
\begin{equation*}
a(x_1,x_2,x_3)=x_3\cdot(x_2\cdot x_1)-x_3\cdot(x_1\cdot
x_2)+x_1\cdot(x_2\cdot x_3)-x_2\cdot(x_1\cdot x_3).
\end{equation*}
The latter sum contains only right normed parenthesized words.
\newline
Thus the non-associative Jacobi relation
\begin{equation*}
0=a(x_1,x_2,x_3)+a(x_3,x_1,x_2)+a(x_2,x_3,x_1)
\end{equation*}
holds, simply because the classical Jacobi relation holds.
\end{itemize}
\end{proof}

\begin{proposition}\label{propmagdegthree}

The spaces of binary and ternary operations of $\Prim\Mag$ can be
described as follows:

The space $\Prim\Mag(2)$ is given by $K\cdot
[x_1,x_2].$

The space $\Prim\Mag(3)$ has dimension 8 and is generated by the operations
\begin{equation*}
\ (x_1,x_2,x_3)_b \text{ and } [[x_1,x_2],x_3]
\end{equation*}
which fulfill $[[x_2,x_1],x_3]=-[[x_1,x_2],x_3]$ and the
non-associative Jacobi relation.

\medskip

Moreover, in terms of representations of the symmetric groups,
\begin{equation*}
\text{\ \ }\Prim\Mag(3) =
 {\rm \vbox{\offinterlineskip
\hbox{\boxtext{\phantom{X }}\boxtext{\phantom{X
}}\boxtext{\phantom{X }}}
 } \rm}\oplus\ 3
\ {\rm
 \vbox{\offinterlineskip \hbox{\boxtext{\phantom{X
}}\boxtext{\phantom{X }}} \hbox{\boxtext{\phantom{X }}}
 } \rm}
\oplus\  {\rm
 \vbox{\offinterlineskip \hbox{\boxtext{\phantom{X
}}} \hbox{\boxtext{\phantom{X }}} \hbox{\boxtext{\phantom{X }}} }
\rm}\ .
\end{equation*}

\medskip
\end{proposition}

\begin{proof}

The proof is very similar to the proof of Proposition (\ref{propminfdegthree}).

The only difference is the following: In step 2 of the proof,
we have to
compare with two instead of three copies of the
regular representation to obtain $\Prim\Mag(3)$ as a representation of
$\Sigma_3$. Counting dimensions yields $8=1+3\cdot 2+1$.

\end{proof}

\begin{example}\label{exopdegfour}
Let $f$ be the tree monomial
\begin{equation*}
(x_1\cdot x_2)\cdot(x_3\cdot x_4)=
\vee^2(\vee^2(x_1.x_2).\vee^2(x_1.x_2))\in A=K\{X\}\subset K\{X\}_{\infty}.
\end{equation*}
Then the constant term of its Taylor expansion with respect to
$x_4$ is the term
\begin{equation*}
a_0^{(4)}(f)=(x_1\cdot x_2)\cdot(x_3\cdot x_4)-((x_1\cdot
x_2)\cdot x_3)\cdot x_4,
\end{equation*}
which is a constant for the derivatives $\partial_4$ and
$\partial_3$.

To get the constant term $a_{\bf 0}(f)\in A_0$, we compute that
\begin{equation*}
\partial_2\bigl(a_0^{(4)}(f)\bigr)
=x_1\cdot(x_3\cdot x_4)-(x_1\cdot x_3)\cdot x_4=-(x_1,x_3,x_4)_b,
\end{equation*}
\begin{equation*}
\begin{split}
\text{and\ \ \  }\partial_1&\bigl( (x_1\cdot x_2)\cdot(x_3\cdot
x_4)-((x_1\cdot x_2)\cdot x_3)\cdot x_4+(x_1,x_3,x_4)_b\cdot
x_2\bigr)\\ &=-(x_2,x_3,x_4)_b.\\
\end{split}
\end{equation*}
Thus $a_{\bf 0}(f)$ is given by
\begin{equation*}
\begin{split}
p&(x_1,x_2,x_3,x_4):=\\
 &(x_1\cdot x_2)\cdot(x_3\cdot x_4)-((x_1\cdot
x_2)\cdot x_3)\cdot x_4+(x_1,x_3,x_4)_b\cdot x_2
+(x_2,x_3,x_4)_b\cdot x_1
\\
&= (x_1,x_3,x_4)_b\cdot x_2 +(x_2,x_3,x_4)_b\cdot x_1 -((x_1\cdot
x_2),x_3,x_4)_b.
\\
\end{split}
\end{equation*}

Let $S$ be one of the six two-leaf tree monomials $x_i\cdot x_j$
with $i,j\in \underline{3},\ i\neq j$. Then it is immediately
checked that the generalized differential operator $\partial_S$
applied to $p(x_1,x_2,x_3,x_4)$ yields zero. (In fact, the
corresponding expression in associative variables is zero, and
$\partial_S$ evaluated on degree 4 elements cannot distinguish
between different bracketings). Therefore, by Lemma
(\ref{critpseudo}), $p(x_1,x_2,x_3,x_4)$ is primitive.

Analogously, using the Taylor expansion of left multiplications,
one can show that the following element is primitive:
\begin{equation*}
\begin{split}
q(x_1,x_2,x_3,x_4) & = \\ (x_1\cdot x_2)&\cdot(x_3\cdot x_4)
-x_1\cdot(x_2\cdot(x_3\cdot x_4))-x_3\cdot(x_1,x_2,x_4)_b
-x_4\cdot(x_1,x_2,x_3)_b.\\
\end{split}
\end{equation*}
By adding the primitive element $[x_4,(x_1,x_2,x_3)_b]$ to
$q(x_1,x_2,x_3,x_4)$ we get the non\-zero primitive element
$(x_1,x_2,x_3\cdot x_4)_b -x_3\cdot(x_1,x_2,x_4)_b
-(x_1,x_2,x_3)_b\cdot x_4.$

We also note that
\begin{equation*}
\begin{split}
(x_1,x_2\cdot x_3, x_4)_b &-x_2\cdot(x_1,x_3,x_4)_b
-(x_1,x_2,x_4)_b\cdot x_3
\\
=&\ q(x_1,x_2,x_3,x_4)+[(x_1,x_3,x_4)_b,x_2]+[(x_2,x_3,x_4)_b,x_1]
\\
 &-p(x_1,x_2,x_3,x_4)
-[(x_1,x_2,x_3)_b,x_4]-[(x_1,x_2,x_4)_b,x_3].\\
\end{split}
\end{equation*}

\end{example}

\begin{corollary}\label{cornotquadtern}
The operad $\Prim\Mag$ (and also the operad $\Prim\Minf$) cannot
be generated by quadratic and ternary operations.
\end{corollary}
\begin{proof}

By Example (\ref{exopdegfour}), $p(x_1,x_2,x_3,x_4)$ is an element
of $\Prim\Mag(4)$.

It cannot be generated with respect to insertion by the binary and
ternary operations given in Proposition (\ref{propminfdegthree}),
because $p(x_1,x_1,x_1,x_1)\neq 0$ and
$(x_1\cdot x_2)\cdot(x_3\cdot x_4)$ occurs in the support
of $p(x_1,x_2,x_3,x_4)$.

\end{proof}

\begin{remark}
It was stated in \cite{lighmag} that the element $q$ is primitive.
In \cite{lisu}, Shestakov and Umirbaev recursively construct
similar primitive elements
\begin{equation*}
P(x^{\bf i},y^{\bf j},z)=(x^{\bf i},y^{\bf j},z)_b -{\sum}'
(x^{\bf i})_{(1)}\cdot (y^{\bf j})_{(1)} \cdot P\bigl((x^{\bf
i})_{(2)}, (y^{\bf j})_{(2)},z\bigr),\ \ x^{\bf i}:=[1]_{\bullet}
x^{\bf i}
\end{equation*}
They also show that this gives a complete set of primitive
operations, and they pose the problem to give an intrinsic
characterization of the set of primitive elements in terms of
primitive operations.

\end{remark}

\end{section}

\begin{section}{The shuffle multiplication}

We are going to describe the graded dual $K\{X\}^{*g}_{\infty}$ of
the $\Minf$-Hopf algebra $K\{X\}_{\infty}$ equipped with
co-addition.

In analogy to the classical case, see Remark (\ref{remcocom}), the
commutative associative binary operation corresponding to
$\Delta_a$ is called planar tree shuffle multiplication (or planar
shuffle product, see \cite{ligepsp}).

\begin{proposition}\label{propdualshuffle}

The vector spaces $K\{X\}_{\infty}$ and $K\{X\}^{*g}_{\infty}$ can
be identified by mapping the basis given by tree monomials $T$ on
the corresponding dual basis elements $\delta_T$. Then the
commutative associative multiplication $\Delta_a^{*g}$ is given by
the binary operation
\begin{equation*}
\begin{split}
\sh: & K\{X\}_{\infty}\otimes K\{X\}_{\infty}\to K\{X\}_{\infty}\\
& T^1 \otimes T^2 \mapsto \sum_{\text{all shuffles } T \text{ of }
T^1 \text{ and } T^2} c_T\ T\text{ \ \ for tree monomials } T^1,
T^2,\\
\end{split}
\end{equation*}
where $c_T\geq 1$ is the number of subsets $I\subseteq {\rm
Le}(T)$ with
\begin{equation*}
{\rm red}(T\vert I)=T^1, {\rm red}(T\vert I^c)=T^2.
\end{equation*}
Especially $\ 1\sh 1=1$.

\end{proposition}

\begin{proof}

The space $K\{X\}_{\infty}\cong K\{X\}^{*g}_{\infty}$ together
with the binary operation $\sh$ given by $\Delta_a^{*g}$ is a
$\Com$-algebra with unit 1, see Lemma (\ref{lemdualhopf}) and
Remark (\ref{remdualhopf}).

The multiplication $\sh$ is characterized by the formula
\begin{equation*}
\langle T^1\sh T^2, T\rangle=\langle T^1\otimes
T^2,\Delta_a(T)\rangle,
\end{equation*}
where $\langle S^1\otimes S^2, T^1\otimes T^2\rangle = \langle
S^1, T^1\rangle\otimes\langle S^2,T^2\rangle$, and $\langle S,
T\rangle = \begin{cases} 1 : S=T \\ 0 : S\neq T\\ \end{cases}$,
for $T,T^1,T^2, S,S^1,S^2$ tree monomials.

By Proposition (\ref{propdelta}),
\begin{equation*}
\Delta_a(T)=\sum_{I\subseteq {\rm Le}(T)}{\rm red} (T\vert
I)\otimes {\rm red}(T\vert I^c).
\end{equation*}

Therefore the number $c_T$ of subsets $I\subseteq {\rm Le}(T)$
with
\begin{equation*}
{\rm red}(T\vert I)=T^1, {\rm red}(T\vert I^c)=T^2
\end{equation*}
is equal to $\langle T^1\sh T^2, T\rangle=\langle T^1\otimes
T^2,\Delta_a(T)\rangle$, and $c_T\geq 1$ if and only if $T$ is a
shuffle of $T^1$ and $T^2$, see Definition (\ref{defleafsplit}).

\end{proof}

\begin{remark}

Obviously, if $T^1$ has $k$ leaves labeled bijectively by
$x_1,\ldots,x_k$, and $T^2$ has $n-k$ leaves labeled bijectively
by $x_{k+1},\ldots, x_n$, then only $c_T=1$ can occur, thus
\begin{equation*}
T^1 \sh T^2 = \sum_{\text{all shuffles } T \text{ of } T^1 \text{
and } T^2} T
\end{equation*}
and we get a generalization of the well-known shuffle
multiplication of permutations $\sh:K\Sigma_k\times
K\Sigma_{n-k}\to K\Sigma_n$.
\end{remark}

\begin{example}\label{exshalltrees}

For $T$ in $\PRTree^n$ a (non-labeled) planar reduced tree with
$n$ leaves, let $T^{12\ldots n}$ be the corresponding tree in
$\PRTrx$ with first leaf labeled by $x_1$, second by $x_2$, and so
on.

For $\sigma\in\Sigma_n$, let $\sigma T^{12\ldots n}$ be the tree
which is obtained by permuting the labels of the leaves, i.e.\
with first leaf labeled by $x_{\sigma(1)}$, second by
$x_{\sigma(2)}$, and so on.

Then the product $x_1\sh x_2 \sh \ldots \sh x_n$ is given by
\begin{equation*}
\sum_{T\in\PRTree^n}\ \ \sum_{\sigma\in\Sigma_n}\ \ \sigma
T^{12\ldots n}
\end{equation*}
as can be shown by induction (the case $n=1$ being trivial) using
the fact that every term of $x_1\sh x_2 \sh \ldots \sh x_n$ occurs
in a unique way as a shuffle of a term of $x_1\sh x_2 \sh \ldots
x_{n-1}$ and $x_n$.

\end{example}

\begin{example}

Let $f=[x_2, x_1]$. Then, by a similar argument as in
(\ref{exshalltrees}), we get
\begin{equation*}
f\sh x_3=\sum_{T\in\PRTree^3}\ T^{321}+T^{231}+T^{213}
-T^{123}-T^{132}-T^{312}.
\end{equation*}

If $f=\vee^3(x_1.x_2.x_3)$, then $f\sh x_4$ is given by
\begin{equation*}
\begin{split}
&\vee^4(x_1.x_2.x_3.x_4)+\vee^4(x_1.x_2.x_4.x_3)+\vee^4(x_1.x_4.x_2.x_3)+
\vee^4(x_4.x_1.x_2.x_3)
\\
 &+ \vee^3(\vee^2(x_1.x_4).x_2.x_3)+\vee^3(x_1.\vee^2(x_2.x_4).x_3)
+\vee^3(x_1.x_2.\vee^2(x_3.x_4)) \\
 &+ \vee^3(\vee^2(x_4.x_1).x_2.x_3)+\vee^3(x_1.\vee^2(x_4.x_2).x_3)
+\vee^3(x_1.x_2.\vee^2(x_4.x_3)) \\
 &+ \vee^2(\vee^3(x_1.x_2.x_3).x_4)+\vee^2(x_4.\vee^3(x_1.x_2.x_3)).\\
\end{split}
\end{equation*}

In this sum, the shuffles $\vee^3(x_1.\vee^2(x_2.x_4).x_3)$ and
$\vee^2(x_4.\vee^3(x_1.x_2.x_3))$ correspond to the trees

$\xymatrix{ {\ \ \bullet_{x_1}}\ar@{-}[dd] & {\ \
\bullet_{x_2}}\ar@{-}[d] & {\ \bullet_{x_4}}\ar@{-}[dl] & {
\bullet_{x_3}}\ar@{-}[ddlll]\\
  &{\circ}\ar@{-}[dl]\\
\circ\ar@{-}[d]& & & \text{ and }\\
\\}
\xymatrix{ {\ \ \bullet_{x_4}}\ar@{-}[dd] & {\ \
\bullet_{x_1}}\ar@{-}[d] & {\ \bullet_{x_2}}\ar@{-}[dl] & {
\bullet_{x_3}}\ar@{-}[dll]\\ & {\circ}\ar@{-}[dl]\\
\circ\ar@{-}[d]\\
\\}
$

\end{example}

\begin{proposition}\label{propdualshmag}

The graded dual of the $\Mag$-Hopf algebra $K\{X\}$ equipped with
co-addition is the quotient of $(K\{X\}_{\infty},\sh)$ with
respect to the projection
\begin{equation*}
\pi: K\{X\}_{\infty}\to K\{X\}
\end{equation*}
\begin{equation*}
 \text{given on monomials } T \text{
by }T\mapsto
\begin{cases}T &: T \text{binary}\\ 0 &: \text{else.}\\
\end{cases}
\end{equation*}

\end{proposition}

\begin{proof}

Under the identification $K\{X\}_{\infty}^{*g}=K\{X\}_{\infty}$ of
Proposition (\ref{propdualshuffle}), the injection $\iota:
K\{X\}\to K\{X\}_{\infty}$ corresponds to the projection $\pi$
defined above.

The injection $\iota$ can be viewed as a morphism of $\Minf$-Hopf
algebras.

It follows that $(K\{X\}_{\infty}^{*g}/\ker\pi, \sh)$, where $\sh$
becomes an operation on binary trees, is the graded dual of
$(K\{X\},\Delta_a)$.

\end{proof}

\end{section}

\begin{section}{An analogon of Poincar\'e-Birkhoff-Witt}

As a vector space (in fact as a coalgebra) the free $\As$-algebra
on $V$ is isomorphic to the free $\Com$-algebra generated
by all Lie polynomials, i.e.\ by the primitive elements.
This is the Poincar\'e-Birkhoff-Witt theorem.

To describe the operads $\Prim\Mag$ and $\Prim\Minf$, we are going
to use an analogon of this theorem.

First, we need to describe the primitive elements as irreducible
elements with respect to the shuffle multiplication, and we
also need a description of the operation $(\vee^2)^*:K\{X\}_{\infty}\to
K\{X\}_{\infty}\otimes K\{X\}_{\infty}$.

\begin{proposition}\label{propshuffortho}

Let $f\in K\{X\}_{\infty}^{(n)}$ be homogeneous of degree $n\geq 1$.
\begin{itemize}
\item[(i)]
If $f$ is primitive, then $\langle f, g_1\sh g_2
\rangle=0$ for all homogeneous $g_1, g_2\in K\{X\}_{\infty}$ of degree $\geq 1$.
\item[(ii)]
If $S$ is a tree monomial of degree $1\leq k \leq n-1$, then
$\partial_S(f)=0$ if and only if $\langle f, S\sh g \rangle=0$ for
all $g\in K\{X\}_{\infty}^{(n-k)}$.
\item[(iii)]
If the homogeneous element $f$ of degree $n$ is
 orthogonal to all shuffle products $S\sh T$, for
$S,T$ tree monomials of degrees $k, n-k$ with $1\leq k <
\frac{n+1}{2}$, then $f$ is primitive.
\end{itemize}
The analogous assertions hold in the binary case $f\in K\{X\}^{(n)}$.
\end{proposition}

\begin{proof}
We consider elements of $K\{X\}_{\infty}$
(the proof for $K\{X\}$ is almost literally the same).
By definition,
\begin{equation*}
\langle \Delta_a(f),  g_1\otimes g_2 \rangle=\langle f, g_1\sh g_2
\rangle.
\end{equation*}
Clearly, if $f$ is primitive,
$\langle f, g_1\sh g_2\rangle=\langle f,g_1\rangle\langle 1,g_2 \rangle
+\langle 1,g_1\rangle\langle f,g_2 \rangle=0$.

If $\partial_S(f)\neq 0$ for some tree monomial $S$ of degree
$1\leq k \leq n-1$, then by definition
\begin{equation*}
0\neq \langle \Delta_a(f),  S\otimes g \rangle=\langle f, S\sh g
\rangle \text{ for some } g\in K\{X\}_{\infty}^{(n-k)}.
\end{equation*}
If $f$ is orthogonal to all shuffle products $S\sh S'$, $S, S'$
tree monomials with $\deg S=k \geq 1$, then
\begin{equation*}
\Delta_a(f)=\sum_{T\in\PRTrx-\{S\}} T\otimes
\partial_T(f)
\end{equation*}
and $\partial_S(f)=0$.

Then assertion (iii) follows by Lemma (\ref{critpseudo}).

\end{proof}

\begin{defn}

Let $A=K\{X\}_{\infty}$ or $A=K\{X\}$.

We denote by
$\nabla_k:A\to A^{\otimes k}, k\geq 2$
(for $A=K\{X\}$, $k=2$ only)
the maps given
by
$(\vee^k)^*:A^{*g}\to (A^{*g})^{\otimes k}$,
see (\ref{lemdualhopf}).

\end{defn}

\begin{lemma}\label{lemcomag}

We consider the unitary $\Com$-algebras $A=(K\{X\}_{\infty},\sh)$ and
 $A=(K\{X\},\sh)$, and we denote by $\bar A$ the augmentation ideal.
 \begin{itemize}
\item[(i)] The $K$-linear map $\nabla_2:A \to A\otimes A$ is a morphism
of $\Com$-algebras.
\item[(ii)] The $K$-linear map $\nabla_2$ is non-coassociative. If $T$ is
an admissibly labeled tree, then
\begin{equation*}
\nabla_2(T)=\begin{cases} T\otimes 1+1\otimes T +T^1\otimes T^2&:
T=\vee^2(T^1.T^2)\\ T\otimes 1+1\otimes T &: \text{ar}_\rho\neq 2\
\ (\rho \text{ the root}).\\ \end{cases}
\end{equation*}
\item[(iii)] If $\nabla'_2(f):=\nabla_2(f)-f\otimes 1-1\otimes f$, then
$\nabla'_2(f)\in \bar A\otimes \bar A$.
\end{itemize}
In other words,
these unitary $\Com$-algebras are co-${\bf D}$ objects,
where ${\bf D}$ is the category of unitary magmas.
\end{lemma}

\begin{proof}

\begin{itemize}
\item[1)] The chosen unit actions on $\Minf$ (and also on $\Mag$) and $\Com$
have the property that the operations on the tensor product $A\otimes A$
are defined component-wise.
We have to check that
\begin{equation*}
\nabla_2\circ\sh=(\sh\otimes\sh)\circ\tau_2\circ(\nabla_2\otimes\nabla_2).
\end{equation*}
By looking at the graded dual, this is equivalent to the equation
\begin{equation*}
\Delta\circ\vee^2=(\vee^2\otimes\vee^2)\circ\tau_2\circ(\Delta\otimes\Delta).
\end{equation*}
The latter equation is fulfilled, because
$(\vee^2\otimes\vee^2)\circ\tau_2$ is $\vee^2_{A\otimes A}$ and
$\Delta$ is a morphism of $\Mag$-algebras.
\item[2)]Since $\nabla_2$ is determined by the equation
\begin{equation*}
\langle \vee^2(f_1.f_2), g\rangle=
\langle f_1\otimes f_2,  \nabla_2(g)\rangle
\end{equation*}
we conclude that, for $T$ an admissibly labeled tree,
\begin{equation*}
\nabla_2(T)=T\otimes 1+1\otimes T +\sum_{T=\vee^2(T^1.T^2)}T^1\otimes T^2.
\end{equation*}
It also follows that $\nabla_2$ is not coassociative.
We note that for $f, g\in A$,
\begin{equation*}
\langle f, g\rangle=
\langle \vee^2(1.f), g\rangle=
\langle 1\otimes f,  \nabla_2(g)\rangle=
\langle f\otimes 1,  \nabla_2(g)\rangle
\end{equation*}
Thus assertion (iii) follows.
\item[3)]
The categorical coproduct for $\Com$-algebras is the tensor
product $\otimes$. By (i) and (iii), $\nabla_2$ provides the
unitary $\Com$-algebras $A=(K\{X\}_{\infty},\sh)$ and
$A=(K\{X\},\sh)$ with the structure of a co-magma object, see
Section (\ref{secthopftheory}), with counit given by the
augmentation map.
\end{itemize}
\end{proof}

\begin{remark}

The vector space $K\{X\}_{\infty}$ (or $K\{X\}$) equipped with
$\nabla_2$ is a non-associative coalgebra in the sense of
(\cite{liac}, \cite{ligr}).

\end{remark}

\begin{thm}\label{thmpbw}
Let either $A=K\{X\}_{\infty}$ or $A=K\{X\}$, and let
$W=\Prim A$ be its space of primitive elements, graded
by $W^{(k)}=\Prim A^{(k)}$.
\begin{itemize}
\item[(i)]
The unitary $\Com$-algebra $(A,\sh)$ is freely generated by $W=\Prim A$ with
respect to the shuffle multiplication.
\item[(ii)] The generating series of  $\Prim\Minf$
is the logarithm $\log(1+t)$ of the generating series
of $\Minf$. Similarly, the generating series of  $\Prim\Mag$
is the logarithm of the generating series
of $\Mag$.
\item[(iii)]
Let $n\geq 2$, and let $B^{(n)}$ be the orthogonal complement
of $\Prim A^{(n)}$ in the vector space $A^{(n)}$.
If $f_1,\ldots,f_r$ is
a basis (consisting of homogeneous elements) of $\oplus_{k=1}^{n-1}\Prim A^{(k)}$,
then the elements
\begin{equation*}
f_{i_1}\sh f_{i_2}\sh\ldots\sh f_{i_k}, 1\leq i_1\leq i_2\leq
\ldots \leq i_k\leq r, \text{ with } \sum_{j=1}^k\deg f_{i_j}=n
\end{equation*}
form a basis of $B^{(n)}$.
\end{itemize}
\end{thm}

\begin{proof}

\begin{itemize}
\item[1)]
We have shown in Proposition (\ref{propshuffortho}) that the
homogeneous primitive elements of $A$ are
 also the homogeneous $\sh$-irreducible elements of $A$,
i.e.\  the elements $f$ not of the form $g_1\sh g_2$ (for $f$, $g_1, g_2$
homogeneous of degree $\geq 1$).
Thus the unitary $\Com$-algebra morphism $\pi: K1\oplus F_{\Com}(W)\to (A,\sh)$
is surjective, and for every proper subspace $U$ of $W$,
$K1\oplus F_{\Com}(U)$ cannot be isomorphic to $(A,\sh)$.
\item[2)] We have to show that $\pi$ is injective, i.e.\ that
 $(A,\sh)$ is a free $\Com$-algebra.
\newline
By Lemma (\ref{lemcomag}), there exists a morphism
            $\nabla_2:A \to A\otimes A$ of $\Com$-algebras,
which provides the unitary $\Com$-algebras $(A,\sh)$ with the structure of a
co-${\bf D}$ object in the category of $\Com$-algebras, where
${\bf D}$ is the category of unitary magmas.
\newline
Over a field $K$ of characteristic 0,
all connected (i.e.\ $A^{(0)}=K$) $\Com$-algebras that are
equipped with the structure of a unital co-magma are free.
This is the Leray theorem, see \cite{lioud}.
Thus $\pi$ is an isomorphism.
\item[3)] We recall from Example (\ref{exascom}) that
the generating series of the operad $\Com$ is $\exp(t)-1$.
Combining this with assertion (i),
we get that
\begin{equation*}
\begin{split}
f^{\Mag}&=\exp\bigl(f^{\Prim\Mag}\bigr)-1,\\
f^{\Minf}&=\exp\bigl(f^{\Prim\Minf}\bigr)-1,\\
\end{split}
\end{equation*}
 Thus
assertion (ii) follows.
\item[4)]
Since the space $\Prim A^{(n)}$ is orthogonal (with respect to
$\langle, \rangle$) to the shuffle products in $A^{(n)}$, see
Proposition (\ref{propshuffortho}), assertion (iii) follows from
assertion (i).
\end{itemize}

\end{proof}

\begin{thm}\label{thmpbwdual}
Let $A=K\{X\}_{\infty}$ or $A=K\{X\}$, and let
$W=\Prim A$ be its space of primitive elements.

Then $(A,\Delta_a)$ is a cofree (co-)nilpotent $\Com$-coalgebra,
in the sense of (\ref{remcofree}) and (\ref{excomcofree}),
(co-)generated by $W$.

\end{thm}

\begin{proof}

We introduced $(A,\sh)$ as the graded dual of $(A,\Delta_a)$. Thus
it is allowed to dualize, and the statement of the theorem is
exactly
 dual to the one of Theorem (\ref{thmpbw})(i).

\end{proof}

\begin{remark}

Recall that $\Prim\As=\Lie$ and $\Prim\Com={\rm Vect}$
with generating series $f^{\Lie}(t)=\sum_{n=1}^{\infty}\frac{(n-1)!}{n!}t^n=
-\log(1-t)$ and $f^{\rm Vect}(t)=t$.

It is easy to check that
$f^{\As}(t)=\frac{1}{1-t}-1=\exp(f^{\Lie}(t))-1$ and
$f^{\Com}(t)=\exp(f^{\rm Vect}(t))-1$, which are the analogues of
Theorem (\ref{thmpbw})(ii) for the operads $\As$ and $\Com$.

Theorem (\ref{thmpbw})(iii) shows that one can define orthogonal projectors
\begin{equation*}
e_n^{(1)},e_n^{(2)},\ldots,e_n^{(n)}
\end{equation*}
     which are similar
to the Eulerian idempotents (see \cite{lilohaus}).
The idempotent $e_n^{(i)}$ projects elements
of $A^{(n)}$ to the subspace of $i$-factor shuffle products of primitive elements.

To obtain an explicit form of Eulerian idempotents for $\Mag$ and
$\Minf$, one may use the non-associative (planar) exponential
series of \cite{ligeplaexp} (which maps primitive elements to
group-like elements).

\end{remark}

\begin{remark}\label{remprimcmg}

We are going to give explicit formulas for the generating series
of $\Prim\Mag$ and $\Prim\Minf$ in the next section.

If $\Mag$ and $\Minf$ are replaced by   $\Cmg$ and the commutative version
$\Cminf$ of $\Minf$,
Theorem (\ref{thmpbw}) remains true, with the exception that we do not get
explicit formulas for the generating series of
the operads of primitives, see Remark (\ref{remsequences}).

For $\mathcal P=\Prim\Cmg$, we get that ${\mathcal P}(2)=0$ and
${\mathcal P}(3)$ is 2-dimensional, generated by the
operation $(x_1,x_2,x_3)_b=(x_1\cdot x_2)\cdot x_3 - x_1\cdot (x_2\cdot x_3)
=-(x_3,x_2,x_1)_b$ subject to the relation
$(x_1,x_2,x_3)_b+(x_2,x_3,x_1)_b+(x_3,x_1,x_2)_b=0$.

Since the coefficient of $t^4$ in $\exp(t+\frac{2t^3}{3!})$ is
$\frac{1}{4!}(1+8)$, we get that $\dim{\mathcal
P}(4)=\dim\Cmg(4)-9=15-9=6$.

\end{remark}

\end{section}

\begin{section}{The generating series  and
representations}

In the following, we describe the series
\begin{equation*}
f^{\Prim\Mag} \text{ and } f^{\Prim\Minf}
\end{equation*}
explicitly. Moreover, we compute
$\Prim\Minf(4)$ and $\Prim\Mag(4)$ as representations of
$\Sigma_4$.

We start with the operad $\Prim\Mag$ and obtain from Theorem
(\ref{thmpbw}) the following corollary:

\begin{corollary}\label{corlogcatalan}

The dimension of ${\Prim\Mag}(n)$ is given by
\begin{equation*}
\dim{\Prim\Mag}(n)=(n-1)! c'_n
\end{equation*}
where $c'_n$ is the $n$-th $\log$-Catalan number, see Example
(\ref{exlogcatalan}).

Thus the sequence ${\Prim\Mag}(n), n\geq 1,$ starts with
\begin{equation*}
\begin{split}
0!\cdot & 1,\ 1!\cdot 1,\ 2!\cdot 4,\ 3!\cdot 13,\ 4!\cdot 46,\
5!\cdot 166,\ 6!\cdot 610,\ \ \ldots\\ \text{ i.e.\ }& \ 1,\ \ \
1,\ \ \ \ 8,\ \ \ \ 78,\ \ \ \ 1104,\ \ \
 19920,\ \ 439200,\ \ \ldots
\end{split}
\end{equation*}

\end{corollary}

\begin{proof}

By Theorem (\ref{thmpbw})(ii), we get that the generating series
$f^{\Prim\Mag}(t)$ is given by
\begin{equation*}
\log\bigl(\frac{3-\sqrt{1-4t}}{2}\bigr).
\end{equation*}

Thus, for $\sum_{n=1}^{\infty} \frac{b_n}{n!} \
t^n:=f^{\Prim\Mag}(t)$,
\begin{equation*}
\begin{split}
\sum_{n=1}^{\infty} \frac{b_n}{(n-1)!} t^n
&=t\cdot\frac{\partial}{\partial t}f^{\Prim\Mag}(t)\\
&=\frac{2t}{3\sqrt{1-4t}-1+4t}
\ =\ \sum_{n=1}^{\infty} c'_n t^n,\\
\end{split}
\end{equation*}
where $c'_n$ is the $n$-th $\log$-Catalan number of Example
(\ref{exlogcatalan}).

\end{proof}

\begin{remark}

We have already seen in Proposition (\ref{propmagdegthree}), that
the representation of $\Sigma_n$ given by $\Prim\Mag(n)$ is not
given by copies of the $(n-1)!$-dimensional representation
$\Lie(n)$. For $n=4$, we obtain the following description:

\end{remark}

\begin{proposition}\label{propmagdegfour}

In terms of representations of the symmetric groups,
\begin{equation*}
\begin{split}
 &\Prim\Mag(4)=
\\
 & 3 {\rm \vbox{\offinterlineskip \hbox{
\boxtext{\phantom{X }}\boxtext{\phantom{X }}\boxtext{\phantom{X
}}\boxtext{\phantom{X }}}
 } \rm}\ \oplus\ 10\
{\rm
 \vbox{\offinterlineskip \hbox{\boxtext{\phantom{X
}}\boxtext{\phantom{X }}\boxtext{\phantom{X }}}
\hbox{\boxtext{\phantom{X }}}
 } \rm}
\oplus\ 6\ {\rm
 \vbox{\offinterlineskip \hbox{\boxtext{\phantom{X
}}\boxtext{\phantom{X }}} \hbox{\boxtext{\phantom{X
}}\boxtext{\phantom{X }}}
 } \rm}
 \oplus\ 10\ {\rm
 \vbox{\offinterlineskip \hbox{\boxtext{\phantom{X
}}\boxtext{\phantom{X }}} \hbox{\boxtext{\phantom{X }}}
\hbox{\boxtext{\phantom{X }}} } \rm}
 \oplus\ 3\ {\rm
 \vbox{\offinterlineskip \hbox{\boxtext{\phantom{X
}}} \hbox{\boxtext{\phantom{X }}} \hbox{\boxtext{\phantom{X }}}
\hbox{\boxtext{\phantom{X }}} } \rm} \ .\\
\end{split}
\end{equation*}
\bigskip

\end{proposition}
\begin{proof}
Let $A=K\{X\}$. We determine the space $A_0^{(4)}$ of degree 4
homogeneous constant elements (as a representation of the general
linear or symmetric group) first.

We recall that $A_0^{(0)}=K, \ A_0^{(1)}=0,$ and
\begin{equation*}
A_0^{(2)}\otimes{\rm
 \vbox{\offinterlineskip \hbox{\boxtext{\phantom{X
}}\boxtext{\phantom{X }}}
 } \rm}
\ = {\rm
 \vbox{\offinterlineskip \hbox{\boxtext{\phantom{X
}}} \hbox{\boxtext{\phantom{X }}}
 } \rm}\otimes\
{\rm
 \vbox{\offinterlineskip \hbox{\boxtext{\phantom{X
}}\boxtext{\phantom{X }}}
 } \rm}
\ = \ {\rm
 \vbox{\offinterlineskip \hbox{\boxtext{\phantom{X }}\boxtext{\phantom{X
}}\boxtext{\phantom{X }}} \hbox{\boxtext{\phantom{X }}}
 } \rm}
 \oplus\ {\rm
 \vbox{\offinterlineskip \hbox{\boxtext{\phantom{X
}}\boxtext{\phantom{X }}} \hbox{\boxtext{\phantom{X }}}
\hbox{\boxtext{\phantom{X }}} } \rm}\ .
\end{equation*}
Furthermore, by Proposition (\ref{propmagdegthree}),
\begin{equation*}
A_0^{(3)}\otimes{\rm
 \vbox{\offinterlineskip \hbox{\boxtext{\phantom{X
}}}
 } \rm}\ \text{ is equal to }
\end{equation*}
\begin{equation*}
{\rm \vbox{\offinterlineskip \hbox{ \boxtext{\phantom{X
}}\boxtext{\phantom{X }}\boxtext{\phantom{X }}\boxtext{\phantom{X
}}}
 } \rm}\ \oplus\ 4\
{\rm
 \vbox{\offinterlineskip \hbox{\boxtext{\phantom{X
}}\boxtext{\phantom{X }}\boxtext{\phantom{X }}}
\hbox{\boxtext{\phantom{X }}}
 } \rm}
\oplus\ 3\ {\rm
 \vbox{\offinterlineskip \hbox{\boxtext{\phantom{X
}}\boxtext{\phantom{X }}} \hbox{\boxtext{\phantom{X
}}\boxtext{\phantom{X }}}
 } \rm}
 \oplus\ 4\ {\rm
 \vbox{\offinterlineskip \hbox{\boxtext{\phantom{X
}}\boxtext{\phantom{X }}} \hbox{\boxtext{\phantom{X }}}
\hbox{\boxtext{\phantom{X }}} } \rm}
 \oplus\ \ {\rm
 \vbox{\offinterlineskip \hbox{\boxtext{\phantom{X
}}} \hbox{\boxtext{\phantom{X }}} \hbox{\boxtext{\phantom{X }}}
\hbox{\boxtext{\phantom{X }}} } \rm} \ .
\end{equation*}

\bigskip

Next we have to compare
\begin{equation*}
\bigoplus_{j=0}^3 A_0^{(j)}\otimes \underbrace{\rm
 \vbox{\offinterlineskip \hbox{\boxtext{\phantom{X
}}\boxtext{$\cdots$}\boxtext{\phantom{X }}}
 } \rm}_{4-j}
\end{equation*}
with $c_4=5$ copies of the regular representation
\begin{equation*}
  {\rm \vbox{\offinterlineskip
\hbox{ \boxtext{\phantom{X }}\boxtext{\phantom{X
}}\boxtext{\phantom{X }}\boxtext{\phantom{X }}}
 } \rm}\ \oplus\ 3\
{\rm
 \vbox{\offinterlineskip \hbox{\boxtext{\phantom{X
}}\boxtext{\phantom{X }}\boxtext{\phantom{X }}}
\hbox{\boxtext{\phantom{X }}}
 } \rm}
\oplus\ 2\ {\rm
 \vbox{\offinterlineskip \hbox{\boxtext{\phantom{X
}}\boxtext{\phantom{X }}} \hbox{\boxtext{\phantom{X
}}\boxtext{\phantom{X }}}
 } \rm}
 \oplus\ 3\ {\rm
 \vbox{\offinterlineskip \hbox{\boxtext{\phantom{X
}}\boxtext{\phantom{X }}} \hbox{\boxtext{\phantom{X }}}
\hbox{\boxtext{\phantom{X }}} } \rm}
 \oplus\ \ {\rm
 \vbox{\offinterlineskip \hbox{\boxtext{\phantom{X
}}} \hbox{\boxtext{\phantom{X }}} \hbox{\boxtext{\phantom{X }}}
\hbox{\boxtext{\phantom{X }}} } \rm} \ .
\end{equation*}

\bigskip

It follows that $A_0^{(4)}$ is given by
\begin{equation*}
 3 {\rm \vbox{\offinterlineskip
\hbox{ \boxtext{\phantom{X }}\boxtext{\phantom{X
}}\boxtext{\phantom{X }}\boxtext{\phantom{X }}}
 } \rm}\ \oplus\ 10\
{\rm
 \vbox{\offinterlineskip \hbox{\boxtext{\phantom{X
}}\boxtext{\phantom{X }}\boxtext{\phantom{X }}}
\hbox{\boxtext{\phantom{X }}}
 } \rm}
\oplus\ 7\ {\rm
 \vbox{\offinterlineskip \hbox{\boxtext{\phantom{X
}}\boxtext{\phantom{X }}} \hbox{\boxtext{\phantom{X
}}\boxtext{\phantom{X }}}
 } \rm}
 \oplus\ 10\ {\rm
 \vbox{\offinterlineskip \hbox{\boxtext{\phantom{X
}}\boxtext{\phantom{X }}} \hbox{\boxtext{\phantom{X }}}
\hbox{\boxtext{\phantom{X }}} } \rm}
 \oplus\ 4\ {\rm
 \vbox{\offinterlineskip \hbox{\boxtext{\phantom{X
}}} \hbox{\boxtext{\phantom{X }}} \hbox{\boxtext{\phantom{X }}}
\hbox{\boxtext{\phantom{X }}} } \rm} \ .
\end{equation*}

In view of Theorem (\ref{thmpbw}), $\Prim A$ is
exactly the part of $A_0^{(4)}$ which is orthogonal to the
shuffle products of two commutators $[x_i,x_j]$.

This is just the classical associative setting: In the
representation of $\Sigma_4$ given by
\begin{equation*}
 {\rm
 \vbox{\offinterlineskip \hbox{\boxtext{\phantom{X
}}} \hbox{\boxtext{\phantom{X }}}
 } \rm}\otimes\
 {\rm
 \vbox{\offinterlineskip \hbox{\boxtext{\phantom{X
}}} \hbox{\boxtext{\phantom{X }}}
 } \rm}
=\ \ {\rm
 \vbox{\offinterlineskip \hbox{\boxtext{\phantom{X
}}\boxtext{\phantom{X }}} \hbox{\boxtext{\phantom{X
}}\boxtext{\phantom{X }}}
 } \rm}
 \oplus\ \ {\rm
 \vbox{\offinterlineskip \hbox{\boxtext{\phantom{X
}}\boxtext{\phantom{X }}} \hbox{\boxtext{\phantom{X }}}
\hbox{\boxtext{\phantom{X }}} } \rm}
 \oplus\ \ {\rm
 \vbox{\offinterlineskip \hbox{\boxtext{\phantom{X
}}} \hbox{\boxtext{\phantom{X }}} \hbox{\boxtext{\phantom{X }}}
\hbox{\boxtext{\phantom{X }}} } \rm} \ ,
\end{equation*}
the $\Sigma_4$-module $\ \ {\rm
 \vbox{\offinterlineskip \hbox{\boxtext{\phantom{X
}}\boxtext{\phantom{X }}} \hbox{\boxtext{\phantom{X }}}
\hbox{\boxtext{\phantom{X }}} } \rm}$ with vector space basis
\begin{equation*}
[[x_1,x_2],[x_3,x_4]],\ [[x_1,x_3],[x_2,x_4]],\
 [[x_1,x_4],[x_2,x_3]]
\end{equation*}
  corresponds to the primitive part, and
$\ \ {\rm
 \vbox{\offinterlineskip \hbox{\boxtext{\phantom{X
}}\boxtext{\phantom{X }}} \hbox{\boxtext{\phantom{X
}}\boxtext{\phantom{X }}}
 } \rm}
 \oplus\ \ {\rm
 \vbox{\offinterlineskip \hbox{\boxtext{\phantom{X
}}} \hbox{\boxtext{\phantom{X }}} \hbox{\boxtext{\phantom{X }}}
\hbox{\boxtext{\phantom{X }}} } \rm}\ $ corresponds to the shuffles.

\medskip

Canceling out the part given by shuffle products, we arrive at the
asserted sum of irreducible representations given by
$\Prim\Mag(4)$.

\end{proof}

\begin{proposition}

Let $\sum_{\bf i}a_{\bf i}\ x_{i_1}.x_{i_2}\ldots x_{i_n}, a_{\bf
i}\in K,$ be a homogeneous primitive element of degree $n$ in
$K\langle x_1,x_2,\ldots,x_n\rangle$. Then

\begin{equation*}
\sum_{\bf i}a_{\bf i}\ (\ldots ((x_{i_1}\cdot x_{i_2})\cdot
x_{i_3})\cdot\ldots )\cdot x_{i_n} \text{ and }
\sum_{\bf i}a_{\bf
i}\ x_{i_1}\cdot (x_{i_2}\cdot (\ldots\cdot(x_{i_{n-1}}\cdot
x_{i_n})\ldots))
\end{equation*}
are homogeneous primitive elements of degree $n$ in
$K\{x_1,x_2,\ldots,x_n\}.$

 Especially, for $n\geq 3$,
there are always (at least) two copies of the representation
$\Lie(n)$ present in
 $\Prim\Mag(n)$ .

\end{proposition}

\begin{proof}
When we apply $\Delta_a$ to the terms of $f=\sum_{\bf i}a_{\bf i}\
(\ldots ((x_{i_1}\cdot x_{i_2})\cdot x_{i_3})\cdot\ldots )\cdot
x_{i_n}$, which are left-normed bracketed words, it is immediately
seen that the result is given by tensor products of left-normed
bracketed words. Therefore we do not need to pay attention to
brackets when we want to check if $f$ is primitive.

Completely similar is the case of right-normed bracketed words.

The actions of the general linear and the symmetric group do also
keep left-normed (right-normed, respectively) brackets in order.
Thus the assertion follows.

\end{proof}

\begin{remark}

With a similar argument, one can show that $\Prim\Cmg$ can be
imbedded into $\Prim\Mag$. More exactly, we can map primitive
elements of the free $\Cmg$-algebra $F_{\Cmg}$ into the
commutative subalgebra of $F_{\Mag}$ generated by the variables
with respect to the multiplication $a\circ b:=\frac{1}{2}(a\cdot b
+ b\cdot a)$.

\end{remark}

\begin{remark}

We consider
 $\lambda=(\lambda_1,\lambda_2,\ldots,\lambda_m),
\lambda_1\geq\ldots\geq\lambda_m$, a partition of $n$ (and the
corresponding Young diagram with $\lambda_1$ boxes in the first
row, \ldots, $\lambda_m$ boxes in the last row).

To obtain a basis of the irreducible  $GL_n$-modules (or
$\Sigma_n$-modules) which are contained in $K\{X\}^{(n)}$ or in
$\Prim K\{X\}^{(n)}$,
 one can consider highest weight vectors
(compare \cite{liw}, \cite{lidre}, \cite{lidrecombi} \S 1.2):

\end{remark}

\begin{defn}
A polynomial $f\in K\{X\}^{(n)}$, of multi-degree
$(\lambda_1,\lambda_2,\ldots,\lambda_m)$ is said to be of highest
weight, if
\begin{equation*}
\partial_{ij}\bigl(f\bigr)=0 \text{ all }1\leq j<i\leq m.
\end{equation*}

\end{defn}

\begin{remark}

Each element of highest weight generates an irreducible
$GL_m$-module in $K\{X\}^{(n)}$, and linearly independent elements
of highest weight generate different copies.

Given a highest weight element $f(x_1,x_2,\ldots,x_m)$ of
multi-degree $\lambda$,
 the complete linearization, i.e.\ the multilinear component of
\begin{equation*}
f(\sum_{i=1}^{\lambda_1} x_i,
\sum_{i=\lambda_1+1}^{\lambda_1+\lambda_2}
x_i,\ldots,\sum_{i=n-\lambda_m}^{n} x_i),
\end{equation*}
is used to generate the corresponding $\Sigma_n$-module.

We consider the representation given by $\Prim\Mag(4)$, see
Proposition (\ref{propmagdegfour}),
 and we may look for a basis of primitive highest weight elements.

\end{remark}

\begin{example}

For the five copies of the trivial representation occurring in the
 $GL$-module
$K\{X\}^{(n)}$ we can take all five monomials of degree 5 in
$K\{x_1\}^{(5)}\subset K\{X\}^{(5)}$ as linearly independent
highest weight vectors.

There are only three linearly independent primitive highest weight
vectors, namely
\begin{equation*}
\begin{split}
&\Phi\bigl( (x_1\cdot x_1)\cdot (x_1\cdot x_1)\bigr) =
p(x_1,x_1,x_1,x_1)\\
 &\Phi\bigl( x_1\cdot((x_1\cdot x_1)\cdot x_1)\bigr) =
x_1\cdot((x_1\cdot x_1)\cdot x_1)-3(x_1\cdot(x_1\cdot x_1))\cdot
x_1+2((x_1\cdot x_1)\cdot x_1)\cdot x_1
\\
 &\Phi\bigl( x_1\cdot(x_1\cdot (x_1\cdot x_1))\bigr) =
x_1\cdot(x_1\cdot (x_1\cdot x_1))-4(x_1\cdot(x_1\cdot x_1))\cdot
x_1+3((x_1\cdot x_1)\cdot x_1)\cdot x_1
\\
\end{split}
\end{equation*}
(see Remark (\ref{rembasisonevar})).

\medskip

To compute the primitive highest weight vectors associated to the
Young diagram
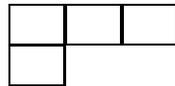
\begin{equation*}
{\rm
 \vbox{\offinterlineskip \hbox{\boxtext{\phantom{X
}}\boxtext{\phantom{X }}\boxtext{\phantom{X }}}
\hbox{\boxtext{\phantom{X }}}
 } \rm}\
\end{equation*}
we have to compute with a basis consisting of 20 monomials, namely
words $x_2 x_1 x_1 x_1$, $x_1 x_2 x_1 x_1$, $x_1 x_1 x_2 x_1$,
$x_2 x_1 x_1 x_1$ with all 5 possible bracketings.

Since the operator $\partial_{21}$ replaces $x_2$ by $x_1$ in
these parenthesized words (and does not change the parentheses),
the coefficient of  $x_2 x_1 x_1 x_1$ is determined by the
coefficients of  $x_2 x_1 x_1 x_1$, $x_1 x_2 x_1 x_1$, $x_1 x_1
x_2 x_1$ (for every type of bracketing).

Under this condition, the requirement to be a constant (i.e.\ a
primitive, in this case) yields five further equations for the
coefficients. A basis for the space of solutions is given by the
following ten primitive highest weight vectors associated to the
Young diagram ${\rm
 \vbox{\offinterlineskip \hbox{\boxtext{\phantom{X
}}\boxtext{\phantom{X }}\boxtext{\phantom{X }}}
\hbox{\boxtext{\phantom{X }}}
 } \rm}\ :$
\begin{equation*}
\begin{split}
f_1=& x_2\cdot(x_1\cdot (x_1\cdot x_1)) - 3x_1\cdot(x_2\cdot
(x_1\cdot x_1))+3x_1\cdot(x_1\cdot (x_2\cdot
x_1))\\
&-x_1\cdot(x_1\cdot (x_1\cdot x_2))
\\
\\
f_2=& (x_1\cdot x_1)\cdot (x_2\cdot x_1) -x_1\cdot(x_1\cdot
(x_2\cdot x_1)) +x_1\cdot(x_1\cdot (x_1\cdot x_2))\\
&- (x_1\cdot x_1)\cdot (x_1\cdot x_2)
\\
\\
f_3=& x_2\cdot(x_1\cdot (x_1\cdot x_1)) -x_2\cdot((x_1\cdot
x_1)\cdot x_1) +x_1\cdot(x_2\cdot (x_1\cdot x_1))\\
&-(x_1\cdot x_2)\cdot(x_1\cdot x_1) +2x_1\cdot((x_1\cdot x_2)\cdot
x_1) -2x_1\cdot(x_1\cdot (x_2\cdot x_1))\\
& - x_1\cdot((x_1\cdot x_1)\cdot x_2) + (x_1\cdot x_1)\cdot
(x_1\cdot x_2)
\\
\\
f_4=& x_2\cdot(x_1\cdot (x_1\cdot x_1)) -x_2\cdot((x_1\cdot
x_1)\cdot x_1) -((x_1\cdot x_1)\cdot x_2)\cdot x_1 \\
&+(x_1\cdot (x_1\cdot x_2))\cdot x_1 -x_1\cdot(x_1\cdot (x_1\cdot
x_2)) +x_1\cdot((x_1\cdot x_1)\cdot x_2)) \\
&-(x_1\cdot (x_1\cdot x_1))\cdot x_2 +((x_1\cdot x_1)\cdot
x_1)\cdot x_2
\\
\\
f_5=& 4 x_2\cdot(x_1\cdot (x_1\cdot x_1)) -3x_2\cdot((x_1\cdot
x_1)\cdot x_1) -3(x_1\cdot x_2)\cdot (x_1\cdot x_1)\\ &
+3((x_1\cdot x_1)\cdot x_2)\cdot x_1  -4 x_1\cdot (x_1\cdot
(x_1\cdot x_2)) +3x_1\cdot ((x_1\cdot x_1)\cdot x_2)\\ &
+3(x_1\cdot x_1)\cdot (x_1\cdot x_2)
 -3((x_1\cdot x_1)\cdot x_1)\cdot x_2\\
\\
f_6=&4 x_2\cdot(x_1\cdot (x_1\cdot x_1)) -3x_2\cdot((x_1\cdot
x_1)\cdot x_1) -(x_2\cdot x_1)\cdot (x_1\cdot x_1)\\ & -2(x_1\cdot
x_2)\cdot (x_1\cdot x_1) +2((x_1\cdot x_2)\cdot x_1)\cdot x_1
-4x_1\cdot(x_1\cdot (x_1\cdot x_2))\\ & +3 x_1\cdot((x_1\cdot
x_1)\cdot x_2)) +3(x_1\cdot x_1)\cdot (x_1\cdot x_2) -2((x_1\cdot
x_1)\cdot x_1)\cdot x_2\\
\\
f_7=& x_2\cdot(x_1\cdot (x_1\cdot x_1)) -x_2\cdot((x_1\cdot
x_1)\cdot x_1) -x_1\cdot(x_2\cdot (x_1\cdot x_1))\\ & +(x_1\cdot
x_2)\cdot (x_1\cdot x_1) +2(x_1\cdot(x_2\cdot x_1))\cdot x_1
-2((x_1\cdot x_2)\cdot x_1)\cdot x_1\\ & + x_1\cdot((x_1\cdot
x_1)\cdot x_2) -(x_1\cdot x_1)\cdot (x_1\cdot x_2) -2
(x_1\cdot(x_1\cdot x_1))\cdot x_2\\ &+2((x_1\cdot x_1)\cdot
x_1)\cdot x_2\\
\\
f_8=& 2x_2\cdot(x_1\cdot (x_1\cdot x_1)) -x_2\cdot((x_1\cdot
x_1)\cdot x_1) -(x_2\cdot x_1)\cdot(x_1\cdot x_1)\\
&-2x_1\cdot(x_2\cdot (x_1\cdot x_1)) +2x_1\cdot((x_2\cdot
x_1)\cdot x_1)
 - x_1\cdot((x_1\cdot x_1)\cdot x_2)\\
 & + (x_1\cdot x_1)\cdot
(x_1\cdot x_2)\\
\end{split}
\end{equation*}
\begin{equation*}
\begin{split}
f_9=& 4x_2\cdot(x_1\cdot (x_1\cdot x_1))- 3x_2\cdot((x_1\cdot
x_1)\cdot x_1) -3(x_2\cdot x_1)\cdot(x_1\cdot x_1)\\ &+2((x_2\cdot
x_1)\cdot x_1)\cdot x_1 -4x_1\cdot(x_1\cdot (x_1\cdot x_2))
+3x_1\cdot((x_1\cdot x_1)\cdot x_2)  \\ &+3(x_1\cdot
x_1)\cdot(x_1\cdot x_2) -2((x_1\cdot x_1)\cdot x_1)\cdot x_2\\
\\
 f_{10}=&
2x_2\cdot(x_1\cdot (x_1\cdot x_1)) -2x_2\cdot((x_1\cdot x_1)\cdot
x_1) -(x_2\cdot x_1)\cdot(x_1\cdot x_1)\\ & +(x_2\cdot(x_1\cdot
x_1))\cdot x_1 -2x_1\cdot(x_1\cdot (x_1\cdot x_2))
+2x_1\cdot((x_1\cdot x_1)\cdot x_2)  \\ &+(x_1\cdot
x_1)\cdot(x_1\cdot x_2) -(x_1\cdot(x_1\cdot x_1))\cdot x_2.
\\
\end{split}
\end{equation*}
In the classical (associative) case, there is only one copy of
this Young diagram present in ${\rm
 \vbox{\offinterlineskip \hbox{\boxtext{\phantom{X
}}\boxtext{\phantom{X }}\boxtext{\phantom{X }}}
\hbox{\boxtext{\phantom{X }}}
 } \rm}
\ \oplus \ {\rm
 \vbox{\offinterlineskip \hbox{\boxtext{\phantom{X
}}\boxtext{\phantom{X }}} \hbox{\boxtext{\phantom{X }}}
\hbox{\boxtext{\phantom{X }}} } \rm}
=
\Lie(4)=\Prim\As(4)$, with highest weight vector
\begin{equation*}
\begin{split}
&x_2.x_1.x_1.x_1 - 3x_1.x_2.x_1.x_1
+3x_1.x_1.x_2.x_1-x_1.x_1.x_1.x_2
\\
 &= [[[x_2,x_1],x_1],x_1].\\
\end{split}
\end{equation*}
Under the projection onto the free $\Cmg$-algebra, the ten
elements $f_i$ yield scalar multiples of only one highest weight
vector (in commuting variables) given by
\begin{equation*}
x_2\cdot(x_1\cdot (x_1\cdot x_1)) +2x_1\cdot(x_1\cdot (x_1\cdot
x_2)) - 3x_1\cdot(x_2\cdot (x_1\cdot x_1)),
\end{equation*}
which is $p(x_2,x_1,x_1,x_1)-p(x_1,x_1,x_1,x_2)$. The image of
$\Phi\bigl( (x_1\cdot x_1)\cdot (x_1\cdot x_1)\bigr) =
p(x_1,x_1,x_1,x_1)$ yields a nontrivial highest weight vector in
commuting variables, too.

Similarly it can be checked that the image of
$p(x_2,x_2,x_1,x_1)-p(x_1,x_2,x_2,x_1)$, see Example
(\ref{exopdegfour}), is a highest weight vector (in commuting
variables)
\begin{equation*}
\begin{split}
&(x_1\cdot x_1)\cdot(x_2\cdot x_2) -2x_1\cdot(x_1\cdot(x_2\cdot
x_2))-2 x_2\cdot(x_2\cdot(x_1\cdot x_1))
\\
 &-(x_1\cdot x_2)\cdot(x_1\cdot x_2)+2x_1\cdot(x_2\cdot(x_1\cdot x_2))
 +2x_2\cdot(x_1\cdot(x_1\cdot x_2))
\\
\end{split}
\end{equation*}
for $\ \ {\rm
 \vbox{\offinterlineskip \hbox{\boxtext{\phantom{X
}}\boxtext{\phantom{X }}} \hbox{\boxtext{\phantom{X
}}\boxtext{\phantom{X }}}
 } \rm}\ $.

The complete linearizations of these three elements generate the
$\Sigma_4$-module

$\Prim\Cmg(4)$, compare Remark (\ref{remprimcmg}).

\end{example}

\bigskip

\begin{remark}

A description of $\Prim\Mag$ by generators and relations would be
interesting. We may take the generators of Proposition
(\ref{propmagdegthree}) together with the operations
$p(x_1,x_2,x_3,x_4)$ and $q(x_1,x_2,x_3,x_4)$ given in Example
(\ref{exopdegfour}) as the first generators.

Since
\begin{equation*}
\begin{split}
([x_1,x_2],x_3,x_4)_b & = p(x_2,x_1,x_3,x_4)-p(x_1,x_2,x_3,x_4),\\
(x_1,x_2,[x_3,x_4])_b & = q(x_1,x_2,x_3,x_4)-q(x_1,x_2,x_4,x_3),\\
\end{split}
\end{equation*}
we may replace $p(x_1,x_2,x_3,x_4)$ and $q(x_1,x_2,x_3,x_4)$ by
\begin{equation*}
\begin{split}
p'(x_1,x_2,x_3,x_4) &:= p(x_1,x_2,x_3,x_4)+p(x_2,x_1,x_3,x_4)
\text{ \ and }
\\
 q'(x_1,x_2,x_3,x_4) &:=
q(x_1,x_2,x_3,x_4)+q(x_1,x_2,x_4,x_3).
\end{split}
\end{equation*}

We get relations, e.g.\ when expressing the 12 elements
$(x_{i_1},[x_{i_2},x_{i_3}],x_{i_4})_b, i_2< i_3,$ by
$p'(x_1,x_2,x_3,x_4)$, $q'(x_1,x_2,x_3,x_4)$, and operations
$[(x_1,x_2,x_3)_b, x_4]$.

\bigskip

We finish with the analogues of Corollary (\ref{corlogcatalan})
and Proposition (\ref{propmagdegfour}) for the case of
$\Prim\Minf$.

\end{remark}

\begin{corollary}
The generating series $f^{\Prim\Minf}(t)$ is given by
\begin{equation*}
\log\Bigl(\frac{5+t-\sqrt{1-6t+t^2}}{4}\Bigr).
\end{equation*}
Thus
\begin{equation*}
\dim\Prim\Minf(n)=(n-1)!C'_n
\end{equation*}
where the integer sequence $C'_n$ is obtained from the
super-Catalan numbers by logarithmic derivation. The sequence
$\dim\Prim\Minf(n), n\geq 1,$ starts with
\begin{equation*}
\begin{split}
0!\cdot &1,\ 1!\cdot 1,\ 2!\cdot 7,\ 3!\cdot 33,\ 4!\cdot 171,\
5!\cdot 901,\ 6!\cdot 4831,\ \ \ldots\\ \text{ i.e.\ }&\ 1,\ \ 1,\
\ \ 14,\ \ \ \ 198,\ \ \ \ 4104,\ \ \
 108120,\ \ 3478320,\ \ \ldots
\end{split}
\end{equation*}
\end{corollary}
\begin{proof}
The series $f^{\Minf}(t)$  is the generating series for the
super-Catalan numbers. Analogously to the proof of Corollary
(\ref{corlogcatalan}) the assertion follows from Theorem
(\ref{thmpbw})(ii).

\end{proof}

\begin{proposition}\label{propminfdegfour}

In terms of representations of the symmetric groups,
\begin{equation*}
\begin{split}
&\Prim\Minf(4)=
\\
 & 8 {\rm \vbox{\offinterlineskip \hbox{
\boxtext{\phantom{X }}\boxtext{\phantom{X }}\boxtext{\phantom{X
}}\boxtext{\phantom{X }}}
 } \rm}\ \oplus\ 25\
{\rm
 \vbox{\offinterlineskip \hbox{\boxtext{\phantom{X
}}\boxtext{\phantom{X }}\boxtext{\phantom{X }}}
\hbox{\boxtext{\phantom{X }}}
 } \rm}
\oplus\ 16\ {\rm
 \vbox{\offinterlineskip \hbox{\boxtext{\phantom{X
}}\boxtext{\phantom{X }}} \hbox{\boxtext{\phantom{X
}}\boxtext{\phantom{X }}}
 } \rm}
 \oplus\ 25\ {\rm
 \vbox{\offinterlineskip \hbox{\boxtext{\phantom{X
}}\boxtext{\phantom{X }}} \hbox{\boxtext{\phantom{X }}}
\hbox{\boxtext{\phantom{X }}} } \rm}
 \oplus\ 8\ {\rm
 \vbox{\offinterlineskip \hbox{\boxtext{\phantom{X
}}} \hbox{\boxtext{\phantom{X }}} \hbox{\boxtext{\phantom{X }}}
\hbox{\boxtext{\phantom{X }}} } \rm} \ .\\
\end{split}
\end{equation*}
\bigskip
\end{proposition}

\begin{proof}

The proof is similar to the proof of Proposition
(\ref{propmagdegfour}). For $A=K\{X\}_{\infty}$,
\begin{equation*}
A_0^{(3)}\otimes{\rm
 \vbox{\offinterlineskip \hbox{\boxtext{\phantom{X
}}}
 } \rm}\ \text{ is equal to }
\end{equation*}
\begin{equation*}
 2 {\rm \vbox{\offinterlineskip
\hbox{ \boxtext{\phantom{X }}\boxtext{\phantom{X
}}\boxtext{\phantom{X }}\boxtext{\phantom{X }}}
 } \rm}\ \oplus\ 7\
{\rm
 \vbox{\offinterlineskip \hbox{\boxtext{\phantom{X
}}\boxtext{\phantom{X }}\boxtext{\phantom{X }}}
\hbox{\boxtext{\phantom{X }}}
 } \rm}
\oplus\ 5\ {\rm
 \vbox{\offinterlineskip \hbox{\boxtext{\phantom{X
}}\boxtext{\phantom{X }}} \hbox{\boxtext{\phantom{X
}}\boxtext{\phantom{X }}}
 } \rm}
 \oplus\ 7\ {\rm
 \vbox{\offinterlineskip \hbox{\boxtext{\phantom{X
}}\boxtext{\phantom{X }}} \hbox{\boxtext{\phantom{X }}}
\hbox{\boxtext{\phantom{X }}} } \rm}
 \oplus\ 2\ {\rm
 \vbox{\offinterlineskip \hbox{\boxtext{\phantom{X
}}} \hbox{\boxtext{\phantom{X }}} \hbox{\boxtext{\phantom{X }}}
\hbox{\boxtext{\phantom{X }}} } \rm} \ .
\end{equation*}

\medskip

We have to compare
\begin{equation*}
\bigoplus_{j=0}^3 A_0^{(j)}\otimes \underbrace{\rm
 \vbox{\offinterlineskip \hbox{\boxtext{\phantom{X
}}\boxtext{$\cdots$}\boxtext{\phantom{X }}}
 } \rm}_{4-j}
\end{equation*}
with $C_4=11$ copies of the regular representation.

It follows that $A_0^{(4)}$ is given by
\begin{equation*}
 8 {\rm \vbox{\offinterlineskip \hbox{ \boxtext{\phantom{X
}}\boxtext{\phantom{X }}\boxtext{\phantom{X }}\boxtext{\phantom{X
}}}
 } \rm}\ \oplus\ 25\
{\rm
 \vbox{\offinterlineskip \hbox{\boxtext{\phantom{X
}}\boxtext{\phantom{X }}\boxtext{\phantom{X }}}
\hbox{\boxtext{\phantom{X }}}
 } \rm}
\oplus\ 17\ {\rm
 \vbox{\offinterlineskip \hbox{\boxtext{\phantom{X
}}\boxtext{\phantom{X }}} \hbox{\boxtext{\phantom{X
}}\boxtext{\phantom{X }}}
 } \rm}
 \oplus\ 25\ {\rm
 \vbox{\offinterlineskip \hbox{\boxtext{\phantom{X
}}\boxtext{\phantom{X }}} \hbox{\boxtext{\phantom{X }}}
\hbox{\boxtext{\phantom{X }}} } \rm}
 \oplus\ 9\ {\rm
 \vbox{\offinterlineskip \hbox{\boxtext{\phantom{X
}}} \hbox{\boxtext{\phantom{X }}} \hbox{\boxtext{\phantom{X }}}
\hbox{\boxtext{\phantom{X }}} } \rm} \ ,
\end{equation*}
for $A=K\{X\}_{\infty}$.

Next we cancel the shuffle products (of commutators) given in the
proof of Proposition (\ref{propmagdegfour}), and the assertion
follows.

\end{proof}

\end{section}

\bigskip
\bigskip
\goodbreak

\begin{section}{Lazard-Lie theory for $\mathcal P$-Hopf algebras}

There is a Lazard-Lie theory
 for complete cogroups over operads, see
 Ginzburg and Kapranov (\cite{ligk}) and Fresse
(\cite{lifra}).
 Given a binary quadratic operad $\mathcal P$,
complete cogroups over $\mathcal P$ are classified by algebras
over the Koszul-dual operad $\mathcal P^{!}$.
 Originally, Lazard \cite{lilaloi} studied
$m$-dimensional formal group laws (for $m\geq 1$)
defined over analyzers $\mathcal A$ by a cohomology theory  $(H^n({\mathcal
A}), \delta^n)_{n\geq 1}$.

It is also possible to define a generalized Lazard-Lie theory for
complete Hopf algebras instead of complete cogroups, see
\cite{lipseudo}.

We are going to sketch the Lazard-Lie theory for $\mathcal P$-Hopf
algebras.
A complex and cohomology groups are defined using the coalgebra structure
of co-addition, in a way which is dual
to the Hochschild $b'$-complex for (not necessarily
unitary) $\As$-algebras (see \cite{lilocyc}).

 The elements of the cohomology group $H^1$, called
pseudo-linear elements by Lazard, are $m$-tuples of primitive
elements in the free $\mathcal P$-algebra on $m$ variables.

\begin{defn}

Let $X=\{x_1,x_2,\ldots,x_m\}$, and let $V=V_X$.
 Let $\mathcal P$ be an operad
equipped with a coherent unit action. Let $n\geq 1$.

We denote by
$\widehat F_{\mathcal P}^{\otimes n}(x_1,\ldots,x_m)$
the complete ${\mathcal P}$-algebra generated by
the $m\cdot n$ elements
\begin{equation*}
x_i^{\otimes n,p}\ \  (\text{for }1\leq p\leq n,\ 1\leq i\leq m)
\end{equation*}
given in Remark (\ref{remtensorpower}).

We consider $\widehat F_{\mathcal P}^{\otimes n}(x_1,\ldots,x_m)$ to
be multigraded by $(\alpha_1,\alpha_2,\ldots,\alpha_n)$ and
denote the space of homogeneous elements of total degree $r=
  \alpha_1+\alpha_2+\ldots+\alpha_n$ by
$F_{\mathcal P}^{\otimes n}(x_1,\ldots,x_m)^{(r)}$, i.e.\
\begin{equation*}
\widehat F_{\mathcal P}^{\otimes n}(x_1,\ldots,x_m)
    =\prod_{r=1}^{\infty} F_{\mathcal P}^{\otimes n}(x_1,\ldots,x_m)^{(r)}.
\end{equation*}
On homogeneous elements, $\alpha_p$ is the degree
 given by
\begin{equation*}
\alpha_p(x_i^{\otimes n,q})=\begin{cases} 1 &: p=q \\ 0 &: \text{else.}\\
  \end{cases}
\end{equation*}

The elements of $\widehat F_{\mathcal P}^{\otimes
n}(x_1,\ldots,x_m)$ are viewed as (generalized) power series over
the "variables" $x_i^{\otimes n,p}$.
    (For $n=1$, elements of $\widehat F_{\mathcal P}^{\otimes 1}(x_1,\ldots,x_m)$ are
just ${\mathcal P}$-power
series in variables $x_1,\ldots,x_m$.)

\end{defn}

\begin{example} \label{exgrplawdata}

Let $A=K 1 \oplus F_{\mathcal P}(x_1,\ldots,x_m)$ or
$A=K 1 \oplus \widehat F_{\mathcal P}(x_1,\ldots,x_m)$ be a free $\mathcal P$-algebra
or a free complete $\mathcal P$-algebra.

Then any morphism $\Delta$ from $A$ to $A\otimes A$ or
$A\hat\otimes A$ is given by an $m$-tuple of elements
$\Delta(x_i)$, $i=1,\ldots,m$, of $\widehat F_{\mathcal
P}^{\otimes 2}(x_1,\ldots,x_m)$.

If, for all $i$, $\Delta(x_i)=x_i^{\otimes 2,1}+x_i^{\otimes
2,2}+$ (terms of higher order), then $A$ together with $\Delta$ is
an augmented $\mathcal P$-bialgebra (or a complete $\mathcal
P$-Hopf algebra in the complete case) if and only if the $m$
elements of $\widehat F_{\mathcal P}^{\otimes 3}(x_1,\ldots,x_m)$
given by
\begin{equation*}
\gamma_i:=(\Delta\hat\otimes\id)\circ\Delta(x_i)-(\Delta\hat\otimes\id)\circ\Delta(x_i)
\end{equation*}
(for $i=1,\ldots,m$) are zero.

\end{example}

\begin{remark}

The tensor product $\otimes $ can be replaced by the
categorical coproduct $\sqcup_{\mathcal P}$ of ${\mathcal P}$-algebras
in the definition above:
One works with (tuples of)
${\mathcal P}$-power series in $m\cdot n$ free variables $x_i^{\sqcup n,p}$.
This is the setting of Lazard's work (c.f.\ \cite{lilaloi})
on $\mathcal P$-formal group laws
(as they are called in \cite{lifra}, \cite{ligk}).

In the special case where $\mathcal P=\Com$, the categorical
  coproduct $\sqcup_{\mathcal P}$ is $\otimes$, and both concepts
  coincide.

\end{remark}

\begin{defn}

We define $\widehat{\mathcal A}_{\mathcal P}^{\otimes n}(x_1,\ldots,x_m)$
to be the vector space of $m$-tuples $f$
with entries $f_j$ in $\widehat F_{\mathcal P}^{\otimes n}(x_1,\ldots,x_m)$
that have the following property:

For all $1\leq j\leq m$, $1\leq p\leq n$, every multihomogeneous
component $g$ of $f_j$ of total degree $\geq 2$ fulfills
$\alpha_p(g)\geq 1$. (In other words, if all occurrences of
$x_1^{\otimes n,p},\ldots,x_m^{\otimes n,p}$ in $g$ are
substituted by 0, the result is 0).

The multigrading is also considered
for $\widehat{\mathcal A}_{\mathcal P}^{\otimes n}(x_1,\ldots,x_m)$.
We denote the homogeneous part of total degree $r$
by ${\mathcal A}_{\mathcal P}^{\otimes n}(x_1,\ldots,x_m)^{(r)}.$

\end{defn}

 \begin{example}
 Let $A=K 1 \oplus \widehat F_{\mathcal P}(x_1,\ldots,x_m)$ be
equipped with a morphism $\Delta: A\to A\hat\otimes A$,
 such that $\Delta(x_i)=x_i^{\otimes 2,1}+x_i^{\otimes 2,2}+$ (terms of higher
order), and let $\gamma_i$ (for $i=1,\ldots,m$) be defined as in
Example (\ref{exgrplawdata}).

Then it can be checked that $\gamma:=(\gamma_1,\ldots,\gamma_r)$
is an element of ${\mathcal A}_{\mathcal P}^{\sqcup 3}(x_1,\ldots,x_m)$.

If $\gamma=0$, the tuple $f=(\Delta(x_1),\ldots,\Delta(x_m))\in
{\mathcal A}_{\mathcal P}^{\sqcup 2}(x_1,\ldots,x_m)$ is called
a quantum $\mathcal P$-formal group law (compare \cite{lih}, \cite{liofqgl}).

Morphisms (coordinate transformations)
are given by $m$-tuples of elements from
 $\widehat F_{\mathcal P}(x_1,\ldots,x_m)$,
 i.e.\ they are given by
 elements of ${\mathcal A}_{\mathcal P}^{\sqcup 1}(x_1,\ldots,x_m)$.
If the coordinate transformation is of the form $x_i\mapsto x_i +$ terms of higher
order (all $i$), one speaks of a strict isomorphism.

\end{example}

\begin{defn}

Let $\Delta_a$ be the co-addition
on
$K\oplus\widehat F_{\mathcal P}(x_1,\ldots,x_m)$
given
by $x_i\mapsto x_i\otimes 1 + 1\otimes x_i, i=1,\ldots,m$.

Let $\Delta'_a$ be given by $\Delta'_a(f)=\Delta_a(f)-f\otimes
1-1\otimes f$.

For $n,i\in \N^*$ with $i<n$, we
 define $\partial_i:
 \widehat F_{\mathcal P}^{\otimes n}(x_1,\ldots,x_m)
\to\widehat F_{\mathcal P}^{\otimes (n+1)}(x_1,\ldots,x_m)$
by
\begin{equation*}
\partial_i(f):=
\underbrace{\id\otimes\ldots\otimes\id}_{i-1}\otimes\Delta'_a
\underbrace{\id\otimes\ldots\otimes\id}_{n-i}.
\end{equation*}

Defined component-wise on $m$-tuples, one sets
\begin{equation*}
\delta=\delta_n: {\mathcal A}_{\mathcal P}^{\otimes n}\to
{\mathcal A}_{\mathcal P}^{\otimes (n+1)},
\text{  }  \delta_n:=\sum_{i=1}^n (-1)^i\partial_i.
\end{equation*}

\end{defn}
\begin{remark}
One checks that $\delta$ restricted to
${\mathcal A}_{\mathcal P}^{\otimes n}(x_1,\ldots,x_m)^{(r)}$
takes in fact values
 in ${\mathcal A}_{\mathcal P}^{\otimes (n+1)}(x_1,\ldots,x_m)^{(r)}$
and that $\delta_n\circ\delta_{n-1}=0$ (all $n$).
In fact, the construction is dual
to the Hochschild $b'$-complex for $\As$-algebras (see \cite{lilocyc}).

\end{remark}

\begin{defn}
Let $\delta_{n,r}:=\delta_n\vert
{\mathcal A}_{\mathcal P}^{\otimes n}(x_1,\ldots,x_m)^{(r)}$.
Let
\begin{equation*}
H^n_r=\ker (\delta_{n,r}) / \im (\delta_{n-1,r}).
\end{equation*}
\end{defn}

\begin{remark}

By the definition of ${\mathcal A}_{\mathcal P}^{\otimes n}(x_1,\ldots,x_m)$,
  $H^n_r=0$ for all $n>r$ always.

Moreover
$H^r_r$, $r\in \N^*$ is in bijection
to anti-symmetric elements of
\begin{equation*}
{\mathcal A}_{\mathcal P}^{\otimes r}(x_1,\ldots,x_m)^{(r)}=
{\mathcal A}_{\mathcal Com}^{\otimes r}(x_1,\ldots,x_m)^{(r)}
\end{equation*}
provided char$(K)=0$ (see \cite{lilaloi}, p.\ 356, and
 \cite{lipseudo}, Proposition 5.9).

The torsion theorem (10.1') of \cite{lilaloi} says in the categorical
case (where $\sqcup$ is used) that
$H^n_r=0$ for all $n\neq r$.

We note that $-\delta_1=\partial_1=\Delta'_a$, thus $H^1_r$
is given by tuples of primitive elements (homogeneous of degree $r$),
and is not zero in general, e.g.\ for ${\mathcal P}=\As$ or
${\mathcal P}=\Mag$. Thus there is no "quantum" version of the torsion
theorem.

For (quantum) group law $r$-chunks (i.e.\ truncated versions
of group laws, truncated after the homogeneous terms of degree $r$)
there is the natural question of extendabilty to $(r+1)$-chunks (and
finally to group laws).

Obstructions for the existence of an extension are found in the
groups $H^3_{r+1}$: a 3-chunk extension exists if and only if the
given 2-chunk defines a Lie bracket.

Obstructions for the uniqueness are elements of
$H^2_{r+1}$.
For example, in the case where ${\mathcal P}=\As$,
the isomorphism classes of $m$-dimensional 3-chunks with trivial
Lie bracket are given by all $(f_1,f_2,\ldots,f_m)$ with
$f_n =$
\begin{equation*}
\begin{split}
&x^{\otimes 2,1}_n + x^{\otimes 2,2}_n + \sum_{(j,i,h)\in {\underline m}^3: i>h}
 \Psi_{j ih}^{(n)}
 \bigl([x^{\otimes 2,1}_h,x^{\otimes 2,1}_i]x^{\otimes 2,2}_j
 +x^{\otimes 2,1}_j[x^{\otimes 2,2}_i,x^{\otimes 2,2}_h]\bigr), \Psi_{j ih}^{(n)}
 \in K.\\
\end{split}
\end{equation*}
The homogeneous part of degree 3 is in $H^2_3$.
For a general normal form of such 3-chunks, see \cite{lih}, \cite{liofqgl}.

It would be interesting to compute the corresponding
obstructions and normal forms for ${\mathcal P}=\Mag$.

\end{remark}

\begin{remark}

Instead of ${\mathcal P}=\Mag$ we can also study the operad
$\Cmg$. Primitive elements of $\Mag$ form the non-associative
Hausdorff series $H(x,y)$ given by $\log(\vee^2(\exp(x),
\exp(y)))$ of \cite{lighmag}, see \cite{lidg} for the definition
of $\exp$ and $\log$. As a motivation to study primitive elements
of $\Cmg$, we can look at the series $\exp(\log(1+x)+\log(1+y))-1$
 in commuting (but non-associative)
variables $x,y$. It starts with the terms
\begin{equation*}
x+y+xy -\frac{(x,x,y)_b-(x,y,y)_b}{3}+\ldots
\end{equation*}
There are projectors from $\Prim\Mag$ to both $\Prim\Cmg$ and
$\Prim\As=\Lie$.

\end{remark}

\end{section}
\bigskip

\goodbreak

\end{chapter}
\vfill\eject

\end{document}